\documentclass[11pt]{amsart}

\usepackage{amsmath, amsthm, amsfonts, amssymb}

\usepackage{enumerate}

\usepackage[bookmarks=false]{hyperref} 

\numberwithin{equation}{section}

\usepackage{xcolor}

\usepackage{verbatim}

\title{Properly proximal von Neumann algebras}
\author{Changying Ding}
\author{Srivatsav Kunnawalkam Elayavalli}
\author{Jesse Peterson}

\address{Department of Mathematics, Vanderbilt University, 1326 Stevenson Center, Nashville, TN 37240, USA}
\email{changying.ding@vanderbilt.edu}
\email{srivatsav.kunnawalkam.elayavalli@vanderbilt.edu}
\email{jesse.d.peterson@vanderbilt.edu}

\thanks{J.P. was supported in part by NSF Grant DMS \#1801125 and NSF FRG Grant \#1853989}

\newtheorem{thm}{Theorem}[section]
\newtheorem{prop}[thm]{Proposition}
\newtheorem{cor}[thm]{Corollary}
\newtheorem{lem}[thm]{Lemma}
\theoremstyle{definition}

\newtheorem{defn/lem}[thm]{Definition/Lemma}
\newtheorem{rem}[thm]{Remark}
\newtheorem{examp}[thm]{Example}

\newcommand{\C}{{\mathbb C}}

\newcommand{\K}{{\mathbb K}}

\newcommand{\N}{{\mathbb N}}

\newcommand{\bS}{{\mathbb S}}

\newcommand{\Y}{{\mathbb Y}}

\newcommand{\X}{{\mathbb X}}

\newcommand{\HH}{{\mathcal H}}

\newcommand{\B}{{\mathbb B}}

\newcommand{\cH}{{\mathcal H}}

\newcommand{\cK}{{\mathcal K}}

\newcommand{\cN}{{\mathcal N}}

\newcommand{\cU}{{\mathcal U}}

\newcommand{\Ad}{\operatorname{Ad}}

\newcommand{\Aut}{\operatorname{Aut}}

\newcommand{\ovt}{\, \overline{\otimes}\,}
\newcommand{\oovt}[1]{\, \overline{\otimes}_{#1}\,}

\newcommand{\ds}{{\sharp\kern-.5pt\sharp}}

\providecommand{\abs}[1]{\lvert#1\rvert}

\newcommand{\actson}{{\, \curvearrowright \,}}
\newcommand{\aactson}[1]{{\, \curvearrowright^{#1} \,}}

\newcommand{\RN}[1]{%
  \textup{\uppercase\expandafter{\romannumeral#1}}%
}

\makeatletter
\DeclareRobustCommand\frownotimes{\mathbin{\mathpalette\frown@otimes\relax}}
\newcommand{\frown@otimes}[2]{%
  \vbox{
    \ialign{##\cr
      \hidewidth$\m@th#1{}_\frown$\kern-\scriptspace\hidewidth\cr
      \noalign{\nointerlineskip\kern-1pt}
      $\m@th#1\otimes$\cr
    }%
  }%
}
\makeatother

\begin{document}
\begin{abstract}
We introduce the notion of proper proximality for finite von Neumann algebras, which naturally extends the notion of proper proximality for groups. Apart from the group von Neumann algebras of properly proximal groups, we provide a number of additional examples, including examples in the settings of free products, crossed products, and compact quantum groups. Using this notion, we answer a question of Popa by showing that the group von Neumann algebra of a nonamenable inner amenable group cannot embed into a free group factor. We also introduce a notion of proper proximality for probability measure preserving actions, which gives an invariant for the orbit equivalence relation.  This gives a new approach for establishing strong ergodicity type properties, and we use this in the setting of Gaussian actions to expand on solid ergodicity results first established by Chifan and Ioana, and later generalized by Boutonnet. The techniques developed also allow us to answer a problem left open by Anantharaman-Delaroche in 1995, by showing the equivalence between the Haagerup property and the compact approximation property for II$_1$ factors. \end{abstract}

\maketitle


\section{Introduction}

A fundamental concept in the theory of von Neumann algebras is the notion of property (Gamma) introduced by Murray and von Neumann \cite{MuVNe43}, which, for a II$_1$ factors, asserts the existence of nontrivial asymptotic central sequences. The absence of this property for the free group factor $L\mathbb F_2$ allowed Murray and von Neumann to distinguish it from the amenable II$_1$ factor $R$, giving the first examples of non-isomorphic separable II$_1$ factors. The subsequent development of the concept played a key role in many of the groundbreaking results advancing the theory of II$_1$ factors, e.g., \cite{Mc69I, Mc69II, Co76, CoWe80}. The asymptotic central sequence algebra $M' \cap M^\omega$ corresponding to a free ultrafilter $\omega \in \beta\mathbb N \setminus \mathbb N$ continues to be an important tool in the recent study of von Neumann algebras, e.g., \cite{Po08, Po10, HoIs16, Oz16, BoChIo17, DeVa18, BaMaOz20}.

A strong refutation of property (Gamma) for $L\mathbb F_2$ was found by Ozawa in \cite{Oz04}, where he showed that if $B \subset L\mathbb F_2$ is any nonamenable von Neumann subalgebra, then the asymptotic commutant $B' \cap (L\mathbb F_2)^\omega$ is not diffuse, i.e., it must possess a minimal projection. In particular, the relative commutant $B' \cap L\mathbb F_2$ is also not diffuse, i.e., $L\mathbb F_2$ is solid (see also \cite{Po07C, Pe09}). The key property that Ozawa was able to exploit was the $C^*$-algebraic property established by Akemann and Ostrand in \cite{AkOs75} that the $*$-homomorphism 
\[
C^*_\lambda \mathbb F_2 \otimes C^*_\rho \mathbb F_2 \ni \sum_{i = 1}^n a_i \otimes x_i \mapsto \sum_{i = 1}^n a_i x_i \in \mathbb B(\ell^2 \mathbb F_2),
\] 
although not continuous with respect to the minimal tensor product, becomes continuous once it is composed with the quotient map onto the Calkin algebra $\mathbb B(\ell^2 \mathbb F_2 ) / \mathbb K(\ell^2 \mathbb F_2)$.   

A compactification $\overline \Gamma$ of a group $\Gamma$ is said to be small at infinity if the right action of $\Gamma$ on itself extends to the trivial action on the boundary $\partial \Gamma$. Every small at infinity compactification is the equivariant image of a universal such compactification, which can be identified (as a $\Gamma \times \Gamma$-space) as the spectrum of the $C^*$-algebra
\[
\mathbb S(\Gamma) = \{ f \in \ell^\infty \Gamma \mid f - R_t(f) \in c_0(\Gamma), \ t \in \Gamma \} \subset \ell^\infty \Gamma.
\]

A group $\Gamma$ is biexact if left action of $\Gamma$ on the spectrum of $\mathbb S(\Gamma)$ is topologically amenable in the sense of \cite{AnDe87}. biexact groups form a robust class of groups that admit a number of different equivalent characterizations (see \cite[Chapter 15]{BrOz08} and the references therein) and Ozawa showed that if $\Gamma$ is biexact, then $L\Gamma$ shares the same solidity property as $L\mathbb F_2$. 

Solidity of a nonamenable von Neumann algebra implies, in particular, that it cannot decompose as the tensor product of two II$_1$ factors. Thus, Ozawa's solidity result for $L\mathbb F_2$ generalizes the result of Ge \cite{Ge96}, where he deduced primeness using free probability techniques. A further generalization was obtained by Ozawa and Popa in \cite{OzPo10I, OzPo10II}, where they show that $L\mathbb F_2$ is strongly solid in the sense that the normalizer of any diffuse amenable sublagebra generates another amenable subalgebra, which also generalizes Voiculescu's result that $L\mathbb F_2$ does not have a Cartan subalgebra \cite{Vo96}. A further generalization showing that every ergodic probability measure preserving (p.m.p.) action of $\mathbb F_2$ gives rise to a crossed product with unique (up to conjugacy) Cartan subalgebra was developed in \cite{PoVa14I, PoVa14II}

To establish strong solidity Ozawa and Popa introduced the notion of a weakly compact action $\Gamma \actson (A, \tau)$ of a group on a tracial von Neumann algebra, which requires the existence of a state $\varphi$ on $\mathbb B(L^2(A, \tau))$ such that $\varphi_{|A} = \tau$, and $\varphi$ is both $A$-central and $\Gamma$-invariant for the conjugation action of $\Gamma$ on $\mathbb B(L^2(A, \tau))$. They then prove two separate results. First, using weak amenability for $\mathbb F_2$, they establish that if $A \subset M \subset L\mathbb F_2$ is an inclusion with $A \subset M$ diffuse and regular, then the action of the normalizer $\mathcal N_{M}(A) \actson A$ is weakly compact. Second, they show that if $\mathcal N_{M}(A) \actson A$ is weakly compact, then $M$ must be amenable. Extending results from \cite{Si11}, the latter result was shown by Chifan and Sinclair to still be the case for all group von Neumann algebras $L\Gamma$ associated with a biexact group $\Gamma$ \cite{ChSi13}.

In order to generalize Ozawa's biexactness methods to a larger class of groups Boutonnet, Ioana and the third author introduced in \cite{BoIoPe21} the notion of proper proximality for a group $\Gamma$ by requiring that the $C^*$-algebra $\mathbb S(\Gamma)$ does not have a state that it invariant under the left action of $\Gamma$. Since a topologically amenable action has an invariant state if and only if the group is amenable, one sees that all nonamenable biexact groups are properly proximal. The class of properly proximal groups though includes many other groups, including nonelementary relatively hyperbolic groups, convergence groups, lattices in non-compact semi-simple Lie groups with finite center, groups admitting a proper cocycle into a nonamenable representation, mapping class groups, many CAT(0) groups, and wreath products $\Sigma \wr \Lambda$ where $\Sigma$ is nontrivial and $\Lambda$ is nonamenable \cite{BoIoPe21, HoHuLe20, DiKuEl21}. Moreover, the class of properly proximal groups is closed under extensions, and under measure-equivalence \cite{BoIoPe21, IsPeRu19}. Contrariwise, natural examples of non-properly proximal groups are given by groups that are inner amenable, i.e., for which there exists a conjugation invariant state $\varphi: \ell^\infty \Gamma \to \mathbb C$ satisfying $\varphi_{| c_0(\Gamma)} = 0$.

The von Neumann algebras associated with properly proximal groups share some of the rigidity properties of those associated with biexact groups. For example, if $\Gamma$ is properly proximal, then $L\Gamma$ does not have property (Gamma) and $L\Gamma$ does not have a diffuse regular von Neumann subalgebra $A \subset L\Gamma$ such that $\mathcal N_{L\Gamma}(A) \actson A$ is weakly compact \cite{BoIoPe21}. Although, unlike for biexact groups, these rigidity properties do not pass to nonamenable von Neumann subalgebras.

It was shown in \cite{IsPeRu19} that the class of properly proximal groups is also stable under $W^*$-equivalence, i.e., if $\Gamma$ and $\Lambda$ satisfy $L\Gamma \cong L\Lambda$, then $\Gamma$ is properly proximal if and only if $\Lambda$ is. This suggests that the notion of proper proximality should be able to be formulated in a purely von Neumann algebraic context. However, the ad hoc argument given in \cite{IsPeRu19} does not give direct evidence on how to achieve this. Based on the analogy between groups and von Neumann algebras (see \cite{Co76II}), a naive approach would be to ask that there does not exist an $M$-central state on the $C^*$-algebra $\{ T \in \mathbb B(L^2M) \mid [T, J x J] \in \mathbb K(L^2M),\ {\rm for\ all\ }x\in M \}$, where $J$ is the conjugation operator. However, Popa showed in \cite{Po87} that the above $C^*$-algebra coincides with $M + \mathbb K(L^2M)$ (cf.\ also \cite{JoPa72}), and it is then easy to see that $M$-central states will always exist. 

Another approach is suggested by Ozawa in \cite{Oza10} where he considers $\K_M$, the space of operators $T \in \B(L^2 M)$ such that both $T$ and $T^*$ are compact when viewed as operators from $M \subset L^2M$ with the uniform norm into $L^2 M$. This forms a hereditary $C^*$-subalgebra that contains $M$ and $JMJ$ in its multiplier algebra $M(\K_M)$, and one may then consider the $C^*$-algebra $S(M)=\{ T \in  M(\K_M) \mid [T, J x J] \in \mathbb \K_M,\ {\rm for\ all\ }x\in M \}$ as an analogue of $\bS(\Gamma)$ in the von Neumann algebra setting. In some cases, this $C^*$-algebra is robust enough to capture the notion of proper proximality. For example, it is shown in \cite[Corollary 4.10, (iii)]{MR17} that for $\mathbb F_2$, if we view the space of continuous function on the Gromov boundary as a $C^*$-subalgebra of $\ell^\infty \mathbb F_2 \subset \B(\ell^2 \mathbb F_2)$, then this is included in $S(L\mathbb F_2)$. As a consequence it follows that $S(L\mathbb F_2)$ does not have an $L\mathbb F_2$-central state. However, for general groups it is unclear why $\mathbb S(\Gamma)$ (or even a suficiently large $C^*$-subalgebra) would embed into $S(L\Gamma)$, since even if we know that the commutator $[T, JxJ]$ is in Ozawa's space of compact operators for $x \in C^*_\lambda(\Gamma)$, there are no general techniques to ensure that this holds also for $x \in L\Gamma$.

We overcome these difficulties by considering instead a strong completion of the compact operators in the sense of Magajna \cite{Ma97, Ma98, Ma00}, which we view as both an operator $M$-bimodule and an operator $JMJ$-bimodule. Using the noncommutative Grothendieck inequality in a similar manner to Ozawa in \cite{Oza10}, we show in Corollary~\ref{cor:M-topology and norms} that this completion coincides with the closure $\mathbb K^{\infty, 1}(M)$ of $\mathbb K(L^2M)$ with respect to the $\|\cdot\|_{\infty, 1}$-norm on $\mathbb B(L^2M)$ given by $\| T \|_{\infty, 1} = \sup_{x \in M, \| x \| \leq 1} \| T \hat{x} \|_1$, where here we consider the natural inclusions $M \subset L^2M \subset L^1M$. (See also \cite{HiPePo22} for some related independent work, where a similar but slightly different closure of the compact operators is being studied.) We may then consider the corresponding ``compactification of $M$'' as
\[
\mathbb S(M) = \{ T \in \mathbb B(L^2M) \mid [T, JxJ] \in \mathbb K^{\infty, 1}(M), \ {\rm for\ all\ }x \in M \}. 
\]

Unlike in the group setting, since $\mathbb K^{\infty, 1}(M)$ is not an ideal in $\mathbb B(L^2M)$, we see that $\mathbb S(M)$ is not a $C^*$-subalgebra of $\B(L^2M)$. However, it is still an operator system that is also an $M$-bimodule, and so we may define the von Neumann algebra $M$ to be properly proximal if there does not exist an $M$-central state $\varphi$ on $\mathbb S(M)$ such that $\varphi_{|M}$ is normal. Owing to the fact that $\mathbb K^{\infty, 1}(M)$ is a strong bimodule in the sense of Magajna, it follows that an operator $T \in \mathbb B(L^2M)$ belongs to $\mathbb S(M)$ if $[T, JxJ] \in \mathbb K^{\infty, 1}(M)$ holds for all $x$ in a set that generates $M$ as a von Neumann algebra. This allows us to pass between $\mathbb S(\Gamma)$ and $\mathbb S(L\Gamma)$ so that we may establish in Theorem~\ref{thm: iff} the equivalence of proper proximality of a group $\Gamma$ and its von Neumann algebra $L\Gamma$. The space $\mathbb S(M)$ also gives a natural von Neumann algebraic setting to study biexactness properties, and this topic will be taken up in a subsequent paper by a subset of the authors \cite{DiPe22}. 

Apart from group von Neumann algebras of properly proximal groups, we also provide several other examples of properly proximal von Neumann algebras in the settings of free products, crossed products, and compact quantum groups. Generalizing results for group von Neumann algebras in \cite{BoIoPe21}, we show that properly proximal von Neumann algebras $M$ never have property (Gamma) (Corollary~\ref{cor:gamma}), nor do they admit diffuse regular subalgebras such that the normalizer action is weakly compact (Theorem~\ref{thm:weaklycompact}). 

Effros established in \cite{Ef75} that if $\Gamma$ is a group such that $L\Gamma$ has property (Gamma), then $\Gamma$ must be inner amenable. Vaes constructed in \cite{Va12} examples of inner amenable groups so that the group von Neumann algebra does not have property (Gamma). Since Ozawa's solidity theorem shows that every nonamenable subalgebra of $L\mathbb F_2$ does not have property (Gamma), Popa asked if it is still true that the group von Neumann algebra of a nonamenable inner amenable group cannot embed into $L\mathbb F_2$ \cite{Po21Q}. Since by Theorem~\ref{thm: iff} the group von Neumann algebra of an inner amenable group is not properly proximal, we obtain a solution to Popa's problem with the following result. 

\begin{thm}[See Theorem~\ref{thm:biexactpropp}]
If $\Gamma$ is a biexact group, then every von Neumann subalgebra of $L \Gamma$ is either properly proximal or has an amenable summand.
\end{thm}

We also introduce proper proximality for trace-preserving group actions $\Gamma \actson (M, \tau)$ by requiring the non-existence of a state $\varphi$ on $\mathbb S(M)$ such that $\varphi_{|M} = \tau$ and $\varphi$ is both $M$-central and $\Gamma$-invariant for the conjugation action of $\Gamma$ on $\mathbb S(M)$. The notion of a properly proximal action is averse to weak compactness. Moreover, it passes to diffuse $\Gamma$-invariant subalgebras, and for actions on separable abelian von Neumann algebras it is an invariant of the orbit equivalence relation (see Proposition~\ref{prop:normalizer}).

As in the setting for groups in \cite{BoIoPe21}, we also produce versions of these results that are relative to certain boundary pieces. (See Section~\ref{sec:relcompact} for the definition of a boundary piece in the setting of von Neumann algebras.) Using a certain normal bidual construction related to the work of Magajna \cite{Ma05}, we are able to show in Proposition~\ref{prop:crossedproduct implies action} and Theorem~\ref{thm:crossedproduct} that if $\Gamma$ is properly proximal, then an action $\Gamma \actson (M, \tau)$ of a countable group on a separable tracial von Neumann algebra is properly proximal if and only if the crossed product $M \rtimes \Gamma$ is properly proximal. Using similar methods, we show in Theorem~\ref{thm:tensorproduct} that a tensor product $M_1 \ovt M_2$ of separable II$_1$ factors is properly proximal if and only if $M_1$ and $M_2$ are properly proximal. 

We also show that proper proximality can also be deduced using Popa's malleable deformations \cite{Po06D, Po06B, Po06C, Po07B, Po08}. Gaussian actions corresponding to orthogonal representations have such malleable deformations as noticed in \cite{Fu07} and studied in \cite{PeSi12}. Using these deformations, we show proper proximality for Gaussian actions associated to orthogonal representations that are nonamenable in the sense of Bekka \cite{Be90}. When the representation is weakly contained in the left regular representation, we also show proper proximality for any nonamenable subequivalence relation. This then gives a generalization of the results of Chifan and Ioana from \cite{ChIo10} and Boutonnet from \cite{Bo12}, establishing, without additional conditions, solid ergodicity (See Section 5 of \cite{Ga10}) for Gaussian actions associated to representations that are weakly contained in the left regular representation.

\begin{thm}\label{thm:soliderg} 
Let $\Gamma$ be a nonamenable group, and $\pi: \Gamma \to \mathcal O(\mathcal H)$ an orthogonal representation on a separable real Hilbert space such that some tensor power of $\pi$ is weakly contained in the left regular representation. Then for every nonamenable subequivalence relation $\mathcal R$ of the equivalence relation associated to the Gaussian action $\Gamma \actson A_{\mathcal H}$, there exists an $\mathcal R$-invariant projection $p \in A_{\mathcal H}$ so that the action $\mathcal N_{L \mathcal R}(pA_{\mathcal H}) \actson pA_{\mathcal H}$ is properly proximal and $\mathcal R_{p^\perp}$ is amenable. In particular, the equivalence relation associated to the Gaussian action $\Gamma \actson A_{\mathcal H}$ is solidly ergodic.
\end{thm}

As a consequence of the tools we develop, we are also able to show the equivalence between the Haagerup property defined in \cite{Ch83}, and the compact approximation property defined in \cite{AnDe95} for finite von Neumann algebras (see also \cite[p.\ 559]{Jo02}). This generalizes Proposition 4.16 in \cite{AnDe95}, where the equivalence was established for group von Neumann algebras by showing that, in this case, both properties are equivalent to the Haagerup property for the group. 

\begin{thm}\label{thm:AD}
Let $M$ be a finite von Neumann algebra with normal faithful trace $\tau$. The following conditions are equivalent:
\begin{enumerate}
\item $M$ has the Haagerup property, i.e., there exists a net $\{ \phi_i \}_i$ of normal completely positive maps from $M$ to $M$ such that 
\begin{enumerate}
\item $\tau \circ \phi_i( x^* x) \leq \tau(x^* x)$ for all $x \in M$;
\item $\lim_i \| \phi_i(x) - x \|_2 = 0$ for all $x \in M$;
\item each $\phi_i$ induces a compact bounded operator on $L^2M$.
\end{enumerate}
\item $M$ has the compact approximation property, i.e., there exists a net $\{ \phi_i \}_i$ of normal completely positive maps from $M$ to $M$ such that 
\begin{enumerate}
\item for all $x \in M$ we have $\lim_i \phi_i(x) = x$ ultraweakly;
\item for all $\xi \in L^2M$ and $i \in I$, the map $x \mapsto \phi_i(x) \xi$ is compact from the normed space $M$ to $L^2M$. 
\end{enumerate}
\item There exists a net $\{ \phi_i \}_i$ of normal completely positive maps from $M$ to $M$ such that
\begin{enumerate}
\item for all $x \in M$ we have $\lim_i \phi_i(x) = x$ ultraweakly;
\item for all $\xi \in L^2M$ and $i \in I$, the map $x \mapsto \phi_i(x) \xi$ is compact from the normed space $M$ to $L^1M$. 
\end{enumerate}
\end{enumerate}
\end{thm}

\subsection*{Acknowledgements} The third author thanks Cyril Houdayer and Sorin Popa for useful comments regarding this project.

\section{Preliminaries}

\subsection{The extended Haagerup tensor product}

If $X \subset \B(\mathcal H)$ is an operator space (for us, all operator spaces will be norm closed) and $I, J$ are sets, then we denote by $\mathbb M_{I, J}(X)$ the set of $I \times J$ matrices having entries in $X$ that represent bounded operators from $\ell^2 J \ovt \mathcal H$ into $\ell^2 I \ovt \mathcal H$. We write $C_I(X) = \mathbb M_{I, 1}(X)$, and $R_J(X) = \mathbb M_{1, J}(X)$, for the column and row matrices respectively.  When no confusion will arise, we will identify each matrix in $\mathbb M_{I, J}(X)$ with the corresponding operator that it represents, so that we view $\mathbb M_{I, J}(X)$ as an operator subspace of $\B(\ell^2 J \ovt \mathcal H, \ell^2 I \ovt \mathcal H)$.
If $a \in \B(\mathcal H, \mathcal K)$, then we let $a^{(J)} \in \B(\ell^2 J \ovt \mathcal H, \ell^2 J \ovt \mathcal K)$ be the diagonal operator $a^{(J)} = 1 \otimes a$.

If $X \subset \B(\mathcal H)$, and $Y \subset \B(\mathcal K)$ are operator spaces, $x \in \mathbb M_{I, J}(X)$, and $y \in \mathbb M_{J, K}(Y)$, then we let $x \odot y$ denote the normal completely bounded operator in $CB^\sigma( \B(\mathcal K, \mathcal H), \mathbb M_{I, K}( \B(\mathcal K, \mathcal H)))$
given by $(x \odot y)(a) = x a^{(J)} y$, for $a \in \B(\mathcal K, \mathcal H)$.

We take the perspective in \cite{Ma97} and define the extended Haagerup tensor product $X \otimes_{eh} Y$ to consist of all maps of the form $x \odot y$ for $x \in R_J(X)$, $y \in C_J(Y)$, for some set $J$. The extended Haagerup tensor product is again an operator space, with the norm on $\mathbb M_n(X \otimes_{eh} Y)$ inherited from $CB^\sigma( \B(\mathcal K, \mathcal H), \mathbb M_n(\B(K, \mathcal H))))$.

If $M \subset \B(\mathcal H)$ and $N \subset \B(\mathcal K)$ are von Neumann algebras, then each element in $M \otimes_{eh} N$ corresponds to an $M'$-$N'$-bimodular map in $CB^\sigma(\B(\mathcal H))$. Moreover, every $M'$-$N'$-bimodular map arises in this way \cite[Theorem 4.2]{BlSm92} so that we have operator space isomorphisms
\[
M \otimes_{eh} N = CB_{M'-N'}^\sigma(\B(\mathcal H)) \cong CB_{M'-N'}(\K(\mathcal H), \B(\mathcal H)).
\]
Note that composition then turns $M \otimes_{eh} N$ into a Banach algebra, which agrees with the usual algebra structure on the algebraic tensor product $M \otimes N^{\rm op}$, i.e., for $a, b \in M$ and $x, y \in N$ we have $( a \odot x ) \circ (b \odot y) = (ab) \odot (yx)$. We also have that $M \otimes_{eh} N$ is naturally a dual Banach space where the weak$^*$-topology on bounded sets is given by ultraweak convergence applied to each compact operator. Multiplication is then weak$^*$-continuous in the second variable, but not in the first variable in general (see the remark before theorem 4.2 in \cite{BlSm92}) due to the fact that the space of compact operators need not be preserved by maps in $M \otimes_{eh} N$.

\subsection{Strong operator bimodules}

If $A$ and $B$ are unital $C^*$-algebras, then an operator $A$-$B$-bimodule consists of unital $*$-homomorphisms $\pi: A \to \B(\mathcal H)$, $\rho: B \to \B(\mathcal H)$, together with operator subspace $X \subset \B(\mathcal H)$ such that $X$ is a $\pi(A)$-$\rho(B)$-bimodule whose bimodule structure is given by composition of operators. We write $a x b$ for the element $\pi(a)x\rho(b)$ whenever $a \in A$, $b \in B$, and $x \in X$. If $A = B$, $\pi = \rho$, and $X$ happens to be on operator system in $\B(\mathcal H)$ that contains $\pi(A)$, then we say that $X$ is an operator $A$-system (see  \cite[p. 215]{Pa02}). We remark that $A$-$B$-bimodules and operator $A$-systems may also be defined abstractly without referring to an explicit Hilbert space $\mathcal H$ (see \cite{ChEfSi87} and \cite[Corollary 15.13]{Pa02} respectively). We say that $X$ is a dual operator $A$-$B$-bimodule (resp.\ dual operator $A$-system) if $X$ may be realized above as an ultraweak closed subspace of $\B(\mathcal H)$. 

If $M$ and $N$ are von Neumann algebras, then we say that an operator $M$-$N$-bimodule (resp.\ $M$-system) is normal if the representations $\pi: M \to \B(\mathcal H)$, and $\rho: N^{\rm op} \to \B(\mathcal H)$ from above may be taken to be normal and faithful. An $M$-$N$ sub-bimodule $X \subset \B(\mathcal H)$ is strong if whenever $a \in R_I(M)$, $x \in \mathbb M_{I, J}(X)$, and $b \in C_J(N)$ for some sets $I$ and $J$, then we have $a x b \in X \subset \B(\mathcal H)$.  

Strong bimodules were introduced by Magajna in \cite{Ma00}, generalizing the notion of a strong module (i.e., a strong $M$-$\mathbb C$ bimodule), which were introduced in \cite{EfRu88} and developed in \cite{Ma97, Ma98}. We summarize here some of Magajna's results from \cite{Ma97, Ma98, Ma00} that we will use involving strong bimodules. 

Strong bimodules are automatically norm closed. Moreover, if $X \subset \B(\mathcal H)$ is strong, then it is strong under any other realization as a concrete normal $M$-$N$ operator bimodule. If $X \subset \B(\mathcal H)$ is an $M$-$N$ bimodule, then denote by $\overline X^{st} $    the smallest strong bimodule containing $X$ as a sub-bimodule. Dual bimodules are always strong, but there are many strong bimodules that are not dual, e.g., in the case when $M = N = \mathbb C$, all closed bimodules are strong.

If $X$ is a normal operator $M$-$N$ bimodule, then the $M$-$N$-topology on $X$ is given by the family of semi-norms
\[
s_\omega^\rho(x) = \inf \{ \omega(a^*a)^{1/2} \| y \| \rho(b^*b)^{1/2} \},
\]
where $\omega$ and $\rho$ are positive normal linear functionals on $M$ and $N$, respectively, and the infimum is taken over all decompositions $x = a^* y b$ where $a \in M$, $b \in N$, and $y \in X$. It is not clear, a priori, that $s_\omega^\rho$ is a seminorm, but this is indeed the case, and follows from an argument similar to the one showing that the Haagerup norm is a norm.

If $X \subset Y$ is an inclusion of normal operator $M$-$N$ bimodules, then the restriction of the $M$-$N$-topology on $Y$ gives the $M$-$N$-topology on $X$. Moreover, if $M$ and $N$ are normally represented on $\B(\mathcal H)$, then an $M$-$N$ sub-bimodule is closed in the $M$-$N$-topology if and only if it is a strong bimodule \cite[Theorem 3.10]{Ma00}. 

If $M_0 \subset M$ and $N_0 \subset N$ are ultraweakly dense $C^*$-subalgebras and $X$ is an $M_0$-$N_0$ operator bimodule, then we let $X^{M_0 \sharp N_0} \subset X^*$ denote the space of bounded linear functionals $\varphi$ such that for each $x \in X$ the map $M_0 \times N_0 \ni (a, b) \mapsto \varphi(a x b)$ extends to a map on $M \times N$ that is separately ultraweakly continuous in $a$ and $b$. 
Note that $X^{M_0 \sharp N_0} $ is a norm closed subspace of $X^*$. 
If $X$ is a normal $M$-$N$ operator bimodule, then $X^{M \sharp N}$ coincides with the linear functionals that are continuous in the $M$-$N$ topology \cite[Theorem 3.7]{Ma00}. The following proposition generalizes this fact to the situation where $X$ can be made a bimodule in multiple ways, e.g., when $M$ is a II$_1$ factor, we may want to view $X = \B(L^2M)$ as both an $M$-$M$ bimodule and a $M'$-$M'$ bimodule. 

If $X$ is a normal $M$-$N$ operator bimodule, then we will often use the notation $X^\sharp$ for $X^{M \sharp N}$ when the bimodule structure is clear from the context.  We remark that Magajna favors the notation $X^\sharp$ for the continuous dual of $X$ in \cite{Ma97, Ma98, Ma00, Ma05}, however no confusion should arise in the present context as we will always denote the continuous dual by $X^*$. 

\begin{prop}\label{prop:quadmodrep}
Let $\mathcal H$ be a Hilbert space and suppose $M, \tilde M, N, \tilde N \subset \B(\mathcal H)$ are von Neumann algebras such that $\tilde M \subset M'$ and $\tilde N \subset N'$. Suppose $X \subset \B(\mathcal H)$ is an operator space that is both an $M$-$N$ bimodule and a $\tilde M$-$\tilde N$ bimodule. Suppose $\varphi \in X^{M \sharp N} \cap X^{\tilde M \sharp \tilde N}$, then there exist Hilbert spaces $\mathcal K_1, \mathcal K_2$, normal and commuting representations of $M$ and $\tilde M$ on $\mathcal K_1$, and of $N$ and $\tilde N$ on $\mathcal K_2$, vectors $\xi \in \mathcal K_1$, $\eta \in \mathcal K_2$, and a completely bounded map $\phi: X \to \B(\mathcal K_2, \mathcal K_1)$ that is both $M$-$N$-bimodular and $\tilde M$-$\tilde N$-bimodular, and such that for all $x \in X$ we have $\varphi(x) = \langle \phi(x) \eta, \xi \rangle$. 
\end{prop}
\begin{proof}
The proof is based on \cite[Theorem 4.2]{Ma98} and \cite[Theorem 3.4]{EfRu88}. Set $\mathcal A = \B(\mathcal H)$ and extend $\varphi$ to a linear functional $\tilde \varphi \in \mathcal A^*$ by Hahn-Banach. Since $\mathcal A$ is a $C^*$-algebra, we may then find a Hilbert space $\mathcal L$, a $*$-representation $\pi: \mathcal A \to \B(\mathcal L)$ and vectors $\xi_0, \eta_0 \in \mathcal L$ such that $\tilde \varphi(a) = \langle \pi(a) \eta_0, \xi_0 \rangle$, for all $a \in \mathcal A$. We may then set $\mathcal L_1$ to be the closure of the subspace spanned by $\pi(X) \eta_0$. Since $\pi$ is a $*$-homomorphism on $M$ and $\tilde M$, it follows that $P_{\mathcal L_1}$ commutes with $\pi(M)$ and $\pi(\tilde M)$, and so we obtain commuting $*$-representations of $M$ and $\tilde M$ on $\mathcal L_1$ by restriction. 

We define $\tilde \pi: X \to \B(\mathcal L, \mathcal L_1)$ by $\tilde \pi(x) = P_{\mathcal L_1} \pi(x)$. This is then a completely bounded $M$-$N$ and $\tilde M$-$\tilde N$-bimodular map. Moreover, if we set $\xi = P_{\mathcal L_1} \xi_0$, then we have 
\[
 \varphi(x) = \langle \pi(x) \eta_0, \xi_0 \rangle = \langle P_{\mathcal L_1} \pi(x) \eta_0, \xi_0 \rangle = \langle \tilde \pi(x) \eta_0, \xi \rangle.
\]

Next, we let $\mathcal L_2$ be the closure of the subspace $\tilde \pi(X^*) \xi$, and note that, as above, $P_{\mathcal L_2}$ commutes with both $N$ and $\tilde N$ and hence we obtain $*$-representations of $N$ and $\tilde N$ on $\mathcal L_2$ by restriction. We define $\phi_0: X \to \B(\mathcal L_2, \mathcal L_1)$ by $\phi_0(x) = \tilde \pi(x) P_{\mathcal L_2}$. This is then a completely bounded $M$-$N$ and $\tilde M$-$\tilde N$-bimodular map, and if we set $\eta = P_{\mathcal L_2} \eta_0$, then as above we have $\varphi(x) = \langle \phi_0(x) \eta, \xi \rangle$. 

Next, we let $\mathcal K_1$ be the closure of the subspace spanned by $\pi(M)\pi(\tilde M) \xi$, and we let $\mathcal K_2$ be the closure of the subspace spanned by $\pi(N) \pi(\tilde N) \eta$. We then have that $\mathcal K_1$ is invariant under $M$ and $\tilde M$, and hence we obtain commuting $*$-representations by restriction. We similarly obtain commuting $*$-representations of $N$ and $\tilde N$ on $\mathcal K_2$. Since $\eta \in \mathcal K_2$ and $\xi \in \mathcal K_1$ if we set 
$\phi(x) = P_{\mathcal K_1} \phi_0(x) P_{\mathcal K_2}$, then $\phi$ is $M$-$N$ and $\tilde M$-$\tilde N$-bimodular and we again have $\varphi(x) = \langle \phi(x) \eta, \xi \rangle$ for all $x \in X$.

We note that since $\varphi \in X^{M \sharp N}$, if $x \in X$ and $\zeta = \sum_{j = 1}^n \pi(a_i) \pi(b_i) \xi$ for $a_1, \ldots, a_n \in M$, and $b_1, \ldots, b_n \in \tilde M$, then the linear functional $M \ni a \mapsto \langle \pi(a) P_{\mathcal K_1} \pi(x) \eta_0, \zeta \rangle = \sum_{j = 1}^n \varphi(a_i^* a b_i^* x)$ is normal. As $P_{\mathcal K_1} \pi(X) \eta_0$ and ${\rm span}\{ \pi(M) \pi(\tilde M) \xi\}$ are both dense in $\mathcal K_1$, it follows that $M$ is normally represented on $\mathcal K_1$. A similar argument shows that $\tilde M$ is also normally represented on $\mathcal K_1$, and that both $N$ and $\tilde N$ are normally represented on $\mathcal K_2$. 
\end{proof}

If we assume $\varphi \in X^{M \sharp N} \cap X^{\tilde M \sharp \tilde N}$ is such that $\| \varphi \| \leq 1$, and we let $\phi$ be as in the previous proposition, then an application of Paulsen's trick (e.g., Lemma 8.1 in \cite{Pa02}) gives a unital completely positive (u.c.p.) $C^*(M, \tilde M) \oplus C^*(N, \tilde N)$-bimodular map $\Phi: \mathcal S_X \to \B(\mathcal K_1 \oplus \mathcal K_2)$, where $S_X$ is the operator $C^*(M, \tilde M) \oplus C^*(N, \tilde N)$-system 
\[
\mathcal S_X = \left\{ \left( \begin{smallmatrix}
a & x \\
y^* & b
\end{smallmatrix} \right) \mid a \in C^*(M, \tilde M), b \in C^*(N, \tilde N), x, y \in X
\right\}.
\]
If $X$ is a subspace of $Y \subset \B(\mathcal H)$ such that $Y$ is also both an $M$-$N$ and $\tilde M$-$\tilde N$ bimodule, then we may use Arveson's extension theorem to extend $\Phi$ to a u.c.p.\ map on $\mathcal S_Y$. Restricting this to the corner $Y$ in $\mathcal S_Y$ then gives the analogue of the Hahn-Banach theorem in our setting (See Proposition 4.5 in \cite{Ma98} for another approach). We record this result here for use later.

\begin{prop}\label{prop:normalhahnbanach}
Let $\mathcal H$ be a Hilbert space and suppose $M, \tilde M, N, \tilde N \subset \B(\mathcal H)$ are von Neumann algebras such that $\tilde M \subset M'$ and $\tilde N \subset N'$. Suppose $Y \subset \B(\mathcal H)$ is an operator space that is both an $M$-$N$ bimodule and a $\tilde M$-$\tilde N$ bimodule, and suppose that $X \subset Y$ is both an $M$-$N$ sub-bimodule and a $\tilde M$-$\tilde N$ sub-bimodule. If $\varphi \in X^{M \sharp N} \cap X^{\tilde M \sharp \tilde N}$, then there exists $\psi \in Y^{M \sharp N} \cap Y^{\tilde M \sharp \tilde N}$ with $\| \psi \| = \| \varphi \|$ such that $\psi_{| X} = \varphi$. 
\end{prop}

The following descriptions of the strong completion of bimodules will be useful for us later.

\begin{prop}\label{prop:strongclosure}
Let $M$ and $N$ be normally represented in $\B(\mathcal H)$ and let $X \subset \B(\mathcal H)$ be a  $M$-$N$ sub-bimodule. Fix $x \in \B(\mathcal H)$. The following are equivalent:
\begin{enumerate}
\item \label{en:strong1} $x$ is in the strong $M$-$N$ bimodule generated by $X$.
\item \label{en:strong2} There exist orthogonal families of projections 
	$\{e_i\}_{i\in I}\subset M$ and $\{f_j\}_{j\in J}\subset N$ 
	with $\sum_{i\in I}e_i=1$ and $\sum_{j\in J}f_j=1$, 
	such that $e_i x f_j\in \overline X$ for any $i\in I$ and $j\in J$.
\item \label{en:strong3} If $\varphi \in \B(\mathcal H)^{M \sharp N}$ is any linear functional such that $\varphi_{| X} = 0$, then we also have $\varphi(x) = 0$. 
\end{enumerate}
\end{prop}
\begin{proof}
It is easy to see that (\ref{en:strong2}) $\implies$ (\ref{en:strong1}), and 
	since the strong $M$-$N$ bimodule generated by $X$ is the closure in the $M$-$N$-topology, 
	we have that (\ref{en:strong1}) and (\ref{en:strong3}) are equivalent by the Hahn-Banach theorem. 

To see (\ref{en:strong1}) $\implies$ (\ref{en:strong2}), 
	note that an $M$-$N$ bimodule is strong if and only if 
	it is a strong left $M$-module and a strong right $N$-module,
	and we may obtain such families $\{e_i\}_{\i\in I}$ and $\{f_j\}_{j\in J}$ by applying \cite[Proposition 2.2]{Ma98} twice.

\end{proof}

We will also need a version of the previous proposition in the setting of Proposition~\ref{prop:quadmodrep}. We omit the proof, which is similar to the proof of Proposition~\ref{prop:strongclosure}.

\begin{prop}\label{prop:quadmodstrong}
Let $\mathcal H$ be a Hilbert space and suppose $M, \tilde M, N, \tilde N \subset \B(\mathcal H)$ are von Neumann algebras such that $\tilde M \subset M'$ and $\tilde N \subset N'$. Suppose $X \subset \B(\mathcal H)$ is an operator space that is both an $M$-$N$ bimodule and a $\tilde M$-$\tilde N$ bimodule. For $x \in \B(\mathcal H)$, the following are equivalent:
\begin{enumerate}
\item $x$ is in the smallest subspace containing $X$ that is both a strong $M$-$N$ bimodule and a strong $\tilde M$-$\tilde N$ bimodule. 
\item There exist orthogonal families of projections 
	$\{e_i\}_{i\in I}\subset M$, $\{\tilde e_i\}_{i\in I}\subset \tilde M$,
	$\{f_j\}_{j\in J}\subset N$ and $\{\tilde f_j\}_{j\in J}\subset \tilde N$
	which all sum to $1$,
	such that $e_i \tilde e_i x \tilde f_j f_j\in \overline X$ for any $i\in I$ and $j\in J$.
\item If $\varphi \in \B(\mathcal H)^{M \sharp N} \cap \B(\mathcal H)^{\tilde M \sharp \tilde N}$ is any linear functional such that $\varphi_{| X} = 0$, then we also have $\varphi(x) = 0$. 
\end{enumerate}
\end{prop}

\subsection{Universal representations of normal $M$-$N$ bimodules}
If $X$ is a normal operator $M$-$N$-bimodule, then by a representation of $X$ we mean an $M$-$N$-bimodular complete contraction $\pi: X \to \B(\mathcal H, \mathcal K)$ where $M$ and $N$ are normally represented in $\B(\mathcal K)$ and $\B(\mathcal H)$ respectively. We let ${\rm Rep}(X)$ denote the space of all representations on some fixed Hilbert space of sufficiently large dimension. Note that direct sums of representations again give representations, and hence we may define the universal representation of $X$ to be $\pi_u = \oplus_{\pi \in {\rm Rep}(X)} \pi$.

By Proposition~\ref{prop:quadmodrep} (or in this case by \cite[Theorem 3.10]{Ma00}), the space $X^{ \sharp }$ coincides with the coefficient linear functionals corresponding to representations.  Although the space $X^{ \sharp }$ need not be an operator $M$-$N$ bimodule, it is still naturally a $N$-$M$ bimodule under the action $(a \cdot \varphi \cdot b)(x) = \varphi(bxa)$. We then have that $(X^{ \sharp })^*$ is a dual $M$-$N$ bimodule. The following proposition is easily deduced. 

\begin{prop}
The universal representation $\pi_u$ induces a weak$^*$ homeomorphic isometric $M$-$N$ bimodular mapping of $( X^{ \sharp })^*$ onto the ultraweak closure $\overline{\pi_u(X)}$.
\end{prop}

We may then give $(X^{ \sharp })^*$ the normal dual operator $M$-$N$-bimodule structure coming from this isomorphism and we note that this satisfies the universal property that if $\pi: X \to Y$ is an $N$-$M$ bimodular complete contraction and $Y$ is a dual normal operator $M$-$N$-bimodule, then there exists a unique normal $M$-$N$-bimodular complete contraction $\tilde \pi: (X^{ \sharp })^* \to Y$ such that $\tilde \pi_{| X} = \pi$. 

If $Y$ is a normal operator $M$-$N$ bimodule containing $X$ as a sub-bimodule, then it follows from Proposition~\ref{prop:normalhahnbanach} (or in this case by \cite[Proposition 4.5]{Ma98}) that the identity map from $X$ to $Y$ extends to an $M$-$N$-bimodular completely isometric normal map from $(X^{ \sharp })^*$ into $(Y^{ \sharp })^*$, and so we will identify $(X^{ \sharp })^*$ as a sub-bimodule of $(Y^{ \sharp })^*$.

An alternative construction for the ``bidual'' is given by Magajna in \cite{Ma05}. There he introduces the bimodule dual of $X$ as $X^\natural = CB_{M-N}(X, \B(L^2 N, L^2M))$. This is an abstract dual normal operator $M'$-$N'$-bimodule and hence the bimodule bidual $X^{\natural \natural}$ gives an abstract dual normal operator $M$-$N$-bimodule. The following result of Magajna shows that these two notions of biduals agree. 

\begin{prop}[Magajna {\cite[Corollary 3.5(iii)]{Ma05}}]
For any normal operator $M$-$N$-bimodule $X$, the identity map on $X$ induces a weak$^*$ homeomorphic completely isometric $M$-$N$ bimodular isomorphism $X^{\natural \natural} \cong (X^{ \sharp })^*$. 
\end{prop}

If $E$ is a normal operator $M$-system, then by a representation of $E$ we mean an $M$-bimodular u.c.p.\ map $\pi: E \to \B(\mathcal H)$ where $M$ is normally represented on $\mathcal H$. We similarly define the space ${\rm Rep}(E)$ and the universal representation $\pi_u = \oplus_{\pi \in {\rm Rep}(E)} \pi$. 

\begin{prop}
The universal representation of a normal operator $M$-system $E$ induces a weak$^*$ homeomorphic isometric $M$-$M$ bimodular mapping of $(E^{ \sharp })^*$ onto the ultraweak closure $\overline{\pi_u(E)}$. In particular, $(E^{ \sharp })^*$ carries the structure of a dual normal operator $M$-system. 
\end{prop}
\begin{proof}
To each $\varphi \in E^{ \sharp }$ with $\| \varphi \| \leq 1$ there exists a normal state $\eta \in M_{*, +}$ and an $M$-bimodular completely bounded map $\psi: E \to \B(L^2(M, \eta))$ such that $\varphi(x) = \langle \psi(x) \hat{1}, \hat{1} \rangle$. Indeed, if $\omega, \rho \in M_{*, +}$ are such that $\abs{ \varphi(a^* x b) } \leq \omega(a^*a)^{1/2} \| x \| \rho(b^*b)^{1/2}$, then setting $\eta = (\omega(1) + \rho(1) )^{-1} ( \omega + \rho )$ we have that for each $x \in E$ the map $M \times M \ni (a^*, b^*) \mapsto \varphi(a^* x b)$ extends to a bounded sesquilinear map on $L^2(M, \eta) \times L^2(M, \eta)$, and hence there is a bounded operator $\psi(x) \in \B(L^2(M, \eta))$ such that $\varphi(x) = \langle \psi(x) \hat{1}, \hat{1} \rangle$ for $x \in M$. The map $E \ni x \mapsto \psi(x)\in\B(L^2(M, \eta))$ is then easily seen to be $M$-bimodular and bounded, hence also completely bounded by Smith's theorem.
 
By Wittstock's theorem, there then exist $M$-bimodular completely positive maps $\tilde \psi_i: E \to \B(L^2(M, \eta))$, $i = 1, 2$ such that $\tilde \psi_1 \pm \Re( \psi)$ and $\tilde \psi_2 \pm \Im( \psi)$ are completely positive. Moreover, if $E$ is concretely realized in $\B(\mathcal H)$, then by Arveson's extension theorem we may assume that $\psi$ and $\tilde \psi$ are defined on $\B(\mathcal H)$.

Considering the corresponding Stinespring dilations associated to these maps, we see that $\varphi$ is implemented by vectors in the $M$-system universal representation $\pi_u$. We therefore see that the universal representations for $E$ as an $M$-system or as an $M$-bimodule coincide and the result follows. 
\end{proof}

By the previous proposition, we have that, as a dual normal operator $M$-system, $(E^{ \sharp })^*$ satisfies the universal property that if $\pi: E \to F$ is an $M$-bimodular u.c.p.\  map where $F$ is a dual normal operator $M$-system, then there is a unique normal $M$-bimodular u.c.p.\  map $\tilde \pi: (E^{ \sharp })^* \to F$ such that $\tilde \pi_{| E} = \pi$.

An $M$-$C^*$-algebra consists of a $C^*$-algebra $A$, together with a faithful unital $*$-homomorphism mapping $M$ into the multiplier algebra $M(A)$. We say that the $M$-$C^*$-algebra is normal if $M(A)$ is normal as an operator $M$-system.

We let $p_{\rm nor} \in M(A)^{**}$ denote the supremum of the support projections of states in $M(A)^{*}$ that restrict to normal states on $M$, so that $M$ may be viewed as a von Neumann subalgebra of $p_{\rm nor} M(A)^{**} p_{\rm nor}$, and we have a canonical identification between $(M(A)^{\sharp})^*$ and the predual of $p_{\rm nor} M(A)^{**} p_{\rm nor}$. 
Equivalently, $p_{\rm nor}\in M^{**}\subset M(A)^{**}$ is the projection such that $M^{**} p_{\rm nor}$ is naturally isomorphic to $M$ and
a state $\varphi\in M(A)^\sharp$ if and only if $\varphi_{\mid M}$ is normal if and only if $\varphi(p_{\rm nor})=1$.
We let $q_A \in \mathcal P( M(A)^{**})$ denote the central projection in $M(A)^{**}$ that gives the identity projection for $A^{**} \subset M(A)^{**}$. We then have that the map $M \ni x \mapsto p_{\rm nor} q_A M(A)^{**} q_A p_{\rm nor} = p_{\rm nor} A^{**} p_{\rm nor}$ gives an embedding of $M$ as a von Neumann subalgebra of $p_{\rm nor} A^{**} p_{\rm nor}$. We therefore may think of $(A^{\sharp})^* \cong p_{\rm nor} A^{**} p_{\rm nor}$ as a von Neumann algebra that contains $M$ as a von Neumann subalgebra (see also \cite{BoCa15} where similar techniques are introduced).

If $A$ is a normal $M$-$C^*$-algebra and $B \subset A$ is a $C^*$-subalgebra that is an $M$-$M$ bimodule such that the natural $*$-homomorphism $M \to M(B)$ is faithful, then $B$ is also a normal $M$-$C^*$-algebra. If we let $q_B$ denote the support of $B^{**}$ in $A^{**}$, then $q_B$ and $p_{\rm nor}$ commute and we obtain an embedding (again as a non-unital von Neumann subalgebra) $(B^{\sharp})^* \cong q_B p_{\rm nor} A^{**} p_{\rm nor} q_B \subset p_{\rm nor} A^{**} p_{\rm nor} \cong (A^{\sharp})^*$.

\begin{rem}
If $A$ is an $M$-$C^*$-algebra, then $(A^{\sharp})^*$ is a von Neumann algebra that contains $M$ as a von Neumann subalgebra. However, one difference between this ``normal enveloping von Neumann algebra'' and the more familiar $A^{**}$ is that while the canonical embedding $A \to (A^{\sharp})^*$ is a complete order isomorphism onto its range (and hence we may think of $A$ as an $M$-subsystem of $(A^{\sharp})^*$), it is not a $*$-homomorphism in general, due to the fact that the projection $p_{\rm nor} \in M(A)^{**}$ need not be central. As a consequence, if $q_Q$ denotes the support of $A^{**}$ in $M(A)^{**}$, then, in general, we have $q_Q ( M(A)^{\sharp} )^* q_Q \cap M(A) \not= A$. 
\end{rem}

In light of Proposition~\ref{prop:quadmodrep}, the same arguments above regarding the universal representations of normal operator $M$-$N$ bimodule, normal $M$-systems, and normal $M$-$C^*$-algebras work in the setting where our object is also a $\tilde M$-$\tilde N$ bimodule, normal $\tilde M$-system, and normal $\tilde M$-$C^*$-algebra respectively. The main example relevant in this work is the case when $\tilde M = M' \subset \B(L^2M)$, in which case we may consider the space $\B(L^2M)^{M \sharp M} \cap \B(L^2M)^{M' \sharp M'}$. For use in the sequel, we end this section by recording these results.  

Suppose that $\mathcal H$ is a Hilbert space and $M, \tilde M, N, \tilde N \subset \B(\mathcal H)$ are von Neumann algebras such that $\tilde M \subset M'$ and $\tilde N \subset N'$. Suppose $X \subset \B(\mathcal H)$ is an operator space that is both an $M$-$N$ bimodule and a $\tilde M$-$\tilde N$ bimodule. Then we may also consider the space of representations in this category, which consists of complete contractions $\phi: X \to \B(\mathcal K)$, such that $M, \tilde M, N$ and $\tilde N$ are normally represented in $\B(\mathcal K)$, with $\tilde M \subset M'$, and $\tilde N \subset N'$, and such that $\phi$ is both $M$-$N$ bimodular and $\tilde M$-$\tilde N$ bimodular. 

\begin{prop}\label{prop:quadmoduniv}
Let $M, N, \tilde M, \tilde N \subset \B(\mathcal H)$ and $X \subset \B(\mathcal H)$ be as above. Then the universal representation induces a weak$^*$-homeomorphic completely isometric $M$-$N$ and $\tilde M$-$\tilde N$ bimodular mapping of $( X^{M \sharp N} \cap X^{\tilde M \sharp \tilde N} )^*$ onto the ultraweak closure $\overline{\pi_u(X)}$. 

Moreover, if $X$ is an operator system and contains $M = N$, and $\tilde M = \tilde N$, then $\overline{\pi_u(X)}$ is a both a dual normal $M$-system and a dual normal $\tilde M$-system.
\end{prop}

\begin{prop}\label{prop:MCalgebra}
Let $M$ and $\tilde M$ be von Neumann algebras and let $A$ be a $C^*$-algebra together with faithful unital $*$-homomorphisms of $M$ and $\tilde M$ into $M(A)$ such that $\tilde M \subset M' \cap M(A)$. Then considering the canonical inclusions $(A^{M \sharp M} \cap A^{\tilde M \sharp \tilde M})^*, (A^{M \sharp M})^*, (A^{\tilde M \sharp \tilde M})^* \subset M(A)^{**}$, we have that 
\[
(A^{M \sharp M} \cap A^{\tilde M \sharp \tilde M})^*=(A^{M \sharp M})^* \cap (A^{\tilde M \sharp \tilde M})^* \subset M(A)^{**}.
\]

Also, if $B \subset A$ is a $C^*$-subalgebra such that $B$ is both an $M$-bimodule and a $\tilde M$-bimodule and such that the corresponding representations of $M$ and $\tilde M$ into $M(B)$ are faithful, then we have a (non unital) inclusion of von Neumann algebras $(B^{M \sharp M} \cap B^{\tilde M \sharp \tilde M})^* \subset (A^{M \sharp M} \cap A^{\tilde M \sharp \tilde M})^*$.
\end{prop}

\section{Relatively compact operators}\label{sec:relcompact}

If $(M, \tau)$ is a finite von Neumann algebra, $\HH$ is a Hilbert space, and $x \in \B(L^2M, \HH)$, then we denote by $\| x \|_{\infty, 2}$ the norm of $x$ when viewed as an operator from $M \subset L^2M$ into $\HH$. If $\HH = L^2(N, \tau)$ for some finite von Neumann algebra $N$, then we let $\| x \|_{\infty, 1}$ denote the norm of $x$ when viewed as an operator from $M \subset L^2M$ into $L^1N \supset L^2N$. Ozawa proved in \cite{Oza10} that for $x \in \B(L^2M, \HH)$ we have
\begin{equation}\label{eq:infty2norm}
\| x \|_{\infty, 2} \leq \inf \{ \| z \| (\| c \|_2^2 + \| d \|_2^2 )^{1/2} \} \leq 4 \| x \|_{\infty, 2},
\end{equation}
where the infimum is taken over all decompositions $x = z \left( \begin{smallmatrix}
J c J \\ 
d
\end{smallmatrix} \right)$, where $z \in \B(L^2M \oplus L^2M, \HH)$ and $c, d \in M$. The following proposition adapts Ozawa's argument for the case of the $\| \cdot \|_{\infty, 1}$-norm.

\begin{prop}\label{prop:infty 1 norm}
Suppose $(M, \tau)$ is a finite von Neumann algebra and $x \in \B(L^2M)$. Then 
\[
\| x \|_{\infty, 1} \leq \inf \left\{ \left( \| a \|_2^2 + \| b \|_2^2 \right)^{1/2} \| z \| \left(\| c \|_2^2 + \| d \|_2^2 \right)^{1/2} \right\} \leq 4 \| x \|_{\infty, 1}
\]
where the infimum is over all $z \in \mathbb M_2(\B(L^2M))$, and $a, b, c, d \in M$ such that $x = \left( \begin{smallmatrix}
J a J \\
b
\end{smallmatrix}
\right)^*  z \left( \begin{smallmatrix}
J c J \\ 
d
\end{smallmatrix} \right)$.
\end{prop}
\begin{proof}
If $x = \left( \begin{smallmatrix}
J a J \\
b
\end{smallmatrix}
\right)^*  z \left( \begin{smallmatrix}
J c J \\ 
d
\end{smallmatrix} \right)$
such that $z \in \mathbb M_2(\B(L^2M))$, and $a, b, c, d \in M$, then for all $e, f \in M$ we have 
\begin{align*}
\abs{ \langle x \hat{e}, \hat{f} \rangle }
& = \left| \left\langle z \left( \begin{smallmatrix}
\widehat{ec^*} \\
\widehat{de}
\end{smallmatrix} \right), 
\left( \begin{smallmatrix}
\widehat{f a^*} \\
\widehat{bf}
\end{smallmatrix} \right) \right\rangle \right| \\
& \leq \| z \| \left( \| ec^* \|_2^2 + \| de \|_2^2 \right)^{1/2} \left( \| fa^* \|_2^2 + \| b f \|_2^2 \right)^{1/2} \\
& \leq \lVert e \rVert \lVert f \rVert \lVert z \rVert \left( \lVert c^* \rVert_2^2 + \lVert d \rVert_2^2 \right)^{1/2} \left( \lVert a^* \rVert_2^2 + \lVert b  \rVert_2^2 \right)^{1/2}.
\end{align*}
Taking the supremum over all $e, f \in M$ such that $\| e \|, \| f \| \leq 1$ then shows that 
\[
\| x \|_{\infty, 1} \leq \left( \| a \|_2^2 + \| b \|_2^2 \right)^{1/2} \| z \| \left(\| c \|_2^2 + \| d \|_2^2 \right)^{1/2}.
\]

For the second inequality, consider the normal bilinear form on $M \times M$ given by $\varphi(a, b) = \langle x \hat{a}, \hat{b} \rangle$ and note that we have $\| x \|_{\infty, 1} = \| \varphi \|$. By the noncommutative Grothendieck Theorem, there exist unit vectors $\xi_1, \xi_2, \eta_1, \eta_2 \in L^2M$ such that 
\[
\abs{ \varphi(a, b) }
\leq \| \varphi \| \left( \| a \xi_1 \|_2^2 + \| \xi_2 a \|_2^2 \right)^{1/2} \left( \| b \eta_1 \|_2^2 + \| \eta_2 b \|_2^2 \right)^{1/2}.
\]

Thinking of $L^2M$ as a subspace of unbounded affiliated operators, we let $p_1 = \chi_{[\| x \|^4, \infty)}(\xi_1 \xi_1^*)$. We then have 
\begin{align*}
\abs{ \varphi(e, f) }
& \leq \abs{ \varphi(e p_1^\perp,  f ) + \varphi(e p_1, f) } \\
& \leq \| \varphi \| \left( \| e p_1^\perp \xi_1 \|_2^2 + \| \xi_2  e \|_2^2 \right)^{1/2} \left( \| f  \eta_1 \|_2^2 + \| \eta_2  f \|_2^2 \right)^{1/2}
+  \| x \| \| e p_1 \|_2 \| f \|_2 \\
& \leq \| \varphi \| \left( \| e p_1^\perp \xi_1 \|_2^2 + \| \xi_2  e \|_2^2  + \| x \|^2 \| e p_1 \|_2^2 \right)^{1/2} \left( \| f  \eta_1 \|_2^2 + \| \eta_2  f \|_2^2  + \| f \|_2^2 \right)^{1/2} \\
&= \| \varphi \| \left( \| e a_0 \|_2^2 + \| \xi_2 e \|_2^2 \right)^{1/2} \left( \| f \tilde \eta_1 \|_2^2 + \| \eta_2 f \|_2^2 \right)^{1/2},
\end{align*}
where $a_0 = ( p_1^\perp \xi_1 \xi_1^* p_1^\perp + \| x \|^4 p_1 )^{1/2}$, and $\tilde \eta_1 = (\eta_1 \eta_1^* + 1 )^{1/2}$. Note that $\| a_0 \| \leq \| x \|^2$ and $\| a_0 \|_2 \leq \| \xi_1 \|_2 \leq 1$, while $\| \tilde \eta_1 \|_2^2 \leq \| \eta_1 \|_2^2 + 1 = 2$. 

Similarly, if we let $p_2 = \chi_{[ \| x \|^4, \infty) }( \xi_2^* \xi_2)$, then we see that 
\[
\abs{ \varphi(e, f) } \leq \| \varphi \| \left( \| e a_0 \|_2^2 + \| b_0 e \|_2^2 \right)^{1/2} \left( \| f \tilde \eta_1 \|_2^2 + \| \tilde \eta_2 f \|_2^2 \right)^{1/2},
\]
where $b_0 = (p_2^\perp \xi_2^* \xi_2 p_2^\perp + \| x \|^4 p_2)^{1/2}$ satisfies $\| b_0 \| \leq \| x \|^2$ and $\| b_0 \|_2 \leq 1$, and where $\tilde \eta_2 = (\eta_2^* \eta_2 + 1)^{1/2}$ satisfies $\| \tilde \eta_2 \|_2 \leq \sqrt{2}$. 

Repeating this argument with $\tilde \eta_1$ and $\tilde \eta_2$, we obtain an inequality of the form
\begin{equation}\label{eq:quad}
\abs{ \varphi(e, f) } \leq \| \varphi \| \left( \| e a \|_2^2 + \| b e \|_2^2 \right)^{1/2} \left( \| f c \|_2^2 + \| d f \|_2^2 \right)^{1/2},
\end{equation}
where $a, b, c, d \in M$ are bounded with $\| a \|_2, \| b \|_2, \| c \|_2, \| d \|_2 \leq \sqrt{2}$. 

If we let $V \subset L^2M \oplus L^2M$ denote the closed subspace spanned by elements of the form $ea \oplus be$, for $e \in M$, and if we let $W \subset L^2M \oplus L^2M$ denote the closed subspace spanned by elements of the form $fc \oplus df$, for $f \in M$, then from (\ref{eq:quad}) we see that $\varphi(e, f)$ defines a bounded quadratic form on $V \times W$ and we therefore obtain a bounded operator $z \in \mathbb M_2(\B(L^2M))$ with $\| z \| \leq \| \varphi \|$ such that for all $e, f \in M$ we have
\[
\varphi(e, f) = \langle z ( ea \oplus be ), ( fc \oplus df ) \rangle.
\]
Hence we have $x = \left( \begin{smallmatrix}
J a J \\
b
\end{smallmatrix}
\right)^*  z \left( \begin{smallmatrix}
J c J \\ 
d
\end{smallmatrix} \right)$.
\end{proof}

\begin{cor}\label{cor:M-topology and norms}
Let $M$ and $N$ be finite von Neumann algebras, and $X \subset \B(L^2N, L^2M)$ a subspace. If $T \in \B(L^2N, L^2M)$, and $\{ T_i \}_i \subset X$ is uniformly bounded such that $\| T_i - T \|_{\infty, 1} \to 0$, then $T$ is in the smallest subspace of $\B(L^2N, L^2M)$ that contains $X$ and is both a strong $M$-$N$ bimodule and a strong $JMJ$-$JNJ$ bimodule.

Also, if $\HH$ is a Hilbert space, $Y \subset \B(L^2N, \HH)$ is a subspace, $T \in \B(L^2M, \HH)$, and $\{ T_i \}_i \subset Y$ is uniformly bounded such that $\| T - T_i \|_{\infty, 2} \to 0$, then $T$ is in the smallest subspace of $\B(L^2N, \HH)$ that contains $Y$ and is both a strong right $N$ module and a strong right $JNJ$ module.
\end{cor}
\begin{proof}
Note there is no loss in generality by assuming that $X$ is an $M$-$N$ bimodule, and a $JMJ$-$JNJ$ bimodule. By Proposition~\ref{prop:quadmodrep} we have 
\[
\B(L^2N, L^2M)^{M \sharp N} \cap \B(L^2N, L^2M)^{JMJ \sharp JNJ} = \B(L^2N, L^2M)^{M \sharp JNJ} \cap \B(L^2N, L^2M)^{JMJ \sharp N}.
\]
Using this fact, we then see that Propositions~\ref{prop:quadmodstrong} and \ref{prop:infty 1 norm} give the result for the case $X \subset \B(L^2N, L^2M)$. The case $Y \subset \B(L^2N, \HH)$ follows similarly by using Ozawa's Lemma in \cite{Oza10} instead of Proposition~\ref{prop:infty 1 norm}.
\end{proof}

Let $M$ be a von Neumann algebra. If $\X$ is a hereditary $C^*$-subalgebra of $\B(L^2M)$, then we denote by $M(\X)$ the multiplier algebra in $\B(L^2M)$.  An $M$-boundary piece is a hereditary $C^*$-subalgebra $\X \subset \B(L^2M)$ such that $M(\X) \cap M$ and $M(\X) \cap JMJ$ are ultraweakly dense in $M$ and $JMJ$ respectively. To avoid pathological examples, we will always assume that $\X \not= \{ 0 \}$, so that we then have $\K(L^2M) \subset \X$. The following are two motivating examples for this definition.

 \begin{examp}[See \cite{BoIoPe21} and Section 15.1 in \cite{BrOz08}]
 Suppose $\Gamma$ is a group and $I \subset \ell^\infty \Gamma$ is a closed ideal that is invariant under the left and right actions of $\Gamma$ and contains $c_0 \Gamma$. Then $\X_I = I \B(\ell^2 \Gamma) I$ gives a boundary piece for $L\Gamma$. Note that $C^*_\lambda \Gamma$ and $C^*_\rho \Gamma$ are contained in the multiplier algebra of $I \B(\ell^2 \Gamma) I$
 \end{examp}

 \begin{examp}\label{examp:bdrysub}
 Suppose $B \subset M$ is a von Neumann subalgebra and let $e_b \in \B(L^2M)$ denote the orthogonal projection onto the space $L^2B \subset L^2M$. Then the algebraic span of operators of the form $x_1 J y_1 J T J y_2 J x_2$, with $x_1, x_2, y_1, y_2 \in M$ and $T \in e_B \B(L^2M) e_B \cong \B(L^2B)$, forms a $*$-subalgebra, and its closure then forms a $C^*$-subalgebra of $\B(L^2M)$, which we denote by $\X_B$. The $C^*$-algebra $\X_B$ clearly contains $M$ and $JMJ$ in its multiplier algebra, and to see that it is hereditary we just note that it can be identified with the hereditary $C^*$-subalgebra of $\B(L^2M)$ that is generated by all operators of the form $x JyJ e_B$ with $x, y \in M$. 
 \end{examp}

We now fix an $M$-boundary piece $\X$ and let $\K_{\X}^L = \K_{\X}^L(M) \subset \B(L^2M)$ denote the $\| \cdot \|_{\infty, 2}$-closure of the closed left ideal $\B(L^2M) \X$, i.e., 
\[
\K_{\X}^L = \overline{ \B(L^2M) \X }^{\| \cdot \|_{\infty, 2}}.
\]

\begin{prop}\label{prop:leftideal}
The space $\K_{\X}^L$ is a closed left ideal in $\B(L^2M)$ such that $M$ and $JMJ$ are contained in the space of right multipliers. 
\end{prop}
\begin{proof}
First notice that $\|\cdot \|_{\infty, 2}\leq \|\cdot\|$ on $\B(L^2M)$ and hence $\K_\X^L$ is closed.
The fact that $\K_\X^L$ is a left ideal follows from the inequality $\|ST\|_{\infty, 2}\leq \|S\| \|T\|_{\infty, 2}$, for any $S\in\B(L^2M)$ and $T\in \K_\X^L$.

To see $M$ is in the right multiplier of $\K^L_\X$, let $x\in M$ and $T\in \K^L_\X$ be given. Then, using Kaplansky's theorem, there exist sequences $\{T_n\}\subset \B(L^2M)\X$ and ${x_n}\subset M(\X)\cap M$ such that $\|T_n-T\|_{\infty, 2}\to 0$, $\|x_n\|\leq \|x\|$ and $\|x_n-x\|_2\to 0$.
Notice that $T_nx_n\in \B(L^2M)\X$ and
\begin{align*}
\|T_nx_n-Tx\|_{\infty,2}&\leq \|(T_n-T)x_n\|_{\infty,2}+\|T(x-x_n)\|_{\infty,2}\\
&\leq \|T_n-T\|_{\infty,2}\|x\|+\|T\|\|x_n-x\|_2 \to 0,
\end{align*}
and thus $Tx\in \K^L_\X$. It similarly follows that $JMJ$ is in the space of right multipliers.
\end{proof}

We let $\K_{\X} = \K_{\X}(M)$ denote the hereditary $C^*$-algebra associated to $\K_{\X}^L$, i.e.,
\begin{equation}\label{eq:cptpiece}
\K_{\X} = ( \K^L_{\X})^* \cdot \K^L_{\X} = ( \K^L_{\X})^* \cap \K^L_{\X}.
\end{equation}
Note that from Proposition~\ref{prop:leftideal} we have that $M$ and $JMJ$ are contained in the multiplier algebra of $\K_{\X}$. We let $\K_{\X}^{\infty, 1} = \K_{\X}^{\infty, 1}(M) \subset \B(L^2M)$ denote the closure of $\X$ in the $\| \cdot \|_{\infty, 1}$ norm, i.e.,
\begin{equation}\label{eq:cpt1}
\K_{\X}^{\infty, 1} = \overline{ \X}^{\| \cdot \|_{\infty, 1}}.
\end{equation}

\begin{prop}\label{prop:strong K_X}
We have $\K_{\X}^{\infty, 1} = \overline{ \K_{\X} }^{\| \cdot \|_{\infty, 1}}$, and this is a self-adjoint closed strong $M$-$M$ and $JMJ$-$JMJ$ bimodule. 
\end{prop}
\begin{proof}
It is clear that $\X\subset \K_\X$ and hence it suffices to show $\overline{ \K_{\X} }^{\| \cdot \|_{\infty, 1}}\subset \K_{\X}^{\infty, 1}$.
For each $T\in \overline{ \K_{\X} }^{\| \cdot \|_{\infty, 1}}$, there exists $\{S_n\}$, $\{R_n\}\subset \K_\X^L$ such that $\|S_n^*R_n-T\|_{\infty, 1}\to 0$. Furthermore, for each $n$ we have sequences $\{S_n^i\}$ and $\{R_n^i\}$ in $\B(L^2M)\X$ that converge to $S_n$ and $R_n$ in $\|\cdot\|_{\infty,2}$, respectively.
By the polarization identity, we have $(S^i_n)^* R_n^i\in \X$ for each $i$ and $n$.
Note that for each $n$ we have
\begin{align*}
	\| (S_n^i)^*R_n^i-T\|_{\infty, 1}&\leq \|((S_n^i)^*-S_n^*)R_n^i\|_{\infty, 1}+\| S_n^*(R_n^n-R_n)\|_{\infty, 1} +\|S_n^*R_n-T\|_{\infty,1}\\
	&\leq \|S_n-S_n^i\|_{\infty,2}\|R_n^i\|_{\infty,2}+\|S_n\|_{\infty, 2}\|R_n^i-R_n\|_{\infty,2}+\|S_n^*R_n-T\|_{\infty,1},
\end{align*}
i.e., $\lim_{i\to\infty} \| (S_n^i)^*R_n^i-T\|_{\infty, 1}\leq \|S_n^*R_n-T\|_{\infty,1}$. Thus for each $n$, we may pick $i(n)$ such that $\| (S_n^{i(n)})^*R_n^{i(n)}-T\|_{\infty, 1}\leq \|S_n^*R_n-T\|_{\infty,1}+2^{-n}$ and therefore $\lim_{n\to\infty}\| (S_n^{i(n)})^*R_n^{i(n)}-T\|_{\infty, 1}=0$.

The fact that it is closed follows from the inequality $\|\cdot\|_{\infty, 1}\leq \|\cdot\|$ on $\B(L^2M)$; since one also has $\|T\|_{\infty,1}=\|T^*\|_{\infty,1}$ for any $T\in \B(L^2M)$, $\K_{\X}^{\infty, 1}$ is self-adjoint.

Since we have now established that $\K_{\X}^{\infty, 1} = \overline{\K_\X}^{\| \cdot \|_{\infty, 1} }$, Proposition~\ref{prop:quadmodstrong} and Corollary~\ref{cor:M-topology and norms} then show that $\K^{\infty, 1}_\X=\overline{\overline{\K_\X}^{_{M\text -M}}}^{_{JMJ \text -JMJ}}$, which is a strong $M$-$M$ and $JMJ$-$JMJ$ bimodule.
\end{proof}

\begin{examp}
The space $\K^{\infty, 1}(M)$ is, in general, considerably larger than $\K(L^2M)$. For example, $\K^{\infty, 1}(M)$ may contain isometries. Indeed, if $\{ p_n \}_{n \in \mathbb N} \subset \mathcal P(M)$ is an infinite partition of unity consisting of non-zero projections, and we let $V \in \B(L^2M)$ be an isometry from $L^2M$ onto the closed span of $\{ p_n \}_n \subset L^2M$, and if we set $q_k = \sum_{n = k}^\infty p_n$, then we have $q_k^\perp V \in \K(L^2M)$ for each $k \geq 1$, and hence by Proposition~\ref{prop:strong K_X} we have $V \in \K^{\infty, 1}(M)$.
\end{examp}

If $\HH$ is a Hilbert space, then we also let $\K_\X^L(M, \HH) =  \overline{ \B(L^2M, \HH) \X }^{\| \cdot \|_{\infty, 2}}$. Note that by considering polar decomposition we have $T \in \K_\X^L(M, \HH)$ if and only if $\abs{ T } \in \K_\X^L(M)$. 

\begin{prop}\label{prop:positivecpt}
Let $M$ be a finite von Neumann algebra, $\X \subset \B(L^2M)$ a boundary piece, and $\HH$ a Hilbert space. If $T \in \B(L^2M, \HH)$, then the following conditions are equivalent:
\begin{enumerate}
\item\label{item:pos1} $T \in \K_\X^L(M, \HH)$.
\item\label{item:pos2} $T^*T \in \K_\X$. 
\item\label{item:pos3} $T^*T \in \K_\X^{\infty, 1}$.
\end{enumerate}
\end{prop}
\begin{proof}
The implications (\ref{item:pos1})$\implies$(\ref{item:pos2})$\implies$(\ref{item:pos3}) are clear.
To see (\ref{item:pos3})$\implies$(\ref{item:pos1}), recall from the proof of Proposition~\ref{prop:strong K_X} that $\K_\X^{\infty, 1}= \overline{\overline{\K_\X}^{_{M\text -M}}}^{_{JMJ\text -JMJ}}$.
It follows from Proposition~\ref{prop:quadmodstrong} that there exist sequences of projections $\{p_n\}_n\subset\mathcal P(M)$ and $\{q_n\}_n\subset \mathcal P(JMJ)$ with $p_n\to 1$ and $q_n\to 1$ such that $p_n q_n T^*T q_n p_n\in \K_\X$.
In particular, we then have $\abs{ T } q_np_n\in\K_\X^L$ (see Proposition II.5.3.2 in \cite{Bl06}), and 
it then follows that $|T|\in\overline{\overline{\K_\X^L}^{_{\C-M}}}^{_{\C-JMJ}}$, which coincides with $\K_\X^L$ by Proposition~\ref{prop:quadmodstrong} and Corollary~\ref{cor:M-topology and norms}.
By the remark before this proposition, we then have $T\in \K_\X^L(M, \HH)$.
\end{proof}

\begin{lem}\label{lem:cptl2}
Let $M$ and $N$ be finite von Neumann algebras and let $\phi: N \to M$ be a bounded map such that $\phi$ has a continuous extension $T_\phi \in \B(L^2N, L^2M)$. Then $T_\phi$ is compact as an operator from $N$ into $L^1M$ if and only if $T_\phi$ is compact as an operator from $N$ into $L^2M$.
\end{lem}
\begin{proof}
Suppose $T_\phi$ is compact as an operator from $N$ into $L^1M$.  Since $T_\phi^*: M \to L^1N$ is also compact, and since $\phi: N \to M$ is bounded, it follows that $T_\phi^* T_\phi$ is compact as an operator from $N$ to $L^1N$. 

Hence for any bounded sequence $\{ x_n \}_{n \in \N} \subset N$ such that $x_n \to 0$ ultraweakly we have 
\[
\| T_\phi( x_n) \|_2^2 = \langle T_\phi^* T_\phi \widehat{x_n}, \widehat{x_n} \rangle \leq \| T_\phi^* T_\phi( \widehat{x_n} ) \|_1 \to 0.
\] 
Therefore $T_\phi$ is compact as an operator from $N$ into $L^2N$. 
\end{proof}

\begin{proof}[Proof of Theorem~\ref{thm:AD}]
The implications (1) $\implies$ (2) $\implies$ (3) are obvious, and we will show (3) $\implies$ (1).  Let $\{ \phi_i \}_i$ be given as in (3). By a standard convexity argument (e.g., \cite[Theorem 2.1]{OzPo10I}), we may assume that $\phi_i$ satisfy $\phi_i(1) \leq 1$ and $\tau \circ \phi_i \leq \tau$. By Kadison's inequality, we then have that that each $\phi_i$ induces a bounded operator on $L^2M$, and by Lemma~\ref{lem:cptl2} each $\phi_i$ is compact as an operator from $M$ to $L^2M$ and hence, as in \cite{Oza10}, each $\phi_i$ is in the $\| \cdot \|_{\infty, 2}$-closure of $\K(L^2M)$. By Corollary~\ref{cor:M-topology and norms}, for each $i$ and $\varepsilon > 0$ there exists a projection $p_{i, \varepsilon} \in \mathcal P(M)$ with $\tau(p_{i, \varepsilon}) > 1 - \varepsilon$ such that the c.p.\ map $\psi_{i, \varepsilon} =  \phi_i \circ{\rm Ad}(p_{i, \varepsilon})$ gives a compact operator on $L^2M$. Letting $i \to \infty$ and $\varepsilon \to 0$ gives a net of c.p.\ maps satisfying (1). (Note that $\phi_i$ will converge pointwise to the identity in $\| \cdot \|_2$ because of the inequality $\| \phi_i(x) - x \|_2^2 \leq 2 ( \| x \|_2^2 - \Re \langle \phi_i(\hat{x}), \hat{x} \rangle )$.
\end{proof}

In Section~\ref{sec:mixhilbert} we will give another perspective on Theorem~\ref{thm:AD} through the use of mixing Hilbert bimodules.

\section{The convolution algebra associated to a von Neumann algebra}\label{sec:conv}

Let $M$ and $N$ be von Neumann algebras. We let $M \otimes_{Con} N^{\rm op}$ be the subset of $M \otimes_{eh} N = CB_{M'-N'}^\sigma(\B(L^2N, L^2M))$ consisting of normal $M'$-$N'$ bimodular completely bounded maps that preserve the space of trace-class operators\footnote{The change in notation from $N$ to $N^{\rm op}$ is to emphasize the fact that we will mainly consider $M \otimes_{Con} N^{\rm op}$ not as a subspace of $M \otimes_{eh} N$, but rather as an algebra that contains the algebraic tensor product $M \otimes N^{\rm op}$ as a subalgebra.}. If $\mu \in M \otimes_{Con} N^{\rm op}$, then we let $\mu^* \in CB_{M'-N'}^\sigma(\B(L^2N, L^2M)) = M \otimes_{eh} N$ denote the adjoint of the operator $\mu$ when restricted to the space of trace-class operators under the usual conjugate linear pairing between $\B(L^2N, L^2M)$ and the space of trace-class operators in $\B(L^2N, L^2M)$, i.e., ${\rm Tr}(\mu(S) T^*) = {\rm Tr}( S (\mu^*(T))^*)$ for each $S, T \in \B(L^2N, L^2M)$ with $S$ trace-class. Note that since $\mu$ is continuous and preserves the trace-class operators, then $\mu$ also preserves the space of compact operators and hence the dual map $\mu^*$ again preserves the trace-class operators, thus the map $M \otimes_{Con} N^{\rm op} \ni \mu \mapsto \mu^* \in M \otimes_{Con} N^{\rm op}$ defines an involutive antilinear algebra-isomorphism.

Note that the involution $\mu \mapsto \mu^*$ extends the usual involution on the algebraic tensor product $M \otimes N^{\rm op}$. We define a norm $\| \cdot \|_{Con}$ on $M \otimes_{Con} N^{\rm op}$ by 
\[
\| \mu \|_{Con} = \max\{ \| \mu \|_{M \otimes_{eh} N}, \| \mu^* \|_{M \otimes_{eh} N} \}.
\] 
We denote by $J$ both the Tomita conjugation operator for $M$ and for $N$. We let $\mathcal J: \B(L^2N, L^2M) \to \B(L^2N, L^2M)$ denote the anti-linear isometry given by $\mathcal J(T) = J T J$. We then have that $M \otimes_{eh} N \ni \nu \mapsto \mathcal J \nu \mathcal J$ gives an isometric anti-isomorphism between $M \otimes_{eh} N$ and $(J M J) \otimes_{eh} (J N J)$. We may therefore view $M \otimes_{Con} N^{\rm op}$ as subspace of the $\ell^\infty$-direct sum 
\[
(M \otimes_{eh} N) \oplus^\infty ( (JMJ) \otimes_{eh} (JNJ))
\]
under the diagonal embedding $M \otimes_{Con} N^{\rm op} \ni \mu \mapsto \mu \oplus {\mathcal J \mu^* \mathcal J}$. It is easy to check that this is a weak$^*$-closed (and hence also norm closed) subspace, so that $M \otimes_{Con} N^{\rm op}$ is a dual Banach space where the weak$^*$-topology on bounded sets is given by point-ultraweak convergence applied to compact operators for a net $\{ \mu_i \}_i$ and its adjoint net $\{ \mu_i^* \}_i$.

Since $M \otimes_{eh} N = CB_{M'-N'}^\sigma(\B(L^2N, L^2M))$ is a Banach algebra under composition, it then follows that $M \otimes_{Con} N^{\rm op}$ is a Banach $*$-algebra. Moreover, since the maps in $\mu \in M \otimes_{Con} N^{\rm op}$ preserve the space of compact operators, it follows that multiplication in $M \otimes_{Con} N^{\rm op}$ is separately weak$^*$-continuous in each variable, and so we view $M \otimes_{Con} N^{\rm op}$ as a dual Banach $*$-algebra.

\subsection{Convolution algebras associated to finite von Neumann algebras}\label{subs:con}

Suppose now that $M$ and $N$ are finite von Neumann algebras with normal faithful traces. Let $\mathcal V$ be a normal Banach $M, N$-bimodule. Following Connes \cite{Co78} we consider the space $\mathcal V^{M \natural N} \subset \mathcal V^*$ consisting of continuous linear functionals $\varphi$ such that there exists $K \geq 0$ for which we have 
\begin{equation}\label{eq:con}
\abs{ \varphi( x^* \xi y) } \leq K \| x \|_2 \| \xi \| \| y \|_2, \ \ \ x \in M, \ \  \xi \in \mathcal V, \ \ y \in N.
\end{equation} 
When no confusion will arise, we will denote $\mathcal V^{M \natural N}$ simply by $\mathcal V^\natural$. 
We define the norm $\| \varphi \|_\natural$ to be the smallest $K \geq 0$ for which (\ref{eq:con}) holds. Note that $\| \varphi \| \leq \| \varphi \|_\natural$ and it then follows that the unit ball in $\mathcal V^\natural$ is compact for the weak-* topology $\sigma(\mathcal V^\natural, \mathcal V)$. Thus, if we define a semi-norm on $\mathcal V$ by $\| \xi \|_\flat = \sup_{ \varphi \in \mathcal V^\natural, \| \varphi \|_\natural \leq 1} \abs{ \varphi( \xi) }$, then with this semi-norm $\mathcal V$ has dense image in a Banach space $\mathcal V_\flat$ such that $\mathcal V^\natural$ is dual to $\mathcal V_\flat$. 

The space $\mathcal V^\natural$ inherits an $N, M$-bimodule structure from $\mathcal V^*$, i.e., if $\varphi \in \mathcal V^\natural$ and $x \in N$, $y \in M$, then $(x \varphi y)(\xi) = \varphi(y \xi x)$. Since $\| x a \|_2 \leq \| x \| \| a \|_2$ and $\| b y \|_2 \leq \| y \| \| b \|_2$ for $a \in N$, $b \in M$, we then have $\| x \varphi y \|_\natural \leq \| x \| \| \varphi \|_\natural \| y \|$ for $\varphi \in \mathcal V^\natural$, and hence also $\| x \xi y \|_\flat \leq \| x \| \| \xi \|_\flat \| y \|$ for $\xi \in \mathcal V_\flat$. Thus, $\mathcal V_\flat$ is a Banach $M, N$-bimodule and $\mathcal V^\natural$ is a dual Banach $N, M$-bimodule. Also, for fixed $\varphi \in \mathcal V^\natural$, $\xi \in \mathcal V$, $x \in M$, and $y \in N$, we have $\abs{ \varphi( \xi y) } \leq \| \varphi \|_\natural \| y \|_2 \| \xi \|$, and $\abs{ \varphi( x \xi ) } \leq \| \varphi \|_\natural \| x \|_2 \| \xi \|$ so that $\mathcal V^\natural$ is a normal dual Banach $N, M$-bimodule.

From (\ref{eq:con}) we see that for each $\varphi \in \mathcal V^\natural$ we obtain a map $\Psi_\varphi: \mathcal V \to \B(L^2N, L^2M)$ such that for all $\xi \in \mathcal V$ we have 
\[
\langle \Psi_\varphi(\xi) \hat y, \hat x \rangle = \varphi( x^* \xi y ).
\] 
Moreover, we see that $\Psi_\varphi$ is $M, N$-bimodular and $\| \Psi_\varphi \| = \| \varphi \|_\natural$. Conversely, to every $M, N$-bimodular bounded map $\Psi$ we obtain a linear functional $\varphi_{\Psi} \in \mathcal V^\natural$ by $\varphi_{\Psi}(\xi) = \langle \Psi(\xi) \hat 1, \hat 1 \rangle$. These maps are inverses of each other so that we have an isometric isomorphism between $\mathcal V^\natural$ and the bimodule dual considered in \cite{Ma05}.
\[
\mathcal V^\natural \cong \B_{M\text-N}(\mathcal V, \B(L^2N, L^2M)).
\] 
We may also view $\B(L^2N, L^2M)$ as an $N, M$-bimodule with the bimodule structure given by $x \cdot T \cdot y = Jy^*J T Jx^*J$. This then induces a bimodule structure on $\B_{M\text-N}(\mathcal V, \B(L^2N, L^2M))$ by post composition. It is then easy to see that $\mathcal V^\natural \ni \varphi \mapsto \Psi_{\varphi}$ is an $N, M$-bimodular map.

Moreover, for a uniformly bounded net $\{ \varphi_i \}_i$ we see that $\varphi_i \to \varphi$ weak$^*$ if and only if $\Psi_{\varphi_i}(\xi) \to \Psi_{\varphi}(\xi)$ in the ultraweak topology for each $\xi \in \mathcal V$. Thus, the map $\mathcal V^\natural \ni \varphi \mapsto \Psi_\varphi \in \B_{M, N}(\mathcal V, \B(L^2N, L^2M))$ gives an isomorphism of dual Banach $N, M$-bimodules 
\[
\mathcal V^\natural \cong \B_{M\text-N}(\mathcal V, \B(L^2N, L^2M)).
\] 
Note that the latter space is naturally a dual normal operator $N, M$-bimodule \cite{EfRu88}, and we endow $\mathcal V^\natural$ with this operator bimodule structure. 

If $\mathcal V$ is an operator $M, N$-bimodule, then by Smith's theorem \cite{Sm91} we have 
\[
\B_{M\text-N}(\mathcal V, \B(L^2N, L^2M)) = CB_{M\text-N}(\mathcal V, \B(L^2N, L^2M)),
\] 
and hence we may rewrite the semi-norm $\| \cdot \|_\flat$ on $\mathcal V$ as 
\begin{equation}\label{eq:ccnorm}
\| \xi \|_\flat = \sup_{\phi: \mathcal V \to \B(L^2N, L^2M), \atop M-N {\rm \ bimodular }, \| \phi \|_{cb} \leq 1} \abs{ \langle \phi(\xi) \hat{1}, \hat{1} \rangle }.
\end{equation}

It then follows from Wittstock's extension theorem that if $\mathcal W \subset \mathcal V$ is an operator $M, N$-sub-bimodule, then we have an isometric inclusion $\mathcal W_\flat \subset \mathcal V_\flat$.

From the discussion above, we obtain a completely isometric isomorphism 
\[
M \otimes_{eh} N = CB_{M'\text-N'}( \K(L^2N, L^2M), \B(L^2N, L^2M) ) \cong \K(L^2N, L^2M)^{M' \natural N'}.
\]
Note that $\K(L^2N, L^2M)$ and $\B(L^2N, L^2M)$ are also natural $M, N$-bimodules. We set 
\[
\B(L^2N, L^2M)^\natural_J = \B(L^2N, L^2M)^{M \natural N} \cap \B(L^2N, L^2M)^{M' \natural N'}
\] 
and 
\[
\K(L^2N, L^2M)^\natural_J = \K(L^2N, L^2M)^{M \natural N} \cap \K(L^2N, L^2M)^{M' \natural N'}.
\] 
We endow these spaces with their natural norms coming from interpolation theory. 

Restricting the isomorphism $M \otimes_{eh} N \cong \K(L^2N, L^2M)^{M' \natural N'}$ from above then gives an isometric isomorphism 
\[
M \otimes_{Con} N^{\rm op} \cong \K(L^2N, L^2M)^\natural_J,
\] 
and one can check easily that $\Psi_{\varphi^*}=(\Psi_\varphi)^*$ and $\varphi_{\Psi^*}=(\varphi_{\Psi})^*$ for any $\varphi\in  \K(L^2N, L^2M)^\natural_J$ and $\Psi\in M\otimes_{Con}N^{\rm op}$, where $\varphi^*(\cdot):=\overline{\varphi(J\cdot J)}$.

The advantage of viewing $M \otimes_{Con} N^{\rm op}$ as a subspace of the dual of $\K(L^2N, L^2M)$ is that we may then use techniques for linear functionals, e.g., Jordan decomposition, when working with elements in $M \otimes_{Con} N^{\rm op}$, see, e.g., Lemma~\ref{lem:approxbd} below or \cite{DaPe20} where similar techniques are developed.

\begin{lem}\label{lem:sharpdense}
Let $M$ and $N$ be tracial von Neumann algebras and suppose we have faithful normal representations $M \subset \B(\mathcal H)$ and $N \subset \B(\mathcal K)$, then $\K(\mathcal K, \mathcal H)^\natural$ is norm dense in $\K(\mathcal K, \mathcal H)^*$. 
\end{lem}
\begin{proof}
If $\xi \in \mathcal H$ and $\eta \in \mathcal K$ are unit vectors such that $\langle a \xi, \xi \rangle \leq K \| a \|_1$ for all $0 \leq a \in M$, and $\langle b \eta, \eta \rangle \leq L \| b \|_1$ for all $0 \leq b \in N$, and if we take $x \in M$ and $y \in N$, then if we denote by $V_{\xi, \eta}$ the rank-1 partial isometry mapping $\eta$ to $\xi$ we see that the rank-1 operator $x V_{\xi, \eta} y$ corresponds to a linear functional in $\K(\cK, \cH)^\natural$, and a simple computation shows that $\| x V_{\xi, \eta} y \|_{\natural} \leq K^{1/2} L^{1/2} \| x \| \| y \|$. The lemma is then established by simply noting that the span of such operators is dense in the space of trace-class operators in $\B(\cK, \cH)$. 
\end{proof}

\begin{lem}\label{lem:sotconvergence}
Let $M$ and $N$ be tracial von Neumann algebras, suppose we have a uniformly bounded net $\{ x_i \}_{i} \subset M$ with $\| x_i \|_2 \to 0$, and a uniformly bounded net $\{ \mu_i \}_i \subset M \otimes_{Con} N^{\rm op}$. If $\mathcal W$ is a normal dual operator $M$-$N$-bimodule, then for all $v \in \mathcal W$ we have weak$^*$-convergence $\mu_i \circ (x_i \otimes 1)(v) \to 0$ and $(x_i \otimes 1) \circ \mu_i(v) \to 0$.
\end{lem}
\begin{proof}
We suppose $M \subset \B(\mathcal H)$, $N \subset \B(\mathcal K)$ and $\mathcal W \subset \B(\mathcal K, \mathcal H)$ is a weak$^*$-closed $M$-$N$-bimodule. 

Suppose $v \in \mathcal W$ and $\varphi \in \K(\mathcal K, \mathcal H)^\natural$, then 
\begin{align*}
\abs{ \varphi( \mu_i \circ ( x_i \otimes 1) (v) ) }
&\leq \| v \| \| \varphi \circ \mu_i \|_\natural \| x_i \|_2 \\
&\leq \| v \| \| \varphi \|_\natural \| \mu_i \|_{Con} \| x_i \|_2
\to 0.
\end{align*}
We similarly have 
\begin{align*}
\abs{ \varphi((x_i \otimes 1) \circ \mu_i (v) ) }
&\leq \| \mu_i(v) \| \| \varphi \|_\natural \| x_i \|_2 \\
&\leq \| v \| \| \mu_i \|_{Con}  \| \varphi \|_\natural \| x_i \|_2
\to 0.
\end{align*}

By Lemma~\ref{lem:sharpdense}  we have that $\K(\mathcal K, \mathcal H)^\natural$ is dense in $\K(\mathcal K, \mathcal H)^*$ and hence we have ultraweak convergence $\mu_i \circ ( x_i \otimes 1) (v) \to 0$ and $(x_i^* \otimes 1) \circ \mu_i(v) \to 0$. 
\end{proof}

\section{Mixing operator bimodules}

\subsection{A relative topology on the convolution algebra}
If $M$ and $N$ are finite von Neumann algebras and if $\X \subset \B(L^2M)$ and $\Y \subset \B(L^2N)$ are boundary pieces for $M$ and $N$ respectively, then we may generalize (\ref{eq:cptpiece}) and (\ref{eq:cpt1}) as in Proposition~\ref{prop:strong K_X} and define the spaces
\[
\K_{\X, \Y}(N, M) = \K_{\X}(M) \B(L^2N, L^2M) \K_{\Y}(N);
\]
\[
\K_{\X, \Y}^{\infty, 1}(N, M) = \overline{ \X \, \B(L^2N, L^2M) \Y }^{\| \cdot \|_{\infty, 1}} 
= \overline{ \K_\X^L(M)^* \B(L^2N, L^2M) \K_\Y^L(N)}^{\| \cdot \|_{\infty, 1}}.
\]
When $\X$ and $\Y$ are the spaces of compact operators on $L^2M$ and $L^2N$ respectively we denote the space $\K_{\X, \Y}^{\infty, 1}(N, M)$ by $\K^{\infty, 1}(N, M)$. 

Realizing the extended Haagerup tensor product as a space of normal completely bounded operators $M \otimes_{eh} N \cong CB_{M'-N'}^\sigma(\B(L^2N, L^2M))$ we may restrict to $\X \, \B(L^2N, L^2M) \Y$ to obtain an isometric embedding into a dual operator space 
\begin{align*}
M \otimes_{eh} N
&\subset  CB_{M'-N'}(\X \, \B(L^2N, L^2M) \Y, \B(L^2N, L^2M)) \\
&\subset CB(\X \, \B(L^2N, L^2M) \Y, \B(L^2N, L^2M)) \\
&= ( ( \X \, \B(L^2N, L^2M) \Y ) \otimes_\pi \B(L^2N, L^2M)_* )^*
\end{align*}
where $ \otimes_\pi$ denotes the operator space projective tensor product. We introduce the $\X$-$\Y$-topology on $M \otimes_{eh} N$ as the restriction to $M \otimes_{eh} N$ of the weak$^*$-topology in $CB(\X \, \B(L^2N, L^2M) \Y, \B(L^2N, L^2M))$, so that a uniformly bounded net $\{ \mu_i \}_i \subset M \otimes_{eh} N$ will converge to $\mu \in M \otimes_{eh} N$ in the $\X$-$\Y$-topology if and only if we have ultraweak convergence $\mu_i(T) \to \mu(T)$ for each $T \in \X \, \B(L^2N, L^2M) \Y$. In the case when $M = N$ and $\X = \Y$ we refer to the $\X$-topology on $M \otimes_{eh} M$. 

Note that by the Krein-\u{S}mulian Theorem if $A \subset M \otimes_{eh} N$ is convex, then $A$ is closed in the $\X$-$\Y$-topology if and only if the intersection of $A$ with each closed ball in $M \otimes_{eh} N$ is closed in the $\X$-$\Y$-topology. As a consequence it follows that the space of bounded linear functionals in $(M \otimes_{eh} N)^*$ that are continuous with respect to the $\X$-$\Y$-topology forms a norm closed set. From this one then deduces from Lemma 1.2(iv) in \cite{StZs79} that a bounded linear functional $\varphi \in (M \otimes_{eh} N)^*$ is continuous in the $\X$-$\Y$-topology if and only if $\varphi$ is continuous in the $\X$-$\Y$-topology when restricted to any bounded set.

If we consider the isometric embedding
\[
M \otimes_{Con} N^{\rm op} \ni \mu \mapsto \mu \oplus {\mathcal J \mu^* \mathcal J} \in (M \otimes_{eh} N) \oplus^\infty ( (JMJ) \otimes_{eh} (JNJ)),
\]
then endowing $M \otimes_{eh} N$ and and $(J M J) \otimes_{eh} (J N J)$ with their $\X$-$\Y$-topologies we may restrict to obtain a topology on $M \otimes_{eh} N$, which we will, by abuse of terminology, also call the $\X$-$\Y$-topology on $M \otimes_{Con} N^{\rm op}$. 

The same analysis applies for $M \otimes_{Con} N^{\rm op}$ as for $M \otimes_{eh} N$, so that, in particular, a bounded linear functional $\varphi \in (M \otimes_{Con} N^{\rm op})^*$ is continuous in the $\X$-$\Y$-topology if and only if $\varphi$ is continuous in the $\X$-$\Y$-topology when restricted to any bounded set. 

We let $P_{\hat 1} \in \B(L^2N, L^2M)$ denote the rank-one partial isometry given by $P_{\hat 1}(\xi) = \langle \xi, \hat{1} \rangle \hat{1}$. If $\mu \in M \otimes_{Con} N^{\rm op}$, $a \in M$ and $b \in N$, then 
\begin{align*}
\langle \mu^*(J a P_{\hat{1}} b J) \hat{1}, \hat{1} \rangle 
&= {\rm Tr}(\mu^*( Ja J P_{\hat{1}} JbJ ) P_{\hat{1}} ) 
= {\rm Tr}(\mu^* (P_{\hat{1}}) JbJ P_{\hat{1}} Ja J) \\
&= \overline{ {\rm Tr}( \mu( aP_{\hat{1}} b)  P_{\hat{1}}  )} 
= \langle \hat{1}, \mu(a P_{\hat 1} b) \hat{1} \rangle.
\end{align*}

Taking spans and using weak$^*$-density it then follows that for all $T \in \B(L^2N, L^2M)$ we have 
\begin{equation}\label{eq:conj}
\langle \mu^*(JTJ) \hat{1}, \hat{1} \rangle = \langle \hat{1}, \mu(T) \hat{1} \rangle,
\end{equation}

If we now have a uniformly bounded net $\{ \mu_i \}_i \subset M \otimes_{Con} N^{\rm op}$, then $\mu_i \to 0$ in the $\X$-$\Y$-topology (in $M \otimes_{Con} N^{\rm op}$) if and only if for all $T \in \X \B(L^2N, L^2M)\Y$ we have $\mu_i(T) \to 0$ ultraweakly, and $\mu_i^*(J T J) \to 0$ ultraweakly. Note however, that since each $\mu_i$ is $M'$-$N'$ bimodular we have that $\mu_i(T) \to 0$ ultraweakly if and only if $\langle \mu_i(T) \hat{1}, \hat{1} \rangle \to 0$. However, by (\ref{eq:conj}) we have that $\mu_i^*( J T J) \to 0$ ultraweakly if and only if $\langle \hat{1}, \mu_i(T) \hat{1} \rangle = \langle \mu_i^*(JTJ) \hat{1}, \hat{1} \rangle \to 0$, and hence this already occurs when $\mu_i(T) \to 0$ ultraweakly. Thus, we see that $\mu_i \to 0$ in the $\X$-$\Y$-topology in $M \otimes_{Con} N^{\rm op}$ if and only if $\mu_i \to 0$ in the $\X$-$\Y$-topology in $M \otimes_{eh} N$, so that, at least on bounded sets, the $\X$-$\Y$-topology on $M \otimes_{eh} N$ restricts to the $\X$-$\Y$ topology on $M \otimes_{Con} N^{\rm op}$.

\subsection{Relatively mixing bimodules}
Given a dual normal operator $M$-$N$-bimodule $\mathcal W \subset \B(\mathcal H)$ we say that a vector $w \in \mathcal W$ is mixing relative to $\X \times \Y$ (or just mixing if $\X = \K(L^2M)$ and $\Y = \K(L^2N)$) if the map $M \otimes_{Con} N^{\rm op} \ni \mu \mapsto \mu(w) \in \mathcal W$ is continuous from the $\X$-$\Y$-topology to the weak$^*$-topology. We note that this is equivalent to the map $M \otimes_{Con} N^{\rm op} \ni \mu \mapsto \mu(w)$ being continuous on uniformly bounded subsets. We let $\mathcal W_{ {\X\text-\Y}{\rm -mix}}$ denote the set of vectors that are mixing relative to $\X \times \Y$ (we denote this space by $\mathcal W_{\rm mix}$ in the case when $\X = \K(L^2M)$ and $\Y = \K(L^2N)$). It is easy to see that this is a norm closed subspace of $\mathcal W$. We say that $\mathcal W$ is mixing relative to $\X \times \Y$ if $\mathcal W_{ {\X\text-\Y}{\rm -mix}} = \mathcal W$.

The bimodule $\mathcal W_{\rm mix}$ need not be weak$^*$-closed in general. However, it will always be closed in the $M$-$N$-topology, so that it is a strong operator $M$-$N$-bimodule. 

\begin{prop}\label{prop:mixstrong}
Let $M$ and $N$ be a finite von Neumann algebras, and $\X \subset \B(L^2M)$ and $\Y \subset \B(L^2N)$ boundary pieces for $M$ and $N$ respectively. Suppose $\mathcal W\subset \B(\cH)$ is a dual normal operator $M$-$N$-bimodule, then $\mathcal W_{{\X\text-\Y}{\rm -mix}}$ is a strong operator $M$-$N$-bimodule. 
\end{prop}
\begin{proof}
First note that since $\X \B(L^2N, L^2M) \Y$ is an $(M(\X) \cap M)$-$(M(\Y) \cap N)$ bimodule it follows that $\mathcal W_{{\X\text-\Y}{\rm -mix}}$ is also an $(M(\X) \cap M)$-$(M(\Y) \cap N)$ bimodule. 

We now fix $x \in M$ and $v \in \mathcal W_{{\X\text-\Y}{\rm -mix}}$, then by Kaplansky's Density Theorem there exists a uniformly bounded net $\{x_i \}_i \subset M(\X) \cap M$ such that $\| x - x_i \|_2 \to 0$. 

If $\{ \mu_j \}_j \subset M \otimes_{Con} N^{\rm op}$ is uniformly bounded such that $\mu_j \to 0$ in the $\X$-$\Y$-topology, then by Lemma~\ref{lem:sotconvergence}, for each $\varphi \in \mathcal W_*$ we have
\[
\lim_{i \to \infty} \sup_{j} \abs{ \varphi( \mu_j( (x - x_i) v ) } = 0.
\]
It then follows that $\lim_{i \to \infty} \varphi( \mu_j( x v ) ) = 0$, and since the net $\{ \mu_j \}_j \subset M \otimes_{Con} N^{\rm op}$ was arbitrary we then have that $xv \in \mathcal W_{{\X\text-\Y}{\rm -mix}}$. Thus $\mathcal W_{{\X\text-\Y}{\rm -mix}}$ is a left $M$-module, and a similar argument shows that it is also a right $N$-module.

We now fix $v \in \mathcal W$ and suppose $\{ p_i \}_i \subset M$ is an increasing sequence of projections which converge to 1 and such that $p_i v \in \mathcal W_{{\X\text-\Y}{\rm -mix}}$ for all $i$. Let $\mu_j \in M \otimes_{Con} N^{\rm op}$ be uniformly bounded such that $\mu_j \to 0$ in the $\X$-$\Y$-topology. Another application of Lemma~\ref{lem:sotconvergence} shows that for each $\varphi \in \mathcal W_*$ we have 
\[
\lim_{i \to \infty} \sup_{j} \abs{ \varphi( \mu_j (p_i^\perp v)) } = 0,
\] 
and it then follows just as above that $\lim_{j \to \infty} \varphi( \mu_j(v) ) = 0$, so that $v \in \mathcal W_{{\X}{\rm -mix}}$. By Theorem 2.1 in \cite{Ma98} it then follows that $\mathcal W_{{\X\text-\Y}{\rm -mix}}$ is a strong left $M$-module. The fact that $\mathcal W_{{\X\text-\Y}{\rm -mix}}$ is a strong right $N$-module follows similarly. 
\end{proof}

\begin{cor}
Under the hypotheses of the previous lemma, multiplication in $M \otimes_{Con} N^{\rm op}$ is separately continuous in the $\X$-$\Y$-topology.
\end{cor}
\begin{proof}
If $\mu \in M \otimes_{Con} N^{\rm op}$ and $\{ \mu_i \}_i \subset M \otimes_{Con} N^{\rm op}$ is uniformly bounded such that $\mu_i \to 0$ in the $\X$-$\Y$-topology, then as $\B(L^2N, L^2M) \ni S \mapsto \mu(S)$ is a normal map it then follows that $\mu( \mu_i(T) ) \to 0$ ultraweakly for all $T \in \X \B(L^2N, L^2M) \Y$, and hence, $\mu \cdot \mu_i \to 0$ in the $\X$-$\Y$-topology. 

We also have that $\mu_i(T) \to 0$ ultraweakly for each $T \in \B(L^2N, L^2M)_{{\X\text-\Y}{\rm -mix}}$, and by the previous proposition $\B(L^2M)_{{\X}{\rm -mix}}$ is an $M \otimes_{Con} M^{\rm op}$-module, which then shows that $\mu_i ( \mu( T) )  \to 0$ ultraweakly for all $T \in \X \B(L^2N, L^2M) \Y$, hence $\mu_i \cdot \mu \to 0$ in the $\X$-$\Y$-topology.
\end{proof}

We recall that $\B(L^2N, L^2M)^{M \sharp N}$ denotes the space of bounded linear functionals $\varphi \in \B(L^2N, L^2M)^*$ such that for each $T \in \B(L^2N, L^2M)$ the map $M \times M \ni (a, b) \mapsto \varphi(aTb)$ is separately normal in each variable. We let 
\[
\B(L^2N, L^2M)^{\sharp}_J = \B(L^2N, L^2M)^{M \sharp N} \cap \B(L^2N, L^2M)^{M' \sharp N'}.
\]

Note that by Theorem 4.2 in \cite{Ma98} if $\mathcal A$ is a unital $C^*$-algebra with $M \subset \mathcal A$, then we have $\varphi \in {\mathcal A}^{\sharp}$ if and only if the maps $L_\varphi: M \to {\mathcal A}^*$, $R_\varphi: M \to {\mathcal A}^*$ defined by $L_\varphi(x) (T) = \varphi(x^*T)$, and $R_\varphi(x)(T) = \varphi(Tx)$ are continuous from the ultrastrong topology on $M$ to the norm topology on ${\mathcal A}^*$. In particular, it follows that if $\varphi \in {\mathcal A}^{\sharp}$, then for each $T \in \mathcal A$ the map $M^2 \ni (x, y) \mapsto \varphi(x^* T y)$ is jointly strong operator topology continuous on bounded sets. 

\begin{lem}\label{lem:jordan}
Suppose $\mathcal A$ is a unital normal $M$-$C^*$-algebra. If $\varphi \in {\mathcal A}^{\sharp}$ is Hermitian with Jordan decomposition $\varphi = \varphi_+ - \varphi_-$, then $\varphi_+, \varphi_- \in  {\mathcal A}^{\sharp}$. 
\end{lem}
\begin{proof}
Suppose first that $\varphi \in {\mathcal A}^{\sharp}$. Take $x_i \in M_+$ increasing with $x_i \to 1$ strongly. For all $T \in \mathcal A$ we then have $\varphi(x_i^{1/2} T x_i^{1/2}) \to \varphi(T)$. Take arbitrary subnets such that $\psi_+$, and $\psi_-$ are weak$^*$ limits of $\mathcal A \ni T \mapsto \varphi_+(x_i^{1/2} T x_i^{1/2})$ and  $\mathcal A \ni T \mapsto \varphi_-(x_i^{1/2} T x_i^{1/2})$ respectively. Then $\varphi = \psi_+ - \psi_-$ and $\| \psi_+ \| + \| \psi_- \| \leq \| \varphi_+ \| + \| \varphi_- \| = \| \varphi \|$ so that we have $\psi_+ = \varphi_+$ and $\psi_- = \varphi_-$ by uniqueness of Jordan decomposition. As the subnets were arbitrary we then have $\varphi_+(x_i) \to \varphi_+(1)$, and $\varphi_-(x_i) \to \varphi_-(1)$. It therefore follows that ${\varphi_+}_{|M}, {\varphi_-}_{|M} \in M_*$ and hence $\varphi_+, \varphi_- \in {\mathcal A}^{\sharp}$.
\end{proof}

Recall that we may view $\K(L^2(M, \tau))^\natural$ as a subspace of $\B(L^2(M, \tau))^\natural$. Indeed, by the discussion in Section \ref{subs:con}, we have natural operator space isomorphisms 
\[
\K(L^2(M, \tau))^\natural\cong CB_{M-M}(\K(L^2(M, \tau)), \B(L^2(M, \tau)) ) \cong CB^\sigma_{M-M}(\B(L^2(M, \tau)))
\] 
and $\B(L^2(M, \tau))^\natural\cong CB_{M-M}(\B(L^2(M, \tau)))$. In the saem way, we may view $\K(L^2(M, \tau))^\natural_J$ as a subspace of $\B(L^2(M, \tau))^\natural_J$.

\begin{lem}\label{lem:approxbd}
Let $(M, \tau)$ be a tracial von Neumann algebra. Given a state $\varphi \in \B(L^2(M, \tau))^\natural$ (resp.\ $\B(L^2(M, \tau))^\natural_J$) there exists a net of states $\mu_i \in \K(L^2(M, \tau))^\natural$ (resp.\ $\K(L^2(M, \tau))^\natural_J$) such that $\mu_i \to \varphi$ weak$^*$ and $\| \mu_i \|_\natural \leq 2 \| \varphi \|_\natural$ (resp.\ $\| \mu_i \|_{Con} \leq 4\| \varphi \|_{Con}$). 
\end{lem}
\begin{proof}
Given a state $\varphi \in \B(L^2(M, \tau))^\natural$, we may find a net of positive functionals $\mu_i^0 \in\K(L^2M)^*$ converging weak$^*$ to $\varphi$, with $\|\mu_i^0\|\leq \|\varphi\|$. Since $\varphi_{| M}$ is normal we have that ${\mu_i}_{| M} \to \varphi_{|M}$ in the weak topology in $M_*$ and hence by taking convex combinations we may assume that $\| {\mu_i}_{| M} - \varphi_{|M} \| \to 0$. If $a_i \in M_+$ is such that $\mu_i(x) = \tau(x a_i)$ for $x \in M$, then setting $p_i = 1_{[0, 2\| \varphi \|_\natural]}(a_i)$ and $\tilde \mu_i( \cdot) = \mu_i(p_i\cdot p_i)/\tau(p_i a_i)$ we have that $\| \tilde \mu_i \|_\natural \leq 2 \| \varphi \|_\natural$. Moreover, since $\| {\mu_i}_{| M} - \varphi_{|M} \| \to 0$ we have $\mu_i(p_i) \to 1$ so that the Cauchy-Schwarz Inequality gives $\| \mu_i - \tilde \mu_i \| \to 0$. Hence, $\tilde \mu_i \to \varphi$ in the weak$^*$-topology. 

In the case when $\varphi \in \B(L^2(M, \tau))_J^\natural$ we may also assume that $\| {\mu_i}_{| JMJ} - \varphi_{| JMJ} \| \to 0$ and repeat the above argument by considering $\tilde{ \mu_i}_{| J M J}$ to construct $\tilde {\tilde {\mu}}_i \in \K(L^2(M, \tau))_J^\natural$ with $\| {\tilde {\tilde \mu }}_i \|_{Con} \leq 4\| \varphi \|_{Con}$. 
\end{proof}

\begin{lem}\label{lem:approxlf}
Let $M$ and $N$ be tracial von Neumann algebras. If $\varphi \in \B(L^2N, L^2M)^{\sharp}_J$, then there is a uniformly bounded net $\{ \mu_i \}_i \subset \K(L^2N, L^2M)^{\natural}_J$ such that $\varphi$ is the weak$^*$ limit of $\{ \mu_i \}_i$. 
\end{lem}
\begin{proof}
We first consider the case when $M = N$. By considering the real and imaginary parts separately, it suffices to consider the case when $\varphi$ is Hermitian. Moreover, by Lemma~\ref{lem:jordan} it then suffices to consider the case when $\varphi$ is positive. So suppose $\varphi \in \B(L^2M)^{\sharp}_J$ is a positive linear functional. For each $\varepsilon > 0$ there exist $K > 0$, and a projections $p \in \mathcal P(M)$,  with $\tau(p), \geq 1 - \varepsilon/2$ such that $\abs{ \varphi(p x p) } \leq K \| x \|_1$ for all $x \in M$.

We may similarly find a projection $q \in \mathcal P(JMJ)$ and $L > 0$ such that $\tau(q) \geq 1 - \varepsilon$ and $\abs{ \varphi(p q y q p) } \leq L \| y \|_1$ for all $y \in JMJ$. Note that if $x \in M$ has polar decomposition $x = v \abs{ x }$, then by Cauchy-Schwarz we have 
\begin{align*}
\abs{ \varphi( q p x p q ) }
&\leq \varphi(  p v \abs{ x }^{1/2} q \abs{ x }^{1/2} v^* p )^{1/2} \varphi(  p \abs{ x }^{1/2} q \abs{ x }^{1/2} p )^{1/2} \\
&\leq K \| v | x | v^* \|_1^{1/2} \| | x | \|_1^{1/2}
= K \| x \|_1. 
\end{align*}

Also, note that since $\varphi \in \B(L^2M)^{ \sharp}_J$ it follows that taking $\varepsilon$ tending to $0$, the corresponding net $q p \varphi p q$ converges weak$^*$ to $\varphi$. Since $q p \varphi p q \in \B(L^2M)^{\natural}_J$ the result then follows from Lemma~\ref{lem:approxbd}.

To prove the general case we set $\tilde M = M \oplus N$ and realize $\K(L^2N, L^2M) = (1_M \oplus 0) \K(L^2\tilde M) (0 \oplus 1_N)$, and $\B(L^2N, L^2M) = (1_M \oplus 0) \B(L^2\tilde M) (0 \oplus 1_N)$. If $\varphi \in \B(L^2N, L^2M)^{M \sharp N}_J$, then the map $\varphi'$ defined by $\varphi'(T) = \varphi( (1_M \oplus 0) T (0 \oplus 1_N) )$ defines an element of $\B(L^2\tilde M)^{ \sharp }_J$, and if $\mu_i \in \K(L^2\tilde M)^{ \sharp }_J$ are such that $\mu_i \to \varphi$ weak$^*$, then the maps ${\mu_i}_{| \K(L^2N, L^2M)}$ are in $\K(L^2N, L^2M)^{ \sharp }_J$ and converge weak$^*$ to $\varphi$. 
\end{proof}

\begin{thm}\label{thm:vanishing}
Let $M$ and $N$ be tracial von Neumann algebras. Suppose $\X \subset \B(L^2M)$ and $\Y \subset \B(L^2N)$ are boundary pieces. Fix $T \in \B(L^2 N, L^2M)$. The following conditions are equivalent:
\begin{enumerate}
\item\label{mixK} $T \in \K^{\infty, 1}_{\X, \Y}(N, M)$. 
\item\label{mixA} $T \in \B(L^2 N, L^2M)_{{\X\text-\Y}{\rm -mix}}$.
\item\label{mixM} $T \in \overline{ \overline{\K_{\X} \B(L^2N, L^2M) \K_{\Y}}^{_{M-N}} }^{_{JMJ-JNJ}}$.
\item\label{mixB} $\varphi(T) = 0$ whenever $\varphi \in \B(L^2 N, L^2M)^{\sharp}_J$ such that $\varphi_{| \X \B(L^2N, L^2M) \Y} = 0$.
\end{enumerate}

Moreover, if $M = N$, $\X = \Y$ with $M, JMJ \subset M(\X)$, and if $T \in M(\X)$, then the above conditions are also equivalent to either of the following conditions:
\begin{enumerate}\setcounter{enumi}{4}
\item\label{mixE} $\varphi(T) = 0$ for all states $\varphi \in \B(L^2 M)^{ \sharp }_J$ such that $\varphi_{|\X} = 0$.
\item\label{mixF} Whenever we have a net of states $\mu_i \in \K(L^2M)^\natural_J$ such that $\{ \mu_i \}_i$ is uniformly bounded in $M \otimes_{Con} M^{\rm op}$ and such that $\mu_i(S) \to 0$ for all $S \in \X$, then we have $\mu_i(T) \to 0$. 
\end{enumerate}
\end{thm}
\begin{proof}
The equivalence between (\ref{mixM}) and (\ref{mixK}) follows from Proposition~\ref{prop:infty 1 norm}. An application of the Hahn-Banach theorem gives the equivalence between (\ref{mixM}) and (\ref{mixB}). By Proposition~\ref{prop:mixstrong} we have that $\B(L^2 N, L^2M)_{{\X\text-\Y}{\rm -mix}}$ is a strong $M$-$N$ bimodule, and it is easy to see that it is also a strong $M'$-$N'$ bimodule. Since $\B(L^2 N, L^2M)_{{\X\text-\Y}{\rm -mix}}$ contains $\X \B(L^2N, L^2M) \Y$ it then follows from Corollary~\ref{cor:M-topology and norms} that $\K_{\X, \Y}^{\infty, 1}(N, M) = \overline{ \X \, \B(L^2N, L^2M) \Y }^{\| \cdot \|_{\infty, 1}} \subset \B(L^2 N, L^2M)_{{\X\text-\Y}{\rm -mix}}$, which shows that (\ref{mixK}) $\implies$ (\ref{mixA}).

To see that (\ref{mixA}) $\implies$ (\ref{mixB}) we suppose that condition (\ref{mixA}) is satisfied and take $\varphi \in \B(L^2M)^{ \sharp }_J$ such that $\varphi_{| \X \B(L^2N, L^2M) \Y} = 0$. By Lemma~\ref{lem:approxlf} there exists a uniformly bounded net $\{ \mu_i \}_i \subset M \otimes_{Con} M^{\rm op}$ such that $\mu_i \to \varphi$ weak$^*$. Since $\varphi_{| \X \B(L^2N, L^2M) \Y} = 0$ we have that $\mu_i(S) \to 0$ for all $S \in \X \B(L^2N, L^2M) \Y$ and hence $\varphi(T) = \lim_{i \to \infty} \mu_i(T) = 0$.

We now suppose that $M = N$, $\X = \Y$ with $M, JMJ \subset M(\X)$, and $T \in M(\X)$. Clearly we have (\ref{mixB}) $\implies$ (\ref{mixE}). Suppose $T \in M(\X)$ satisfies (\ref{mixE}) and let $\varphi \in \B(L^2M)_J^{ \sharp }$ be such that $\varphi_{| \X} = 0$. Taking the real and imaginary part of $\varphi$ separately we may assume that $\varphi$ is Hermitian. We then consider the Jordan decomposition $\varphi_{| M(\X )} = \varphi_+ - \varphi_-$. By Lemma~\ref{lem:jordan} we have that $\varphi_{+}$, and $\varphi_-$ each restrict to a normal state on both $M$ and $JMJ$. Since $\X$ is an ideal in $M(\X)$ it follows that $(\varphi_\pm)_{| \X} = 0$. If we extend $\varphi_+$ and $\varphi_-$ to arbitrary positive linear functionals $\psi_{+}$ and $\psi_-$ on $\B(L^2M)$ respectively that satisfy $\| \psi_+ \| = \| \varphi_+ \|$ and $\| \psi_- \| = \| \varphi_- \|$, then we have $\psi_\pm \in \B(L^2M)^{\sharp}_J$ and hence by (\ref{mixE}) we have $\psi_\pm(T) = 0$ showing that $\varphi(T) = \varphi_+(T) - \varphi_-(T) = \psi_+(T) - \psi_-(T) = 0$. 

Finally, the equivalence between (\ref{mixE}) and (\ref{mixF}) follows from Lemma~\ref{lem:approxbd}.
\end{proof}

\begin{rem}\label{commutative operations}
Given a normal operator $M$-$N$ and $JMJ$-$JNJ$ bimodule $X$, we note that the $M$-$N$ closure operation for $X$
	and $JMJ$-$JNJ$ closure operation for $X$ commute.
Indeed, this follows from Proposition \ref{prop:strongclosure} (\ref{en:strong2}).             
\end{rem}

\begin{rem}\label{rem:mix}
From the previous theorem it follows that $\B(L^2 N, L^2M)_{{\X\text-\Y}{\rm -mix}}$ defines the same subspace whether we view $\B(L^2 N, L^2M)$ as an $M$-$N$ bimodule, or as an $JMJ$-$JNJ$ bimodule. 
\end{rem}

Analogous to Proposition~\ref{prop:positivecpt}, we also have the following version for pairs of boundary pieces.

\begin{prop}\label{prop:positivecptbi}
Let $M$ and $N$ be tracial von Neumann algebras. Suppose $\X \subset \B(L^2M)$ and $\Y \subset \B(L^2N)$ are boundary pieces. If $T \in \B(L^2N, L^2M)$, then the following are equivalent:
\begin{enumerate}
\item\label{item:posbi1} $T \in \K_{\X, \Y}(N, M)$.
\item\label{item:posbi2} $T^*T \in \K_{\Y}(N)$ and $T T^* \in \K_{\X}(M)$.
\item\label{item:posbi3} $T^*T \in \K_{\Y}^{\infty, 1}(N)$ and $T T^* \in \K_{\X}^{\infty, 1}(M)$. 
\end{enumerate} 
\end{prop}
\begin{proof}
The implications (\ref{item:posbi1})$\implies$(\ref{item:posbi2})$\implies$(\ref{item:posbi3}) are clear. If (\ref{item:posbi3}) holds, then by Proposition~\ref{prop:positivecpt} we have $T^*T \in \K_{\Y}(N)$ and $T T^* \in \K_{\X}(M)$. Considering the polar decomposition $T = V | T |$ we then have $T = | T^* |^{1/2} V | T |^{1/2} \in \K_{\X, \Y}(N, M)$.
\end{proof}

\subsection{Mixing Hilbert bimodules}\label{sec:mixhilbert}

Let $M$ and $N$ be tracial von Neumann algebras and let $\mathcal H$ be an $M$-$N$ correspondence \cite{Co80, Po86}, i.e., $\mathcal H$ is a Hilbert space equipped with commuting normal representation of $M$ and $N^{\rm op}$, which we write as $\mathcal H \ni \xi \mapsto x \xi y$ for $x \in M$ and $y \in N$. We recall here some facts regarding correspondences that we will need in this section. We refer the reader to \cite[Chapter 13]{AnPo18} for proofs and further background. 

A vector $\xi \in \mathcal H$ is left-bounded if the map $L_\xi: N \to \mathcal H$ given by $L_\xi(x) = \xi x$ is bounded from $\| \cdot \|_2$, and hence defines a right $N$-modular operator in $\B_N(L^2N, \mathcal H)$. Moreover, if $L \in \B_N(L^2N, \mathcal H)$, then it is easy to check that $\xi = L (\hat{1}) \in \mathcal H$ is left-bounded and $L = L_\xi$. 

A vector $\eta \in \mathcal H$ is similarly defined to be right-bounded if $M \ni x \mapsto x\eta$ is bounded in $\| \cdot \|_2$, and hence defines a left $M$-modular bounded operator $R_\xi \in \B(L^2M, \mathcal H)$. A vector $\xi \in \mathcal H$ is bounded if it is both left and right-bounded. The space of bounded vectors is dense in $\mathcal H$ \cite{Po86}. 

If $\xi \in \mathcal H$ is left-bounded, then we obtain a u.c.p.\ map $\theta_\xi: M \to \B(L^2N)$ by $\theta_\xi(x) = L_\xi^* x L_\xi$. Moreover, since $L_\xi$ is right $N$-modular it follows that $\theta_\xi(x) \in JNJ' = N$, for each $y \in N$. If $\xi$ is a bounded vector, then we have $M \ni x \mapsto \theta_\xi(x) \in N$ extends to a bounded operator $T_\xi$ in $\B(L^2M, L^2N)$, and we may compute this operator explicitly as $T_\xi = L_\xi^* R_\xi$. 

The contragredient correspondence is given by the conjugate Hilbert space $\overline{\mathcal H}$ equipped with the $N$-$M$ bimodule structure
\[
y \cdot \overline \xi \cdot x = \overline{x^* \xi y^*},
\]
for $x \in M, y \in N$ and $\xi \in \mathcal H$. In this case the bounded vectors for $\overline {\mathcal H}$ are vectors of the form $\overline \xi$ where $\xi \in \mathcal H$ is a bounded vector, and we have $T_{\overline{\xi}} = T_\xi^*$. 

Suppose now that $\mathcal H$ and $\mathcal K$ are $M$-$N$ and $N$-$P$ correspondences respectively, and let $\mathcal H^L$ denote the space of left-bounded vectors. If $\xi_1, \xi_2 \in \mathcal H^L$, then $L_{\xi_2}^* L_{\xi_1} \in JNJ' \cap \mathbb B(L^2N) = N$, and so we have a well-defined sesquilinear form on $\mathcal H^L \otimes_{\rm alg} \mathcal K$ satisfying  
\[
\langle \xi_1 \otimes \eta_1, \xi_2 \otimes \eta_2 \rangle 
= \langle ( L_{\xi_2}^* L_{\xi_1} ) \eta_1, \eta_2 \rangle
\]
for $\xi_1, \xi_2 \in \mathcal H^L$, and $\eta_1, \eta_2 \in \mathcal K$. The separation and completion then gives a Hilbert space $\mathcal H \oovt{N} \mathcal K$, which is the Connes fusion of $\mathcal H$ and $\mathcal K$, and is naturally a $M$-$P$ correspondence satisfying $x (\xi \otimes \eta) y = (x \xi) \otimes (\eta y)$ for $x \in M$, $y \in P$, $\xi \in \mathcal H^L$, and $\eta \in \mathcal K$. If $\xi \in \mathcal H$ and $\eta \in \mathcal K$ are bounded, then the elementary tensor $\xi \otimes \eta \in \mathcal H \oovt{N} \mathcal K$ is also bounded and we have $T_{\xi \otimes \eta} = T_\eta T_\xi$.

In the case $\X = \K(L^2M)$ and $\Y = \K(L^2M)$ both of the conditions (\ref{item:mix1}) and (\ref{item:mix4}) below have been suggested as appropriate notions of a mixing bimodule \cite{BaFa11, PeSi12, OkOzTo17}. We show here that, in fact, these conditions are equivalent. 

\begin{thm}
Let $M$ and $N$ be tracial von Neumann algebras and suppose that $\X \subset \B(L^2M)$ and $\Y \subset \B(L^2N)$ are boundary pieces for $M$ and $N$ respectively such that $M, JMJ \subset M(\X)$ and $N, JNJ \subset M(\Y)$. Let $\mathcal H$ be an $M$-$N$ correspondence. The following conditions are equivalent:
\begin{enumerate}
\item\label{item:mix1} The set of bounded vectors $\xi \in \mathcal H$ such that $T_\xi \in \Y \B(L^2M, L^2N) \X$ is $M$-$N$ cyclic for $\mathcal H$. 
\item\label{item:mix2} The dual operator $M$-$N$ bimodule $\B_N(L^2N, \mathcal H)$ is mixing relative to $\X \times \Y$. 
\item\label{item:mix3} Every bounded vector $\xi \in \mathcal H$ satisfies $T_\xi \in \K_{\Y, \X}(M, N)$. 
\end{enumerate}
Moreover, in the case $\X = \K(L^2M)$ and $\Y = \K(L^2N)$ the above conditions are also equivalent to 
\begin{enumerate}\setcounter{enumi}{3}
\item\label{item:mix4} For every sequence $u_n \in \mathcal U(N)$ such that $u_n \to 0$ ultraweakly, and for all $\xi, \eta \in \mathcal H$ we have
\begin{equation}\label{eq:mix}
\lim_{n \to \infty} \sup_{x \in M, \| x \| \leq 1} | \langle x \xi u_n, \eta \rangle | = 0. 
\end{equation}
\end{enumerate}
\end{thm}
\begin{proof}
Given a bounded vector $\xi \in \mathcal H$, $x, c \in M$ and $y, a, b \in N$ we compute
\[
\langle (x T_\xi^* y) \widehat{ab}, \hat{c} \rangle
= \langle \xi yab , x^* c \xi \rangle
= \langle ( x L_\xi y ) \hat{a}, c \xi b^* \rangle.
\]
It then follows that for all $\mu \in M  \otimes_{Con} N^{\rm op}$ we have 
\[
\langle \mu(T_\xi^*) \widehat{ab}, \hat{c} \rangle 
=\langle \mu( L_\xi ) \hat{a}, c \xi b^* \rangle.
\]
Since the range of $\mu( L_\xi )$ is contained in the closure of the span of vectors of the form $c \xi b^*$ for $c \in M$ and $b \in N$ it then follows that $T_\xi^* \in \B(L^2 N, L^2M)_{{\X\text-\Y}{\rm -mix}}$ if and only if $L_\xi \in \B(L^2 N, \mathcal H)_{{\X\text-\Y}{\rm -mix}}$. The set $\mathcal H_0 := \{ \xi \in \mathcal H \mid L_\xi \in \B_N(L^2 N, \mathcal H)_{{\X\text-\Y}{\rm -mix}} \}$ forms a $M$-$N$ bimodule and so assuming that (\ref{item:mix1}) holds we conclude that this subspace is dense. 

Suppose now that $\xi_n \in \mathcal H_0$ are such that $\xi_n \to \xi \in \mathcal H$ with $\xi$ a left-bounded vector. Let $\varphi_n \in N_*$ be the positive linear functional defined by $\varphi_n(a) = \langle (\xi - \xi_n) a, (\xi - \xi_n) \rangle$. Then we have $\| \varphi_n \|_1 \to 0$ and so by passing to a subsequence we may produce an increasing sequence of projections $p_k \in \mathcal P(N)$ such that $\tau(p_k) \to 1$ and for each fixed $k$ we have 
\[
\|  L_{\xi p_k} - L_{\xi_n p_k}  \| = \sup_{a \in N, \| a \|_2 \leq 1} \| (\xi - \xi_n) p_k a \| = \sup_{a \in N, \| a \|_2 \leq 1} \varphi_n( p_k a^*a p_k ) \to 0.
\]
Since $\xi_n p_k \in \mathcal H_0$ for each $n, k \geq 0$ and since $\B_N(L^2 N, \mathcal H)_{{\X\text-\Y}{\rm -mix}}$ is a strong right $N$-module by Proposition~\ref{prop:mixstrong}, it then follows that $\xi p_k \in \mathcal H_0$ for each $k \geq 0$ and also $\xi \in \mathcal H_0$. Thus, 
\[
\B_N(L^2 N, \mathcal H)_{{\X\text-\Y}{\rm -mix}} = \{ L_\xi \mid \xi \in \mathcal H \ {\rm is \ bounded } \} = \B_N(L^2N, \mathcal H).
\]

We now suppose (\ref{item:mix2}) holds, then we see from above that the $T_\xi^* \in \B(L^2N, L^2M)_{{\X\text-\Y}{\rm -mix}}$ for each bounded vector $\xi$. Since $M, JMJ \subset \X$ and $N, JNJ \subset \Y$ we then conclude that for every bounded vector $\xi \in \mathcal H$ and $\varepsilon > 0$, there exist projections $p \in \mathcal P(M)$, and $q \in \mathcal P(N)$ so that $\tau(p), \tau(q) > 1-\varepsilon$ and $T_{p \xi q} = p (JpJ) T_\xi (JqJ) q \in \Y \B(L^2M, L^2N) \X$. Thus, the set of bounded vectors $\xi \in \mathcal H$ such that $T_{\xi} \in \Y \B(L^2M, L^2N) \X$ is dense in $\mathcal H$.

We have thus shown the equivalence between (\ref{item:mix1}) and (\ref{item:mix2}). We have also shown that, in fact, these conditions imply also that the set of bounded vectors $\xi \in \mathcal H$ such that $T_\xi \in \Y \B(L^2M, L^2N) \X$ is dense in $\mathcal H$. In particular, it then follows that there is a dense subset $\mathcal H_1 \subset \mathcal H$ such that $T_{\xi \otimes \overline \eta} = T_{\eta}^* T_\xi \in \mathbb X \B(L^2M) \mathbb X$ for all $\xi, \eta \in \mathcal H_1$. Since these elementary tensors span a dense subset of $\mathcal H \oovt{N} \overline{\mathcal H}$ we see that $\mathcal H \oovt{N} \overline{\mathcal H}$ satisfies the hypotheses of (\ref{item:mix1}), so that by the equivalence between (\ref{item:mix1}) and (\ref{item:mix2}) it follows that if $\xi \in \mathcal H$ is any bounded vector, then we have $T_\xi^* T_\xi = T_{\xi \otimes \overline \xi} \in \B(L^2M)_{{\X}{\rm -mix}} = \K_\X^{\infty, 1}(M)$. 

Considering instead $\overline {\mathcal H} \oovt{M} \mathcal H$ we see that we also have $T_\xi T_\xi^* \in \K_\Y^{\infty, 1}(N)$ and by Proposition~\ref{prop:positivecptbi} it then follows that $T_\xi \in \K_{\Y, \X}(M, N)$, establishing (\ref{item:mix3}). That we also have (\ref{item:mix3}) $\implies$ (\ref{item:mix2}) is clear.

We now suppose that $\X = \K(L^2M)$ and $\Y = \K(L^2N)$. Note that clearly (\ref{item:mix2}) $\implies$ (\ref{item:mix4}) since for any uniformly bounded sequence $x_n \in M$ we have $x_n \otimes (u_n^*)^{\rm op} \to 0$ weak$^*$ when viewed as an element in $M \otimes_{Con} N^{\rm op}$. 

We now suppose that (\ref{item:mix4}) holds. The result is trivially satisfied if $N$ is completely atomic, and so by restricting to the orthogonal projection of the completely atomic part we may assume that $N$ is diffuse. Note that Equation (\ref{eq:mix}) holds also if instead of $u_n \in \mathcal U(N)$ we consider $y_n \in N$ contractions such that $y_n \to 0$ ultraweakly. Indeed, by the polarization identity it is enough to check this for $\xi = \eta$, and if Equation (\ref{eq:mix}) does not hold, then there exist $c > 0$ and $x_n \in M$ contractions such that $| \langle x_n \xi y_n, \xi \rangle | \geq c$. By considering separately the real and imaginary parts of $x_n$ and $y_n$ we may assume that $x_n = x_n^*$ and $y_n = y_n^*$ so that $y \mapsto \langle x \xi y, \xi \rangle$ is a Hermitian linear functional. 

We let $v_n \in \{ y_n \}' \cap N$ be a sequence of self-adjoint unitaries such that $v_n (1 - y_n^2)^{1/2} \to 0$ ultraweakly (note that this is possible since $\{y_n \}' \cap N$ is diffuse for each $n$). Setting $u_n = y_n + iv_n(1 - y_n^2)^{1/2} \in \mathcal U(N)$ we then have that $u_n \to 0$ ultraweakly and $| \langle x_n \xi u_n, \xi \rangle | \geq  | \langle x_n \xi y_n, \xi \rangle | \geq c$, which would then contradict (\ref{item:mix4}). 

Thus, we have established that for each bounded vector $\xi \in \mathcal H$ and for each uniformly bounded sequence $y_n \in N$ we have 
\[
\lim_{n \to \infty} \sup_{x \in M, \| x \| \leq 1} | \langle x \xi y_n, \xi \rangle | = 0,
\]
i.e., the operator $T_\xi$ is compact as an operator from $N$ into $L^1M$. By Lemma~\ref{lem:cptl2} we then have that $T_\xi$ is compact as an operator from $N$ into $L^2M$, which by \cite{Oza10} agrees with the $\| \cdot \|_{\infty, 2}$-closure of $\K(L^2N, L^2M)$ in $\B(L^2N, L^2M)$. Cutting down by projections as above then shows that a dense set of bounded vectors $\xi \in \mathcal H$ have the property that $T_\xi \in \K(L^2N, L^2M)$ showing that (\ref{item:mix1}) holds. 
\end{proof}

\begin{examp}\label{examp:correbdry}
If $\mathcal H$ is an $M$-$N$ correspondence, then we let $\X_{\mathcal H} \subset \B(L^2M)$ denote the hereditary $C^*$-algebra generated by operators of the form $T_\xi^* T_\eta$ where $\xi, \eta \in \mathcal H$ are bounded vectors. We similarly let $\Y_{\mathcal H} \subset \B(L^2N)$ denote the hereditary $C^*$-algebra generated by operators of the form $T_\eta T_\xi^*$. Then $\X_{\mathcal H}$ and $\Y_{\mathcal H}$ give the smallest boundary pieces of $M$ and $N$ respectively so that $\mathcal H$ is mixing relative to $\X_{\mathcal H} \times \Y_{\mathcal H}$. 

Note that if $B \subset M$ is a von Neumann subalgebra and we consider $L^2M$ as an $M$-$B$ correspondence, then the corresponding boundary piece $\X_{L^2M}$ is the one described in Example~\ref{examp:bdrysub}.
\end{examp}

\section{Properly proximal von Neumann algebras}\label{sec:propprox}

Let $M$ be a von Neumann algebra, $\X \subset \B(L^2M)$ a boundary piece, and suppose that $\mathcal V$ is a dual normal $M$-bimodule. 
A point $\xi \in {\mathcal V}$ is properly proximal relative to $\X$ (or just properly proximal if $\X = \K(L^2M)$) if $[ x, \xi ] \in {\mathcal V}_{{\X}{\rm -mix}}$ for all $x \in M$. We denote the set of properly proximal points by ${\mathcal V}_{{\X}{\rm -prox}}$ (or $\mathcal V_{\rm prox}$ if $\X = \K(L^2M)$). Note that ${\mathcal V}_{{\X}{\rm -mix}} \subset {\mathcal V}_{{\X}{\rm -prox}}$, and both are closed subspaces of ${\mathcal V}$. While ${\mathcal V}_{{\X}{\rm -prox}}$ need not be an $M$-bimodule, it is $M$-convex in the sense that if $x_1, \ldots, x_n \in M$ satisfy $\sum_{i = 1}^n x_i^* x_i = 1$, then $\sum_{i =1}^n x_i^* \xi x_i \in {\mathcal V}_{{\X}{\rm -prox}}$ for all $\xi \in {\mathcal V}_{{\X}{\rm -prox}}$. Indeed, for all $y \in M$ we have 
\[
[y, \sum_{i = 1}^n x_i^* \xi x_i ] = \sum_{i = 1}^n y [x_i^*, \xi] x_i + [y, \xi] + \sum_{i = 1}^n x_i^* [ x_i, \xi] y \in {\mathcal V}_{{\X}{\rm -mix}}.
\]

In practice, explicitly verifying that $[x, \xi] \in {\mathcal V}_{{\X}{\rm -mix}}$ for each $x \in M$ may be cumbersome. The following lemma is a key tool in this respect. This should be compared with Lemma 6.7 in \cite{IsPeRu19}.

\begin{lem}
\label{lem:proxmix}
Let $M$ be a tracial von Neumann algebra, let $\mathcal V$ be a dual normal $M$-bimodule, and fix $\xi \in \mathcal V$. Then $\{ x \in M \mid [x, \xi], [x^*, \xi] \in \mathcal V_{{\X}{\rm -mix}} \}$ is a von Neumann subalgebra of $M$. In particular, $\xi \in \mathcal V_{{\X}{\rm -prox}}$ if and only if $\{ x \in M \mid [x^*, \xi], [x, \xi] \in {\mathcal V}_{{\X}{\rm -mix}} \}$ is ultraweakly dense in $M$.
\end{lem}
\begin{proof}
It is easy to check that $M_0 = \{ x \in M \mid [x^*, \xi], [x, \xi] \in {\mathcal V}_{{\X}{\rm -mix}} \}$ is a unital $*$-subalgebra, so we will only show that it is ultraweakly closed. If $x$ is in the ultraweak closure of $M_0$, then we may take a net $\{ x_i \}_i \subset M_0$ that converges ultraweakly to $x$. Since $[x_i, \xi] - [x, \xi] = (x_i - x) \xi - \xi(x_i - x)$ converges to $0$ in the $M$-$M$ topology, and since ${\mathcal V}_{{\X}{\rm -mix}}$ is a strong $M$-$M$ bimodule it then follows that $[x, \xi] \in {\mathcal V}_{{\X}{\rm -mix}}$. We similarly have $[x^*, \xi] \in {\mathcal V}_{{\X}{\rm -mix}}$, and so $x \in M_0$. 
\end{proof}

We let $\bS_{\X}(M) \subset \B(L^2M)$ be the set of operators that are properly proximal relative to $\X$ when we view $\B(L^2M)$ as an $M$-bimodule under the actions $x \cdot T \cdot y = Jy^*J T J x^* J$, i.e.,
\[
\bS_{\X}(M) = \{ T \in \B(L^2M) \mid [T, J x J] \in \K_{\X}^{\infty, 1} {\rm \ for \ all \ } x \in M \}.
\]
Note that $\bS_{\X}(M)$ is an operator system that contains $M$. In the case when $\X = \K(L^2M)$ we write $\bS(M)$ instead of $\bS_{\K(L^2M)}(M)$.

\begin{thm}\label{thm:propprox}
Let $M$ be a tracial von Neumann algebra, $\X \subset \B(L^2M)$ a boundary piece, and let $B \subset M$ be a von Neumann subalgebra. The following are equivalent:
\begin{enumerate}
\item\label{item:prox1} There exists a $B$-central state $\varphi$ on $\bS_{\X}(M)$ such that $\varphi_{| M }$ is normal.
\item\label{item:prox5} There exists a non-zero projection $p \in \mathcal Z(B)$ and a $B$-central state $\varphi$ on $\bS_\X(M)$ such that $\varphi_{|pMp} = \frac{1}{\tau(p)}\tau$. 
\item\label{item:prox4} There exists a non-zero projection $p \in \mathcal Z(B)$ and an $Mp$-bimodular u.c.p.\ map $\Phi: \bS_\X(M) \to \langle pMp, e_{Bp} \rangle$. 
\item\label{item:prox2}  If $E$ is any normal operator $M$-system such that there exists a state $\varphi_0 \in (E^\natural)_{{\X}{\rm -prox}}$ with ${\varphi_0}_{|M} = \tau$, then there exists an $B$-central state $\varphi$ on $E$ such that $\varphi_{|M}$ is normal. 
\end{enumerate}
\end{thm}
\begin{proof}
The equivalences between (\ref{item:prox1}), (\ref{item:prox5}), and (\ref{item:prox4}) are standard. 
To see that (\ref{item:prox2}) $\implies$ (\ref{item:prox1}) simply observe that the state $\varphi_0(T) = \langle T \hat 1, \hat 1 \rangle$ is in $(\bS_{\X}(M)^\natural)_{{\X}{\rm -prox}}$, which follows from the remark after Theorem~\ref{thm:vanishing}. Conversely, to see that (\ref{item:prox1}) $\implies$ (\ref{item:prox2}) note that if $E$ is a normal operator $M$-system and $\varphi_0 \in (E^\natural)_{{\X}{\rm -prox}}$ with ${\varphi_0}_{|M} = \tau$, then $\varphi_0$ corresponds to an $M$-bimodular u.c.p.\ $\Phi_0: E \to \B(L^2M)$ such that $\varphi_0(T) = \langle \Phi_0(T) \hat{1}, \hat{1} \rangle$, for $T \in E$. Since $\varphi_0$ is a properly proximal point it follows that the range of $\Phi_0$ is contained in $\bS_\X(M)$. Indeed, if we view $\B(L^2M)$ as an $M$-$M$ bimodule with the bimodule structure given by $x \cdot T \cdot y = Jy^*J T Jx^*J$, then for $T \in E$, $a, b, x \in M$, and if $\mu_i \in M \otimes_{Con} M^{\rm op}$ is uniformly bounded such that $\mu_i \to \mu$ in the $\mathbb X$-topology, then since $[x, \varphi_0] \in (E^\natural)_{{\X}{\rm -mix}}$ we see that
\[
\langle \mu_i(  x \cdot \phi_0(T) - \phi_0(T) \cdot x ) \hat{a}, \hat{b} \rangle
=  \mu_i( [x, \varphi_0]) (b^* T a )
\to \langle \mu(  x \cdot \phi_0(T) - \phi_0(T) \cdot x ) \hat{a}, \hat{b} \rangle.
\]
Since $a, b \in M$ are arbitrary it follows that $[\phi_0(T), Jx^*J] = x \cdot \phi_0(T) - \phi_0(T) \cdot x \in \B(L^2M)_{{\X}{\rm -mix}}$. By Theorem~\ref{thm:vanishing} and Remark~\ref{rem:mix} it then follows that $\phi_0(T) \in \mathbb S_{\X}(M)$ for each $T \in E$. 

If we have a $B$-central state $\varphi: \bS_\X(M) \to \mathbb C$ with $\varphi_{|M}$ normal, then $\varphi \circ \Phi_0$ gives an $B$-central state on $E$ with ${\varphi \circ \Phi_0}_{| M}$ normal. 
\end{proof}

We say that the inclusion $B \subset M$ is properly proximal relative to $\X$ (or just properly proximal if $\X = \K(L^2M)$) if it fails to satisfy the conditions in the previous theorem. We say that $M$ is properly proximal if the inclusion $M \subset M$ is properly proximal. We remark that by condition (\ref{item:prox4}) above proper proximality for $M$ is independent of the given trace. 

\begin{rem}
If $\X$ and $\Y$ are boundary pieces with $\X \subset \Y$, then we have $\K_\X^{\infty, 1}(M) \subset \K_\Y^{\infty, 1}(M)$ and hence $\bS_\X(M) \subset \bS_\Y(M)$. Thus, if $B \subset M$ is properly proximal relative to $\X$, then it is also properly proximal relative to $\Y$. In particular, if $B \subset M$ is properly proximal, then the inclusion is properly proximal with respect to any boundary piece. Similarly, if $B \subset Q \subset M$, and $B \subset M$ is properly proximal relative to $\X$, then $Q \subset M$ is also properly proximal relative to $\X$. 

Also note that if $B \subset Q \subset M$ and we identify $L^2Q = e_Q L^2M$, then $\overline{ e_Q \K_\X(M) e_Q } \subset \B(L^2Q)$ gives a boundary piece for $Q$. Indeed, $\overline{ e_Q \K_\X(M) e_Q }$ is a $C^*$-algebra, and if $0 \leq a \leq b \in \overline{ e_Q \K_\X(M) e_Q }$, then we may write $a = b^{1/4} c b^{1/4}$ for $c \in \mathbb B(L^2Q)$ \cite[II.3.2.5]{Bl06}, and taking an approximate identity $\{ p_i \}_i \subset \overline{ e_Q \K_\X(M) e_Q }$ we then have $a = \lim_{i \to \infty} p_i a p_i \in \overline{ e_Q \K_\X(M) e_Q }$, showing that $\overline{ e_Q \K_\X(M) e_Q }$ is hereditary. Moreover, since $\K_\X(M)$ contains both $M$ and $JMJ$ in its multiplier algebra we see that $\overline{ e_Q \K_\X(M) e_Q }$ contains both $Q$ and $JQJ$ in its multiplier algebra. 

We then see that $e_Q \K_\X^{\infty, 1}(M) e_Q \subset \K_{\overline{ e_Q \K_\X(M) e_Q }}^{\infty, 1}(Q)$ and it follows that $e_Q \bS_\X(M) e_Q \subset \bS_{\overline{ e_Q \K_\X(M) e_Q }}(Q)$. We therefore see that if $B \subset M$ is properly proximal relative to $\X$, then $B \subset Q$ is properly proximal relative to $\overline{ e_Q \K_\X(M) e_Q }$. In particular, if $B \subset M$ is properly proximal, then $B \subset Q$ is also properly proximal (and $B$ is also properly proximal). 

The same argument above shows that if $p \in \mathcal Z(M)$, then by identifying $L^2 (zM) = z L^2M$ we have that $\overline{p^\perp \mathbb X p^\perp}$ is a boundary piece for $zM$. 
\end{rem}

Recall from \cite{BoIoPe21} that if $\Gamma$ is a group, then a boundary piece $X$ for $\Gamma$ is a non-empty closed $\Gamma \times \Gamma$-invariant subspace of $\beta \Gamma \setminus \Gamma$. We have a bijective correspondence between boundary pieces for $\Gamma$ and $\Gamma \times \Gamma$-invariant proper closed ideals $c_0 \Gamma \subset I \subset \ell^\infty \Gamma$, which is given by  $I = \{ f \in \ell^\infty \Gamma = C(\beta\Gamma) \mid f_{| X} = 0 \} \subset \ell^\infty \Gamma$. A group $\Gamma$ is properly proximal relative to $X$ (or relative to $I$) if there is no left $\Gamma$-invariant state on $C(X)^{R(\Gamma)}=(\ell^\infty\Gamma/ I)^{R(\Gamma)}$. 

We set $\bS_X(\Gamma) = \bS_I(\Gamma) =\{f\in\ell^\infty\Gamma\mid f-R_tf\in I{\rm \ for\ any}\ t\in\Gamma\}$, which is a $C^*$-subalgebra of $\ell^\infty \Gamma$. The group $\Gamma$ acts on $\bS_X$ by left-translation, and $\Gamma$ is properly proximal relative to $X$ if and only if there is no $\Gamma$-invariant state on $\bS_X(\Gamma)$.

\begin{thm}\label{thm: iff}
Let $\Gamma$ be a group with a $\Gamma \times \Gamma$-invariant proper closed ideal $c_0 \Gamma \subset I \subset \ell^\infty \Gamma$, and let $\X$ be a boundary piece for $L\Gamma$. We let $I_\X \subset \ell^\infty \Gamma$ the closed ideal generated by $E_{\ell^\infty \Gamma}( \X)$, and we let $\X_I$ denote the $L\Gamma$-boundary piece $\overline{ I \B(\ell^2 \Gamma) I }^{\lVert \cdot \rVert}$ generated by $I$. The following statements are true:
\begin{enumerate}
\item If $L\Gamma$ is properly proximal relative to $\X$, then $\Gamma$ is properly proximal relative to $I_\X$. 
\item $\Gamma$ is properly proximal relative to $I$ if and only if $L\Gamma$ is properly proximal relative to $\X_I$. 
\end{enumerate}
\end{thm}
\begin{proof}
First assume that $\Gamma$ is not properly proximal relative to $I_\X$, i.e., there exists a $\Gamma$-left invariant state $\varphi$ on $\bS_{I_\X}(\Gamma)$. As $E_{\ell^\infty \Gamma}: \B(\ell^2 \Gamma) \to \ell^\infty \Gamma$ is continuous from the $\lVert \cdot \rVert_{\infty, 1}$-topology to the norm-topology it follows that $E_{\ell^\infty \Gamma}$ maps $\K_\X^{\infty, 1}(L\Gamma)$ into $I_\X$, and consequently, $E_{\ell^\infty \Gamma}$ maps $\bS_\X(L\Gamma)$ into $\bS_{I_\X}(\Gamma)$, since if $T \in \bS_\X(L\Gamma)$, then $ E_{\ell^\infty \Gamma}(T)-R_t (E_{\ell^\infty \Gamma}(T) )= E_{\ell^\infty \Gamma}(T-\rho_t T\rho_t^*)\in I_{\X}$, for each $t \in \Gamma$. 

We may therefore consider the state $\psi:=\varphi\circ E_{\ell^\infty \Gamma}: \bS_{\X}(L\Gamma) \to \mathbb C$.
Observe that for any $x\in L\Gamma$, $E_{\ell^\infty \Gamma}(x)=\tau(x)$ and hence $\psi_{\mid L\Gamma}=\tau$.
Moreover, note that $E_{\ell^\infty \Gamma}$ is left $\Gamma$-equivariant, i.e., $L_s (E_{\ell^\infty \Gamma}(T))=E_{\ell^\infty \Gamma}(\lambda_s T\lambda_s^*)$ for any $s\in\Gamma$ and $T\in \bS_{\X}(L\Gamma)$
and thus $\psi(\lambda_s T)= \psi( \lambda_s T \lambda_s \lambda_s^* ) = \psi(T\lambda_s)$.
Finally, using the fact $\psi$ is normal on $L\Gamma$ and the Cauchy-Schwarz inequality, we conclude that $\psi$ is $L\Gamma$-central.

Next, suppose $L\Gamma$ is not properly proximal relative to $\X_I$ and let $\psi$ be an $L\Gamma$-central state on $\bS_{\X_I}(L\Gamma)$.
We claim that $\bS_I(\Gamma) \hookrightarrow \bS_{\X_I}(L\Gamma)$ by viewing $f\in \bS_I(\Gamma)$ as a multiplier $M_f$.
Indeed, for any $f\in \bS_I(\Gamma)$ and $t\in\Gamma$, we have $[M_f, \rho_t]=M_{f-R_t(f)} \rho_t\in \X_I$, and hence $M_f \in \bS_{\X_I}(L\Gamma)$ by Lemma~\ref{lem:proxmix}. 
Since $\lambda_t M_f\lambda_t^*=M_{L_t (f)}$ for $t \in \Gamma$ it follows that $\psi$ gives a $\Gamma$-invariant state on $\bS_I(\Gamma)$. The ``only if'' direction of the second statement follows from the first statement upon noticing that $E_{\ell^\infty}(I \B(\ell^2 \Gamma) I) = I \ell^\infty \Gamma I = I$. 
\end{proof}

We recall that an von Neumann subalgebra $N \subset M$ is co-amenable if there exists a conditional expectation from $\langle M, e_N \rangle$ to $M$ (see, e.g., \cite{MoPo03}).

\begin{lem}\label{lem:co-amenable}
	Let $(M,\tau)$ be a finite von Neumann algebra with a faithful normal trace $\tau$,  $N\subset M$ a von Neumann subalgebra and $E$ an $M$-system. Suppose there exists an $N$ central state $\varphi$ on $E$ such that $\varphi_{\mid M} = \tau$.  If $N\subset M$ is co-amenable, then there exists an $M$-central state $\psi$ on $E$ with $\psi_{\mid M} = \tau$.
\end{lem}
\begin{proof}
	As in Section~\ref{subs:con}, there exists an $M$-bimodular u.c.p.\  map $\Phi: E\to \B(L^2M)$ such that $\varphi(T)=\langle \Phi(T)\hat 1,\hat 1\rangle$ for any $T\in E$. 
	Since $\varphi$ is $N$-central, it is clear that $\Phi: E\to \langle M, e_N\rangle$.
	Denote by $\mathcal E$ a conditional expectation $\langle M, e_N\rangle\to M$ given by co-amenability, then $\tau\circ \mathcal E\circ \Phi$ is an $M$-central state on $E$ that restricts to $\tau$ on $M$.
\end{proof}

\begin{prop}\label{prop:co-amenable}
	Let $(M,\tau)$ be a finite factor and $N\subset M$ a co-amenable subfactor. If $M$ is properly proximal, then so is $N$.
\end{prop}
\begin{proof}
	Suppose $N$ is not properly proximal, then there exists an $N$-central state $\varphi$ on $\bS(N)$ with $\varphi_{\mid N} = \tau$.
	Notice that $\Ad(e_N):\bS(M)\to\bS(N)$ and hence $\psi:=\Ad(e_N)\circ\varphi$ is an $N$-central state on $\bS(M)$ with $\psi_{\mid N} = \tau$.  Now it follows from Proposition~\ref{lem:co-amenable} that there exists an $M$-central state on $\bS(M)$ and thus $M$ is not properly proximal.
\end{proof}

\begin{prop}
Let $M_1$ and $M_2$ be diffuse tracial von Neumann algebras. Then $M= M_1* M_2$ is properly proximal.
\end{prop}
\begin{proof}
We will show that the usual paradoxical decomposition that proves nonamenability for $M_1 * M_2$ also works to show proper proximality. For $i=1,2$, let $\cH_i= L^2(M_i)$ and $\cH_i^0=\mathcal H_i\ominus \C \hat{1}=\overline{M_i^0\hat{1}}^{\|\cdot\|_{2}}$, where $M_i^0$ denotes the kernel of the trace.
Recall that $\cH = L^2(M_1 * M_2)$ decomposes as
 \[
 \cH= \C \hat{1} \oplus \bigoplus_{n\geq 1}\left( \bigoplus_{i_1\neq i_2\neq \cdots\neq i_{n}} \cH_{i_1}^0\otimes \cdots \otimes\cH_{i_n}^0 \right).
 \]
Set 
 \[
 \cH_\ell(1)=\C \hat{1} \oplus \bigoplus_{n\geq 1}\left( \bigoplus_{\substack{i_1\neq i_2\cdots \neq i_{n}\\
 i_1\neq 1}} \cH_{i_1}^0\otimes \cdots \otimes\cH_{i_n}^0 \right),
 \]
 and let $P\in \B(\cH)$ be the orthogonal projection onto $\cH_\ell(1)^\perp$.

If $z \in M_2^0$, then as $JzJ$ and $Jz^*J$ preserve the space $\mathcal H_\ell(1)$ we have $[P, JzJ] = 0$. Also, if $z \in M_1^0$, then as $JzJ$ and $Jz^*J$ preserve the subspace $\mathcal H_\ell(1) \ominus \mathbb C \hat{1}$ we have $[P, JzJ] = [JzJ, {\rm Proj}_{\mathbb C \hat{1}}]$ is finite-rank. By Lemma~\ref{lem:proxmix} we then have $P \in \bS(M)$.

Similarly, if we by $Q$ the projection onto $H_\ell(2)^\perp$, where 
\[
 \cH_\ell(2)=\C \hat{1} \oplus \bigoplus_{n\geq 1}\left( \bigoplus_{\substack{i_1\neq i_2\cdots \neq i_{n}\\
 i_1\neq 2}} \cH_{i_1}^0\otimes \cdots \otimes\cH_{i_n}^0 \right),
 \]
 then $Q\in \bS(M)$. Since $M_2$ is diffuse, we may choose orthogonal trace-zero unitaries $u_1, u_2\in \cU(M_2)$ so that we then have $u_1^*Pu_1 + u_2^*P u_2 \leq Q$.
Similarly we may choose a trace-zero unitary $v\in \cU(M_1)$ and obtain that $v^*Qv\leq P$.

If there were an $M$-central state $\varphi$ on $\bS(M)$, we would then have $2\varphi(P)\leq \varphi (Q)\leq \varphi(P)$ and hence $\varphi(P + Q) = 0$. Since $1 - P - Q$ is the projection onto $\mathbb C\hat{1}$ and $\varphi$ is $M$-central we have $\varphi(1-P-Q)=0$, which then gives a contradiction. 
\end{proof}

\begin{thm}\label{thm:ppcommute}
Let $M$ be a finite von Neumann algebra, $\X \subset \B(L^2M)$ an $M$-boundary piece and let $\mathcal G \subset \mathcal U(M)$ be a subgroup. Suppose there exists a state $\varphi \in \B(L^2M)^*$ such that 
\begin{enumerate}
\item $\varphi$ restricts to the canonical traces on $M$ and $JMJ$.
\item $\varphi \circ {\rm Ad}(u) = \varphi \circ {\rm Ad}(Ju^*J)$ for all $u \in \mathcal G$.
\item $\varphi_{| \X} = 0$,
\end{enumerate}
then $\varphi_{| \bS_{\X}(M)}$ is $\mathcal G''$-central, so that the inclusion $\mathcal G'' \subset M$ is not properly proximal.
\end{thm}
\begin{proof}
Note that since $\varphi_{|M} = \tau$ is normal, the set of elements $x \in M$ such that $[ x, \varphi_{| \bS_{\X}(M)} ] = 0$ forms a von Neumann subalgebra of $M$, thus it suffices to show that $\mathcal G$ is contained in this set. Also, as $\varphi\in \B(L^2M)^{M\sharp M}_J$ and $\varphi_{| \X} = 0$ we have $\varphi_{|\K_{\X}^{\infty, 1}(M)} = 0$.

If $T \in \bS_{\X}(M)$ and $u \in \mathcal G$, then $ ( J u^* J ) T (J u J) - T = [J u^* J, T] (JuJ) \in \K_{\X}^{\infty, 1}(M)$. Therefore,
\[
\varphi \circ {\rm Ad}(u) (T) 
= \varphi \circ {\rm Ad}(Ju^*J)(T)
= \varphi( T).
\]
\end{proof}

\begin{cor}\label{cor:gamma}
Let $M$ be a tracial von Neumann algebra with property (Gamma). Then $M$ is not properly proximal.
\end{cor}
\begin{proof}
Suppose $M$ has property (Gamma). Let $u_n \in \mathcal U(M)$ such that $u_n \to 0$ ultraweakly and $\| [x, u_n] \|_2 \to 0$ for all $x \in M$. Let $\varphi$ be any weak$^*$-limit point of the states $\B(L^2M) \ni T \mapsto \langle T \hat{u_n}, \hat{u_n} \rangle$. Then it is easy to see that $\varphi$ satisfies the hypotheses of Theorem~\ref{thm:ppcommute}, for $\mathcal G = \mathcal U(M)$ and $\X = \K(L^2M)$.
\end{proof}

\begin{lem}\label{lem:regularmix}
Let $M$ be a finite von Neumann algebra and $Q \subset M$ a regular von Neumann subalgebra. If $u_n \in \mathcal U(M)$ is such that $\| E_Q(a u_n b) \|_2 \to 0$ for all $a, b, \in M$, then for all $S, T \in \B(L^2M)$ of the form $a JbJ e_Q JcJ d$ with $a, b, c, d \in M$ we have $\| S^* u_n T \|_{\infty, 1} \to 0$. 
\end{lem}
\begin{proof}
Suppose $\{ u_n \}_n \subset \mathcal U(M)$ is given as above. Since $\| T x JyJ \|_{\infty, 1} \leq \| T \|_{\infty, 1} \| x \| \| y \|$, it is enough to check that  $\| S^* u_n T \|_{\infty, 1} \to 0$ when $S$ and $T$ are each of the form $a JbJ e_Q$, with $a, b \in M$. 

Also, note that if $a, b, c, d \in M$, then for $x \in M$ we have 
\[
\| e_Q a JbJ x JcJ d e_Q \|_{\infty, 1} 
= \sup_{y, z \in (Q)_1} | \tau(a x d y c^* b^* z) |
\leq \| b c \|_2 \| ax d \|_2.
\]
It therefore follows that by taking spans and using density in $\| \cdot \|_2$, to prove the lemma it suffices to show that $\| S^* u_n T \|_{\infty, 1} \to 0$ when $S$ and $T$ are each of the form $e_Q a JbJ$, where $a \in M$ and $b \in \mathcal N_M(Q)$. Finally, note that for $b \in \mathcal N_M(Q)$ we have $e_Q JbJ = Jb J b e_Q b^*$, and from this we then see that it suffices to consider the case when $S$ and $T$ are each of the form $e_Q a$ for $a \in M$. This is easily verified, for if $a_1, a_2 \in M$, then 
\begin{align*}
\lVert e_Q a_1^* u_n  a_2 e_Q \rVert_{\infty, 1}
&= \lVert E_Q(a_1^* u_n c_a) \rVert_1  \to 0.
\end{align*}
\end{proof}

For the next theorem, recall from Example~\ref{examp:bdrysub} that if $Q \subset M$ is a von Neumann subalgebra, then we denote by $\X_Q$ the $M$-boundary piece consisting of the norm closed span of all operators of the from $x_1 J y_1 J T J y_2 J x_2$, with $x_1, x_2, y_1, y_2 \in M$ and $T \in e_Q \B(L^2M) e_Q$.

\begin{thm}\label{thm:weaklycompact}
Let $M$ be a finite von Neumann algebra, and $Q \subset M$ a regular von Neumann subalgebra such that $M$ is properly proximal relative to the $M$-boundary piece $\X_Q$. If $P \subset M$ is a weakly compact regular von Neumann subalgebra, then $P \prec_M Q$. In particular, if $M$ is properly proximal, then $M$ has no diffuse weakly compact regular von Neumann subalgebras. 
\end{thm}

\begin{proof}
By the weak compactness of $P\subset M$, there exists a state $\varphi: \B(L^2M)\to \mathbb C$ satisfying the following properties:

(i) $\varphi$ is the canonical normal trace on $M$ and $JMJ$;\\
(ii) $\varphi(xT)=\varphi(Tx)$ for all $T\in \B(L^2M)$, $x\in P$;\\
(iii) $\varphi \circ {\rm Ad}(u)(T)=\varphi \circ {\rm Ad}(Ju^*J)(T)$ for all $T\in \B(L^2M)$ and $u\in \mathcal N_M(P)$.

If $P \not\prec_M Q$, then by Lemma~\ref{lem:regularmix} for each $T \in \B(L^2M)$ of the form $a^* Jb^*J e_Q JbJ a$ with $a, b \in M$ there exists a sequence $\{ u_n \}_n \subset \mathcal U(P)$ such that $\lVert T^* u_m^*u_n T \rVert_{\infty, 1}, \lVert T^* u_n^* u_m T \rVert_{\infty, 1} < 2^{-n}$ whenever $n > m$, and such that $\| T \widehat{u_n^*} \|_2 < 2^{-n}$ for each $n \geq 1$. We then have 
\begin{align*}
\lVert \frac{1}{N} \sum_{n = 1}^N u_n^* T u_n \rVert_{\infty, 2}^2 
&\leq \frac{1}{N^2} \sum_{n, m = 1}^N \lVert u_m^* T u_m u_n^* T u_n \rVert_{\infty, 1} \\
&\leq \frac{2}{N^2} \sum_{1 \leq m < n \leq N} 2^{-n} +  \frac{1}{N^2} \sum_{n = 1}^N \lVert u_n^* T^2 u_n \rVert_{\infty, 1} \\
&\leq \frac{2}{N} + \frac{\| T \|^2}{N} \to 0. 
\end{align*}
Since $\varphi$ is continuous in the $\| \cdot \|_{\infty, 2}$-norm on bounded sets it then follows from (ii)  that 
\[
\varphi(T) = \lim_{N \to \infty} \frac{1}{N} \sum_{n = 1}^N \varphi(u_n^* T u_n) = 0.
\] 
If we now have $a, b \in M$ and $T \in e_Q \B(L^2M) e_Q$ with $T \geq 0$, then 
\[
\varphi(a^* Jb^*J T Jb J a) \leq \| T \| \varphi( a^* Jb^*J e_A JbJ a ) = 0
\] 
and hence $\varphi(a^* Jb^*J T Jb J a) = 0$. By polarization it then follows that $\varphi(a J b J T J c J d) = 0$ for all $a, b, c, d \in M$ and $T \in e_Q \B(L^2M) e_Q$. Since the span of such elements is norm dense in $\X_Q$ it follows that $\varphi_{|\X_Q} = 0$, and hence $M$ is not properly proximal relative to $\X_Q$ by Theorem~\ref{thm:ppcommute}.
\end{proof}

\subsection{Proper proximality relative to the amenable boundary piece}

We see from Theorem~\ref{thm: iff} that a group $\Gamma$ is properly proximal if and only if the group von Neumann algebra $L\Gamma$ is properly proximal. In this section give an application of the development of boundary pieces for von Neumann algebras by showing that this also holds for proper proximality relative to a canonical ``amenable'' boundary piece. 

Let $\Gamma$ be a group, and let $\pi:\Gamma \to \mathcal U(\mathcal H)$ be a universal representation that is weakly contained in the left regular representation, i.e., $\pi$ is the restriction of the universal representation of $C^*_\lambda \Gamma$. Associated to this representation is a boundary piece $X_{\rm amen}$ as described in \cite{BoIoPe21}, where a net $(t_i)_i$ has a limit in $X_{\rm amen} \subset \beta \Gamma$ if and only if $\pi(t_i) \to 0$ ultraweakly, i.e., $t_i \to 0$ in the weak-topology in $C^*_\lambda \Gamma$. Alternatively, we can view the corresponding ideal $I_{\rm amen} \subset \ell^\infty \Gamma$ as the ideal generated by the set $B_r(\Gamma)$ of all matrix coefficients $\varphi_{\xi, \eta}(t) = \langle \pi(t) \xi, \eta \rangle$ for $\xi, \eta \in \mathcal H$. Note that since weak containment is preserved under tensor products it follows from Fell's absorption principle that $B_r(\Gamma)$ is an ideal in the Fourier-Stieltjes Algebra $B(\Gamma)$. In particular, $B_r(\Gamma)$ is a self-adjoint subalgebra of $\ell^\infty \Gamma$, and so $f \in I_{\rm amen}$ if and only if $| f | \leq g$ for some $g \in \overline{B_r(\Gamma)}$. 

If $\Gamma$ is a nonamenable group, then the following lemma shows that we often have $c_0 \Gamma \subsetneq I_{\rm amen} \subsetneq \ell^\infty \Gamma$.  

\begin{lem}\label{lem:amenmix}
Suppose $\Sigma < \Gamma$ is a subgroup. Then $\Sigma$ is amenable if and only if $1_{\Sigma} \in I_{\rm amen}$. 
\end{lem}
\begin{proof}
If $\Sigma$ is amenable, then the quasi-regular representation $\ell^2 (\Gamma/\Sigma)$ is weakly contained in the left regular representation and $1_{\Sigma}$ is a matrix coefficient. Conversely, if $1_{\Sigma} \in I_{\rm amen}$, then there exists $\varphi \in B_r(\Gamma)$, say $\varphi(t) = \langle \pi(t) \xi, \eta \rangle$ so that $1_{\Sigma} \leq g$, where $g$ is some element in $\overline {B_r(\Gamma)}$ with $\| g - \varphi \|_\infty < 1/2$. Hence for all $t \in \Sigma$ we have $\Re( \langle \pi(t) \xi, \eta \rangle ) \geq 1/2$. If we let $\xi_0$ denote the minimal norm element in the closed convex hull of $\{ \pi(t) \xi \}_{t \in \Sigma}$, it then follows that $\xi_0$ is a $\Sigma$-invariant vector and is non-zero since $\Re( \langle \xi_0, \eta \rangle ) \geq 1/2$. It therefore follows that the trivial representation for $\Sigma$ is weakly contained in $\pi \prec \lambda$ and hence $\Sigma$ is amenable. 
\end{proof}

We now fix a tracial von Neumann algebra $M$, and by analogy with above we consider a universal $M$-$M$ correspondence $\mathcal H$ that is weakly contained in the coarse correspondence $L^2M \ovt L^2M$. Note that we may assume that as an $M$-$M$ correspondence we have $\mathcal H \cong \overline{\mathcal H}$ and we then have a boundary piece $\X_{\rm amen} = \X_{\mathcal H}$ as defined in Example~\ref{examp:correbdry}

Similar to Lemma~\ref{lem:amenmix}, we have the following lemma.

\begin{lem}\label{lem:amenbdry}
Let $M$ be a tracial von Neumann algebra and $A \subset M$ a von Neumann subalgebra such that $A$ does not have an amenable direct summand, and $\mathcal G \subset \mathcal U(A)$ a subgroup that generates $A$ as a von Neumann algebra. Then there exists a net $\{ u_i \}_i \subset \mathcal G$ such that $u_i \otimes u_i^{\rm op} \to 0$ in the $\X_{\rm amen}$-topology.  
\end{lem}
\begin{proof}
Fix a universal $M$-$M$ correspondence $\mathcal H$ that is weakly contained in the coarse correspondence. Note that just as in the case for groups the span $B_0$ of operators of the form $T_\xi$ for some bounded vector $\xi \in \mathcal H$ forms a $*$-subalgebra of $\B(L^2M)$. Also, $M$ and $JMJ$ are contained in the multiplier algebra of $B_0$ in $\B(L^2M)$. Thus if we denote by $B = \overline{B_0}$, then we have $\X_{\rm amen} = B \B(L^2M) B$. In particular, a net $\{u_i \}_i \subset \mathcal G$ satisfies $u_i \otimes u_i^{\rm op} \to 0$ in the $\X_{\rm amen}$-topology if and only if $u_i T_{\xi_1} S T_{\xi_2} u_i^* \to 0$ ultraweakly for any $S \in \B(L^2M)$ and $\xi_1, \xi_2 \in \mathcal H$ bounded vectors. Moreover, by the polarization identity we see that this is if and only if for each bounded vector $\xi \in \mathcal H$, $S \in \B(L^2M)$, and $a \in M$ we have $\| S T_{\xi} JaJ \widehat{u_i^*} \|_2^2 \to 0$, and for this it suffices to consider only the case when $S = 1$. 

Therefore, if no such net of unitaries existed, then there would exist bounded vectors $\xi_1, \ldots, \xi_n \in \mathcal H$, $a_1, \ldots, a_n \in M$ and $c > 0$ such that for all $u \in \mathcal G$ we would have $\sum_{k = 1}^n \| T_{\xi_k} Ja_k J \hat{ u } \|_2^2 \geq c$. We let $\xi = \oplus_{k = 1}^n \xi_k \in \oplus_{k = 1}^n \mathcal H \oovt{M} \overline {\mathcal H}$ and we let $\tilde \xi = \oplus_{k = 1}^n a_k^* \xi_k a_k$. We also let $\mathcal C$ denote the closed convex set generated by vectors of the form $u^* \xi u$, for $u \in \mathcal G$. Then for each $u \in \mathcal G$ we have $\langle u^* \xi u, \tilde \xi \rangle = \sum_{k = 1}^n \| T_{\xi_k} J a_k J \hat{u} \|_2^2 \geq c$, and hence for any $\eta \in \mathcal C$ we have $\langle \eta, \tilde \xi \rangle \geq c > 0$. Hence, if we take $\eta \in \mathcal C$ to be the vector of minimal norm, then $\eta$ is a non-zero vector that is invariant under conjugation by $\mathcal G$. Since $\mathcal G$ generates $A$ it then follows that $\eta$ is a non-zero $A$-central vector. Since $\oplus_{k = 1}^n \mathcal H \oovt{M} \overline {\mathcal H}$ is weakly contained in the coarse correspondence this would then show that $Ap$ is amenable where $p \in \mathcal P(\mathcal Z(A))$ is the non-zero support projection for $\xi$.
\end{proof}

\begin{thm}
Let $\Gamma$ be a discrete group, then $I_{\rm amen} \subset \ell^\infty \Gamma$ is the closed ideal generated by $E_{\ell^\infty \Gamma}( \X_{\rm amen} )$ and $\X_{\rm amen} \subset \B(L^2 (L\Gamma))$ contains the boundary piece generated by $I_{\rm amen}$. 
\end{thm}
\begin{proof}
Let $\mathcal H$ be a universal $M$-$M$ correspondence that is weakly contained in the coarse correspondence, and let $\xi, \eta \in \mathcal H$ be bounded vector. Then $\mathcal H \oovt{M} \overline{\mathcal H}$ is also weakly contained in the coarse correspondence. In particular if we consider the representation $\pi: \Gamma \to \mathcal U(\mathcal H \oovt{M} \overline{\mathcal H})$ given by conjugation, then $\pi$ is weakly contained in the conjugation action associated to the coarse correspondence $L^2M \ovt L^2M$, which is easily seen to be a multiple of the left regular representation. For $t \in \Gamma$ we compute 
\[
E_{\ell^\infty \Gamma}( T_\xi^* T_\eta )(t)
= E_{\ell^\infty \Gamma}(T_{\eta \otimes \overline{\xi}})(t)
= \langle T_{\eta \otimes \overline{\xi}} \delta_t, \delta_t \rangle
= \langle u_t^* (\eta \otimes \overline{\xi}) u_t, \eta \otimes \overline{\xi} \rangle.
\]
Since $\pi \prec \lambda$ we then have $E_{\ell^\infty \Gamma}(T_\xi^* T_\eta) \in I_{\rm amen}$, and it then follows that $E_{\ell^\infty \Gamma}( \X_{\rm amen}) \subset I_{\rm amen}$. 

On the other hand, if $\pi: \Gamma \to \mathcal U(\mathcal H)$ is a representation that is weakly contained in the left regular representation, then we may consider the $L\Gamma$-$L\Gamma$ correspondence $\ell^2 \Gamma \ovt \mathcal H$ where the first copy of $L\Gamma$ acts as $\lambda_t \otimes \pi(t)$ (which is conjugate to $\lambda_t \otimes 1$ by Fell's absorption principle), and the second copy of $L\Gamma$ as as $\rho_t \otimes 1$. As is well known, if $\pi \prec \lambda$, then $\ell^2 \Gamma \ovt \mathcal H$ is weakly contained in the coarse correspondence. Also, if $\xi \in \mathcal H$, then it is easy to check that $\delta_e \otimes \xi$ is a bounded vector and $T_{\delta_e \otimes \xi}$ is the diagonal multiplication operator corresponding to $\Gamma \ni t \mapsto \langle \pi(t) \xi, \xi \rangle$. It then follows that $I_{\rm amen} \subset \X_{\rm amen}$.

Moreover, if $\varphi(t) = \langle \pi(t) \xi, \xi \rangle$ and $M_\varphi$ denotes the diagonal multiplication operator, then since $E_{\ell^\infty}(M_\varphi) = \varphi$ we see that $I_{\rm amen}$ is, in fact, equal to $E_{\ell^\infty \Gamma}(\X_{\rm amen})$. 
\end{proof}

\begin{cor}
Let $\Gamma$ be a discrete group, then $L\Gamma$ is properly proximal relative to $\X_{\rm amen}$ if and only if $\Gamma$ is properly proximal relative to $I_{\rm amen}$. 
\end{cor}
\begin{proof}
This follows from the previous theorem, together with Theorem~\ref{thm: iff}. 
\end{proof}

If $\Gamma$ is properly proximal, then clearly $\Gamma$ is properly proximal relative to $I_{\rm amen}$, but also if $\Gamma$ has a normal amenable subgroup $\Sigma \lhd \Gamma$ such that $\Gamma/ \Sigma$ is properly proximal, then it follows from \cite{BoIoPe21} that $\Gamma$ is also properly proximal relative to $I_{\rm amen}$. Examples of groups and von Neumann algebras that are not properly proximal relative to the amenable boundary piece are infinite direct products of nonamenable groups or infinite tensor products of nonamenable II$_1$ factors, which follows from the following proposition.

\begin{prop}
Suppose $M$ is a tracial von Neumann algebra and $B_n \subset M$ is a decreasing sequence of von Neumann subalgebras such that each $B_n$ has no amenable summand and such that $\cup_n (B_n' \cap M)$ is ultraweakly dense in $M$. Then $M$ is not properly proximal relative to $\X_{\rm amen}$. 
\end{prop}
\begin{proof}
This is similar to the proof of Corollary~\ref{cor:gamma}. Given any asymptotically central net $\{ u_i \}_i \subset \mathcal U(M)$ we may consider a state $\varphi$ on $\mathbb S_{\X_{\rm amen}}(M)$ which is a weak$^*$-limit point of the vector states $\mathbb S_{\X_{\rm amen}}(M) \ni T \mapsto \langle T \widehat{ u_i }, \widehat{u_i} \rangle$. Since each $B_n$ has no amenable summand, it follows from Lemma~\ref{lem:amenbdry} that we may find an asymptotically central net $\{ u_i \}_i \subset \mathcal U(M)$ so that $\varphi$ vanishes on $\X_{\rm amen}$, and it then follows that $\varphi$ satisfies the hypotheses of Theorem~\ref{thm:ppcommute}.
\end{proof}

\section{Biexact groups and properly proximal von Neumann algebras}\label{sec:biexact}
 
 Throughout this section we fix a group $\Gamma$ and we fix a trace-preserving action $\Gamma \actson B$ of $\Gamma$ on a tracial von Neumann algebra $B$. 
 
If $A$ is a $C^*$-algebra and $\Gamma \actson A$, then we denote by $A \rtimes_r \Gamma$ the reduced $C^*$-crossed product, and by $A \rtimes_f \Gamma$ the full $C^*$-crossed product.

\begin{thm}\label{thm:biexactpropp}
Suppose $\Gamma$ is biexact and $B$ is abelian. Then for every von Neumann subalgebra $P \subset B\rtimes \Gamma$, either the inclusion $P \subset B\rtimes \Gamma$ is properly proximal relative to $\X_B$ or else $P$ has an amenable direct summand.
\end{thm}
\begin{proof}
Suppose the inclusion $P \subset B \rtimes \Gamma$ is not properly proximal relative to $\X_B$, and let $\phi: \bS_{\X_B}(B \rtimes \Gamma) \to \langle p( B \rtimes \Gamma )p, e_{Pp} \rangle$ be a $p(B \rtimes \Gamma)p$-bimodular u.c.p.\  map, where $p \in \mathcal Z(P)$ is a non-zero central projection.

If we consider the $\Gamma$-equivariant diagonal embedding $\ell^\infty \Gamma \subset \B(\ell^2\Gamma) \subset \B(L^2B) \ovt \B(\ell^2 \Gamma)$, then we see that $\ell^\infty \Gamma$ commutes with $B$, and that $c_0 \Gamma$ is mapped to $\X_B$. Restricting to $\bS(\Gamma)$ then gives a $\Gamma$-equivariant embedding into $B' \cap \bS_{\X_B}(B \rtimes \Gamma)$. We therefore obtain a $*$-homomorphism $(B \otimes_{\rm min} \bS(\Gamma)) \rtimes_f \Gamma \to \B( L^2B \ovt \ell^2 \Gamma )$ whose image is contained in $\bS_{\X_B}(B \rtimes \Gamma)$. Composing this $*$-homomorphism with the u.c.p.\ map $\phi$ then gives a u.c.p.\ map $\tilde \phi: (B\otimes_{\rm min} \bS(\Gamma)) \rtimes_f \Gamma \to \langle p(B \rtimes \Gamma)p, e_{Pp} \rangle$ such that $\tilde \phi(t) = pu_tp$ for all $t \in \Gamma$. 

Since $\Gamma$ is biexact, the action $\Gamma \actson \bS(\Gamma)$ is topologically amenable, and since $B$ is abelian the action $\Gamma \actson B \otimes_{\rm min} \bS(\Gamma)$ is then also topologically amenable. Hence, $(B\otimes_{\rm min} \bS(\Gamma)) \rtimes_f \Gamma =  (B\otimes_{\rm min} \bS(\Gamma)) \rtimes_r \Gamma$ is a nuclear C$^*$-algebra.

We set $\varphi(\cdot):=\frac{1}{\tau(p)} \langle \tilde \phi(\cdot) \hat p, \hat p\rangle$ and note that for $x \in B \rtimes_r \Gamma$ we have $\varphi(x) = \frac{1}{\tau(p)} \tau(pxp)$. Since $B \rtimes_r \Gamma$ is ultraweakly dense in $B \rtimes \Gamma$ an argument similar to Proposition 3.1 in \cite{BoCa15} then gives a representation $\pi_\varphi:(B\otimes_{\rm min} \bS(\Gamma)) \rtimes_r \Gamma\to \B(\mathcal H_\varphi)$, a state $\tilde \varphi\in \B(\mathcal H_\varphi)_*$ with $\varphi = \tilde \varphi \circ \pi_\varphi$, and a projection $q\in  \pi_\varphi((B\otimes_{\rm min} \bS(\Gamma)) \rtimes_r \Gamma)''$ with $\tilde \varphi(q) = 1$ such that there is a normal unital $*$-homomorphism $\iota:B\rtimes \Gamma\hookrightarrow q \pi_\varphi((B\otimes_{\rm min} \bS(\Gamma)) \rtimes_r \Gamma)''q$.

Since $(B \otimes_{\rm min} \bS(\Gamma) ) \rtimes_r \Gamma$ is nuclear, we have that $\pi_\varphi((B\otimes_{\rm min} \bS(\Gamma)) \rtimes_r \Gamma)''$ is injective, and so there is a u.c.p.\ map $\tilde \iota: \langle B\rtimes \Gamma, e_B\rangle\to q \pi_\varphi((B\otimes_{\rm min} \bS(\Gamma)) \rtimes_r \Gamma)''q$ that extends $\iota$.
Notice that $\psi:=\tilde \varphi \circ \tilde \iota$ is then a $Pp$-central state on $\langle B\rtimes \Gamma, e_B\rangle$ and $\psi(x) = \frac{1}{\tau(p)} \tau(pxp)$ for all $x \in B \rtimes \Gamma$. 
By Theorem 2.1 in \cite{OzPo10I} we conclude that $Pp \oplus \mathbb C(1-p)$ is amenable relative to $B$, and since $B$ is abelian we then have that $Pp$ is amenable.
\end{proof}

\section{Properly proximal actions and crossed products}

Let $M$ be a separable finite von Neumann with normal faithful trace $\tau$ and let $\Gamma$ be a group that acts on $M$ by trace-preserving automorphisms. Note that the action $\Gamma \aactson{\alpha} (M, \tau)$ is unitarily implemented by the unitaries $\alpha^0_t( \hat{x} ) = \widehat{\alpha_t(x)}$. Conjugation by these unitaries then implements an action of $\Gamma$ on $\B(L^2(M, \tau))$, which preserves the space $\bS(M)$. 

If $\X \subset \B(L^2M)$ is a boundary piece such that conjugation by unitaries in the Koopman representation of $\Gamma$ preserves $\X$, then the action of $\Gamma$ on $\B(L^2M)$ will also preserve $\bS_\X(M)$. We say that the action $\Gamma \actson M$ is properly proximal relative to $\X$ if there does not exist a non-zero $\Gamma$-invariant projection $p \in \mathcal Z(M)$, such that $\Gamma \actson Mp$ is weakly compact, and such that there does not exist a $\Gamma$-invariant $M$-central state on $\bS_{\X}(M)$ that vanishes on $\K_{\X}^{\infty, 1}(M)$ and is normal when restricted to $M$.

Note that the same proof as in Corollary~\ref{cor:gamma} shows that if $\Gamma \actson M$ is properly proximal, then $(M \rtimes \Gamma)' \cap M^\omega$ cannot be diffuse. In particular, if $M$ is abelian and the action is ergodic, then the the action $\Gamma \actson M$ is strongly ergodic. Also note that if $\Gamma \actson \mathcal Z(M)$ is ergodic, then the trace $\tau$ is the unique $\Gamma$-invariant $M$-central normal state on $M$.

The definition of proper proximality for actions is slightly more technical than for von Neumann algebras. This is done to accommodate the proof of Lemma~\ref{lem:approxcentral} below. We will show in Proposition~\ref{prop:vanishonpiece} that proper proximality for the trivial action is equivalent to proper proximality for the von Neumann algebra.

Note that if $M$ is abelian and if $\varphi$ is a state on $\bS_{\X}(M)$ that restricts to the trace on $M$, then $\varphi$ vanishes on $\K_{\X}^{\infty, 1}(M)$ if and only if $\varphi$ vanishes on $\X$. Indeed, if $T \in \K_{\X}^{\infty, 1}(M)$, then for each $\varepsilon > 0$ there exists a projection $p \in \mathcal P(M)$ with $\tau(p) > 1-\varepsilon$ such that $p (J p J) T (J p J) p \in \X$. Since $M$ is abelian we have $J p J = p$ and by Cauchy-Schwarz $\varphi(T) \leq 2\varepsilon \| T \| + \varphi(p T p)$. Since $\varepsilon > 0$ is arbitrary this then shows that if $\varphi$ vanishes on $\X$, then it must also vanish on $\K_{\X}^{\infty, 1}(M)$. Thus, in the case when $M$ is abelian, the following proposition can give an easier criterion to check whether or not the action is properly proximal.

\begin{prop}\label{prop:vanish}
Let $(M, \tau)$ be a finite von Neumann algebra, suppose $\Gamma \actson (M, \tau)$ is a trace-preserving action, and $\X \subset \B(L^2M)$ is a $\Gamma$-invariant boundary piece for $M$.
If there exists a $\Gamma$-invariant, $M$-central state $\varphi: \bS_\X(M) \to \C$ such that $\varphi_{|M}$ is normal and $\varphi_{| \X} \not= 0$, then there exists a non-zero $\Gamma$-invariant projection $p \in \mathcal Z(M)$ so that $\Gamma\actson M p$ is weakly compact. 
\end{prop}
\begin{proof}
First notice that since $\bS_\X(M)=\bS_{\K_\X}(M)$, we may assume $\X=\K_\X$ and thus $M \subset M(\X)$.
Let $\varphi: \bS_\X(M) \to \C$ be a $\Gamma$-invariant, $M$-central state that satisfies $\varphi_{| M } = \tau$ and does not vanish on $\X$. Take a $\Gamma$-approximate unit $\{ A_i \}_i\subset \X$ that is quasi-central in $M(\X)$. 
Since $M$ is contained in the multiplier of $\X$ and since $\{ A_i \}_i$ is quasi-central in $M(\X)$, 
we obtain a non-zero $\Gamma$-invariant and $M$-central positive linear functional $\tilde \varphi$ on $\B(L^2M)$ with $\tilde \varphi_{| M}$ normal by setting
\[
\tilde \varphi(\cdot) = \lim_{i \to \omega} \varphi( A_i^{1/2} \cdot A_i^{1/2}),
\]
for some free ultrafilter $\omega$. Since $\tilde \varphi_{| M}$ is $\Gamma$-invariant, if we let $p$ be the support of $\tilde \varphi_{| M}$, then $p$ is $\Gamma$-invariant, and $\tilde \varphi_{|Mp}$ is faithful, so that $\Gamma \actson Mp$ is weakly compact.
\end{proof}

\begin{prop}\label{prop:vanishonpiece}
Let $M$ be a tracial von Neumann algebra, and $\X \subset \B(L^2M)$ a boundary piece. Then $M$ is not properly proximal relative to $\X$ if and only if there either exists a central projection $p \in \mathcal Z(M)$, such that $Mp$ is amenable, or else such that there exists an $M$-central state
$\psi: \bS_{\X} (M) \to \mathbb C$ such that $\psi_{| M}$ is normal and $\psi_{| \K_{\X}^{\infty, 1}(M)} = 0$. 
\end{prop}
\begin{proof}
Suppose $M$ is not properly proximal relative to $\X$ and that $M$ has no amenable direct summand. Let $\varphi: \mathbb S_\X(M) \to \mathbb C$ be an $M$-central state such that $\varphi_{|M}$ is normal. A simple standard argument shows that we may then produce a new $M$-central state (which we still denote by $\varphi$) and a central projection $p \in \mathcal Z(M)$ such that $\varphi(x) = \frac{1}{\tau(p)}\tau(xp)$ for $x \in M$. Since $M$ has no amenable direct summand, Proposition~\ref{prop:vanish} shows that $\varphi_{| \X} \not=0$.

Extend $\varphi$ to a state, still denoted by $\varphi$, on $\B(L^2 M)$. Since $\varphi_{|\mathcal Z(Mp)}= \frac{1}{\tau(p)} \tau$ by the Dixmier property we may produce a net of u.c.p.\ maps $\alpha_i: \B(L^2M) \to \B(L^2M)$ that are convex combinations of conjugation by unitaries in $JM J$ such that $\lim_{i \to \infty} | \varphi \circ \alpha_i(x) - \frac{1}{\tau(p)}\tau(x) | \to 0$ for each $x \in JMJp$. We may then let $\tilde \varphi$ be a weak$^*$-limit point of $\{ \varphi \circ \alpha_i \}_i$. Since each $\alpha_i$ is $M$-bimodular and preserves $\mathbb S_{\X}(M)$ (see the remark at the beginning of Section~\ref{sec:propprox}), it follows that $\tilde \varphi_{| \mathbb S_{\X}(M)}$ is again $M$-central. Also, since each $\alpha_i$ preserves $\X$ it follows that $\tilde \varphi$ again vanishes on $\X$. Since $\tilde \varphi$ restricts to the canonical traces on both $Mp$ and $JMJp$ it then follows from continuity that $\tilde \varphi$ also vanishes on $\K_{\X}^{\infty, 1}(M)$. 
\end{proof}

It was shown in \cite{IsPeRu19} that proper proximality for a group $\Gamma$ is an invariant for the orbit equivalence relation associated to any free probability measure-preserving action. Here we show a counterpart for proper proximality of an action. 

\begin{prop}\label{prop:normalizer}
		Let $\Gamma\aactson{\alpha} (M,\tau)$ be a trace-preserving action of $\Gamma$ on a separable finite von Neumann algebra $M$ with a normal faithful trace $\tau$. Suppose $\X \subset \B(L^2(M, \tau))$ is a $\Gamma$-boundary piece. Then $\Gamma\actson M$ is properly proximal relative to $\X$ if and only if $\mathcal{N}_{M\rtimes\Gamma}(M)\actson M$ is properly proximal relative to $\X$. 
		In particular, proper proximality for a free probability measure-preserving action is an an invariant of the orbit equivalence relation. 
\end{prop}
\begin{proof}
Since $\Gamma \actson M$ is weakly compact if and only if $\mathcal N_{M \rtimes \Gamma}(M) \actson M$ is weakly compact \cite[Proposition 3.4]{OzPo10I}, we may assume that there is no $\Gamma$-invariant non-zero projection $p \in \mathcal Z(M)$ such that $\Gamma \actson Mz$ is weakly compact. 

	Suppose $\Gamma\actson M$ is not properly proximal relative to $\X$ and let $\varphi$ be a $\Gamma$-equivariant and $M$-central state on $\bS_\X(M)$ with $\varphi_{\mid M}$ normal and $\varphi_{| \K_{\X}^{\infty, 1}(M)} = 0$.  We claim that $\varphi$ is $\cN_{M\rtimes\Gamma}(M)$-equivariant. Indeed, take any $u\in \cN_{M\rtimes\Gamma}(M)$ and by \cite{Dy59} there exists a partition of unity $\{p_t\}_{t \in \Gamma}\subset \mathcal P( \mathcal Z(M))$ and unitaries $\{ v_t \}_{t \in \Gamma} \subset \mathcal U( M )$ so that $\{ \alpha_{t^{-1}}(p_{t}) \}_{t \in \Gamma}$ also forms a partition of unity and $u=\sum_{t \in \Gamma} p_t v_t u_t$.
	Let $\alpha_u^0 \in\B(L^2M)$ be a unitary given by $\alpha_u^0(\hat x)=\widehat {uxu^*}$. Notice that $\alpha_u^0 =\sum_{t \in \Gamma} p_t v_t Jv_t J\alpha_{t}^0$. For $F \subset \Gamma$ finite we also set $\alpha_{u, F}^0 = \sum_{t \in F} p_t v_t Jv_t J\alpha_{t}^0$. 
	We may extend $\varphi$ to $\B(L^2M)$, still denoted by $\varphi$, and compute that for any $T\in \B(L^2M)$ and $F \subset \Gamma$ finite
	\begin{align*}
		 |\varphi(\alpha_u^0 T ( \alpha_u^0)^*)-\varphi(\alpha_{u, F}^0 T ( \alpha_{u, F}^0)^*)|
		&\leq  |\varphi((\alpha_u^0 -\alpha_{u, F}^0)T ( \alpha_u^0 )^*)|+|\varphi(\alpha_{u, F}^0 T(\alpha_u^0-\alpha_{u, F}^0 )^*|\\
		&\leq  |\varphi((\alpha_u^0-\alpha_{u, F}^0)(\alpha_u^0-\alpha_{u, F}^0)^*)|^{1/2}|\varphi(\alpha_u^0T^*T ( \alpha_u^0 )^*)|^{1/2} \\
		& \qquad +  |\varphi(\alpha_{u, F}^0 TT^* (\alpha_{u, F}^0)^*)|^{1/2}| \varphi((\alpha_u^0-\alpha_{u, F}^0)(\alpha_u^0-\alpha_{u, F}^0)^*)|^{1/2}\\
		&\leq   2\|T\|\varphi(\sum_{t \not\in F} p_t)^{1/2}.
	\end{align*}		
	
	Notice that for $S\in\bS_\X(M)$, we have $\alpha_{u, F}^0 S (\alpha_{u, F}^0)^*=(\sum_{t \in F}p_t)\alpha_u^0 S(\alpha_u^0)^*(\sum_{t \in F} p_t)\in \bS_\X(M)$ and since $\varphi_{|\bS_\X(M)}$ is $M$-central we have
	\begin{align*}
	\varphi(\alpha_{u, F}^0 S(\alpha_{u, F}^0)^*)&=\sum_{t \in F}\varphi(p_tv_tJv_tJ\alpha_{t}^0 S (\alpha_{t}^0)^* Jv_t^*J v_t)
	=\sum_{t \in F}\varphi(p_tv_t\alpha_{t}^0 S (\alpha_{t}^0)^* v_t^*)+\varphi(K_t)\\
	&=\sum_{t \in F}\varphi(\alpha_{t}^0 \alpha_{t^{-1}}(p_t) S (\alpha_{t}^0)^*)+\varphi(K_t)
	=\varphi(\sum_{t \in F}\alpha_{t^{-1}}(p_t) S)+\varphi(\sum_{t \in F}K_t),
	\end{align*}
	where $K_t=p_tv_t[Jv_tJ,\alpha_{t}^0 S( \alpha_{t}^0)^*] Jv_t^*J v_t\in \K_\X^{\infty,1}(M)$. 
	Since $\varphi$ vanishes on $\K^{\infty,1}_\X(M)$ we then have 
	\[
	|\varphi(S)-\varphi(\alpha_u S\alpha_u^*)|\leq 2\|S\|\varphi(\sum_{t \not\in F} p_t)^{1/2}+\|S\|\varphi(1-\sum_{t \in F}\alpha_{t^{-1}}(p_t))^{1/2}
	= 3 \| S \| \varphi( \sum_{t \not\in F} p_t )^{1/2}.
	\]
	Since $F \subset \Gamma$ was an arbitrary finite set and $\varphi_{|M}$ is normal we conclude that $\varphi_{| \bS_\X(M)}$ is $\cN_{M\rtimes\Gamma}(M)$-invariant.
\end{proof}

The following lemma is a simple application of the Hahn-Banach Theorem.

\begin{lem}\label{lem:paradoxdecomp}
Let $E$ be an operator system. Suppose $\Gamma \aactson{\alpha} E$ is an action of $\Gamma$ on $E$ by complete order isomorphisms. If $E$ does not admit a $\Gamma$-invariant state, then there exist $F \subset \Gamma$ finite, $\{T_t\}_{t \in F}\subset E$, and $T_0\in E_+$ such that $\|1 + T_0 + \sum_{t \in F} (T_t - \alpha_t(T_t) ) \|<1/4$.
\end{lem}

If $\X \subset \B(L^2M)$ is a boundary piece, then we set $\K_\X(M)_J^\sharp = \K_\X(M)^{M \sharp M} \cap \K_\X(M)^{JMJ \sharp JMJ}$. Recall from Proposition~\ref{prop:MCalgebra} that $(\B(L^2 M)_J^{{\sharp}} )^*$ is a von Neumann algebra containing $M$ and $JMJ$ as von Neumann subalgebras, and $( \K_\X(M)_J^\sharp )^*$ is a corner of this von Neumann algebra with the support projection of $( \K_\X(M)_J^\sharp )^*$ commuting with both $M$ and $JMJ$. The following technical lemma is adapted from Theorem 4.3 in \cite{BoIoPe21}.

\begin{lem}\label{lem:approxcentral}
Let $M$ be a separable finite von Neumann algebra with normal faithful trace $\tau$. 
Suppose that $\Gamma \aactson{\alpha} M$ is trace-preserving action by a group $\Gamma$ such that $\Gamma \actson \mathcal Z(M)$ is ergodic and let $\X \subset \B(L^2M)$ be a $\Gamma$-boundary piece

Then $\Gamma \actson M$ is properly proximal relative to $\X$ if and only if there is no $M$-central $\Gamma$-invariant state $\varphi$ on
\[
\widetilde{\bS}_{\X}(M) := \left\{ T \in \left(\B(L^2 M)_J^{{\sharp}} \right)^* \mid [T, a] \in \left( \K_\X(M)_J^\sharp \right)^* \ {\rm for \ all \ } a \in JMJ \right\}
\] 
such that $\varphi_{| M } = \tau$.
\end{lem}
\begin{proof}
First note that as $\bS_{\X}(M) = \bS_{\K_{\X}}(M)$ and $\widetilde{\bS}_{\X}(M) = \widetilde{\bS}_{\K_{\X}}(M)$ we may assume that $\X = \K_{\X}$. In particular, we may assume that both $M$ and $JMJ$ are contained in the multiplier algebra of $\X$.  
Also, note that by considering the natural map from $\B(L^2M)$ into $(\B(L^2M)_J^\sharp )^*$ we get $\Gamma$-equivariant $M$-bimodular u.c.p.\ maps $\B(L^2M) \to \left( \K_\X(M)_J^\sharp \right)^*$, and $\bS_{\X}(M) \to \widetilde{\bS}_{\X}(M)$. Thus, if there exists an $M$-central $\Gamma$-invariant state $\varphi$ on $\widetilde{\bS}_{\X}(M)$ such that $\varphi_{|M} = \tau$, then either $\varphi_{| \left( \K_\X(M)_J^\sharp \right)^*} \not= 0$, in which case restricting to $\B(L^2M)$ shows that $\Gamma \actson M$ is weakly compact \cite[Proposition 3.2]{OzPo10I} (and hence not properly proximal),
 or else $\varphi_{| \left( \K_\X(M)_J^\sharp \right)^*} = 0$, in which case restricting $\varphi$ to $\bS_\X(M)$ shows that $\Gamma \actson M$ is not properly proximal.

We now show the converse. So suppose that $\widetilde{\bS}_{\X}(M)$ has no $M$-central $\Gamma$-invariant state that restricts to the trace on $M$. Set $\mathcal G = \mathcal N_{M \rtimes \Gamma}(M)$ and continue to denote by $\alpha$ the trace-preserving action of $\mathcal G$ on $M$ given by conjugation. Since $\Gamma \actson \mathcal Z(M)$ is ergodic it follows from \cite{dy63} that the trace is the unique $\mathcal G$-invariant state on $M$, and hence there exists no $\mathcal G$-invariant state on $\widetilde{\bS}_{\X}(M)$. 

Let $M_0 \subset M$ be an ultraweakly dense $*$-subalgebra that is generated, as an algebra, by a countable set of contractions $S_0 = \{ x_0 = 1, x_1, \ldots \}$.
By Lemma~\ref{lem:paradoxdecomp} there exists a finite set $F \subset \mathcal G$, $\{ T_g \}_{g \in F} \subset  \widetilde{\bS}_{\X}(M)$, and $T_0 \in \widetilde{\bS}_{\X}(M)_+$, so that $\| 1 + T_0 + \sum_{g \in F} ( T_g - \alpha_g(T_g)  ) \|<1/4$. By enlarging $F$ (possibly with repeated elements) we may assume $\| T_g \| \leq 1$ for each $g \in F$. 

We let $I$ denote the set of all pairs $(A, \kappa)$ where $A \subset \B(L^2M)^\sharp_J$ is a finite set and $\kappa > 0$. We define a partial ordering on $I$ by $(A, \kappa) \prec (A', \kappa')$ if $A \subset A'$ and $\kappa' \leq \kappa$. 

For each $T_g$, we may find a net $\{T_g^i\}_{i\in I}\subset \B(L^2 M)$ whose weak$^*$ limit is $T_g$; moreover, by Goldstein's theorem, we may assume $\|T_g^i\|\leq 1$ for all $i\in I$. Since an operator $T \in \left( {\B(L^2M)}_J^{{ \sharp }} \right)^*$ is a positive contraction if and only if $\varphi(T) \in [0, 1]$ for all states $\varphi \in {\B(L^2M)}_J^{{ \sharp }}$, and since the set of positive contractions is convex, it then follows from the Hahn-Banach theorem that any positive contraction in $\left( {\B(L^2M)}_J^{{ \sharp }} \right)^*$ is in the ultraweak closure of positive contractions in $\B(L^2M)$. We may therefore also find a net of positive operators $\{ T_0^i \}_{i \in I} \subset \B(L^2M)$ so that $\| T_0^i \| \leq \| T_0 \|$, and $T_0^i \to T_0$ weak$^*$.  

We also choose a net $\{ S_i \}_{i \in I} \in  \B(L^2M)$ such that $\| S_i \| < 1/4$, for each $i \in I$, and $S_i \to 1 + T_0 + \sum_{g \in F} ( T_g - \alpha_g(T_g)  )$ weak$^*$. 

For any given $g \in \{ 0 \} \cup F$ and $n\in \mathbb N$, we have $[T_g, Jx_nJ]\in \left( ( \K_{\X} )_J^{ \sharp } \right)^*$ and so we may choose a net $\{K_{g,n}^i\}_{i\in I}\subset \K_{\X}$ such that $\lim _{i \to \infty} K_{g,n}^i-[T_g^i, Jx_nJ]=0$ weak$^*$.
Furthermore, since ${\B(L^2M)}_J^{{ \sharp }}$ is the space of linear functionals that are continuous for both the $M$-$M$ and $JMJ$-$JMJ$-topologies we may pass to convex combinations and assume that there exists a net $\{ e_i^{n, g} \}_{i \in I}, \subset \mathcal P(M)$ such that $e_i^{n, g} \to 1$ strongly, and 
\begin{equation}\label{eq:0}
\lim_{i \to \infty} \|e_i^{n, g} Je_i^{n, g}J (K_{g,n}^i-[T_g^i, Jx_nJ] ) Je_i^{n, g}J e_i^{n, g} \|=0.
\end{equation}
Since we also have that $1 + T_0^i + \sum_{g \in F} ( T_g^i - \alpha_g(T_g^i)  ) - S_{i}$ converges weak$^*$ to $0$, we may assume that, in addition to (\ref{eq:0}), we also have
\begin{equation}
\left\| e_i^{n, g} Je_i^{n, g}J \left( 1 + T_0^i + \sum_{g \in F} ( T_g^i - \alpha_g(T_g^i)  ) - S_{i} \right) Je_i^{n, g}J e_i^{n, g} \right\| \to 0.
\end{equation}

We fix $\varepsilon > 0$ to be chosen later. For each $n\in \mathbb N$, we may then choose $i(n)\in I$ such that setting $f_n = \wedge_{g \in F, m \leq n} e_{i(n)}^{m, g}$, we have $\tau(f_{n}) > 1 - \varepsilon 2^{-n - 1}$, and 
\begin{equation}\label{eq:1}
\left\|f_n Jf_{n}J \left(K_{g,m}^{i(n)}-[T_g^{i(n)}, Jx_mJ] \right) Jf_{n}J f_{n} \right\|< 2^{-n},
\end{equation}
for all $m\leq n$ and $g \in F$, and also
\begin{equation}\label{eq:1S}
\left\| f_{n} J f_{n} J \left( 1 + T_0^{i(n)} + \sum_{g \in F} ( T_g^{i(n)} - \alpha_g(T_g^{i(n)})  ) - S_{i(n)} \right) Jf_{n}J f_{n} \right\| < 2^{-n}/16.
\end{equation}
Setting $p_{n} = \wedge_{k \geq n} f_{k}$ we then have that $\{ p_{n} \}_n$ is an increasing sequence, and $\tau(p_{n}) > 1 - \varepsilon 2^{-n }$, for all $n \geq 1$. 

By Lemma 13.7 in \cite{Da88} there exists $\delta_n > 0$ so that when $K \in  {\X}_+$, and $z \in M( {\X})$ are contractions satisfying $\| [K, z] \| < \delta_n$, then we have $\| [K^{1/2}, z ] \| < 2^{-n}/16(6 | F | + 4 + \| T_0 \|) < 2^{-n}$. 

If we consider the unital $C^*$-subalgebra $B$ generated by 
\[
\{ K_{g, m}^{i(n)}, Jx_nJ, a_k, p_n, Jq_nJ \mid n \geq 1, g \in F, m \leq n \},
\] 
then $B$ is a separable $C^*$-algebra and $B \cap {\X}$ is a closed ideal. It then follows from \cite{Arv77} that there exists an increasing sequence $\{F_n\}_{n \geq 1} \subset {\X}_+ \cap B$ such that that the following conditions are satisfied:
\begin{enumerate}[(a)]
\item\label{item:aproxu} $\|F_{n} K_{g,m}^{i(n + 1)}-K_{g,m}^{i(n+ 1)}\|< 2^{-n}$ for all $g \in \{ 0 \} \cup F$ and $m\leq n + 1$;
\item $\| [ F_n, J x_m J] \| < \delta_n/2$ for $m \leq n$;
\item\label{item:coma} $\| F_n - \alpha_g(F_n) \|< \delta_n/2$ for all $g \in F$;
\item\label{item:com} $\| [F_n, p_{m}] \|, \| [F_n, Jp_{m}J] \| \leq \delta_n/2$ for all $m \leq n$. 
\end{enumerate}

Set $E_1 = F_1^{1/2}$, $E_n=(F_n-F_{n-1})^{1/2}$ for $n \geq 2$, and define $T^0_g=\sum_n E_n T^{i(n)}_g E_n$ 
for $g \in \{ 0 \} \cup F$, where the sum is SOT convergent as $\|T_g^i\|\leq \|T_g\|$ for each $g$. We claim that $T^0_g\in \bS_{\X}(M)$. Since $ \K_{\X}^{\infty, 1}(M)$ is a strong $JMJ$-$JMJ$ bimodule it suffices to show that $[T_g^0, Jx_mJ]\in  \K_{\X}^{\infty, 1}(M)$ for any given $x_m\in S_0$.
We compute
\begin{equation}
\begin{split}
[T^0_g, Jx_m J]=& \sum_nE_n T^{i(n)}_g [E_n, Jx_mJ]+\sum_n [E_n, Jx_m J] T^{i(n)}_g E_n
\\
& +\sum_n E_n [T^{i(n)}_g, Jx_mJ]E_n.
\end{split}
\end{equation}
The finite sums in the summations $\sum_nE_n T^{i(n)}_k [E_n, Jx_mJ]$ and $\sum_n [E_n, Jx_m J] T^{i(n)}_g E_n$ belong to $ {\X}$ since ${\X}$ is hereditary. The summations $\sum_nE_n T^{i(n)}_g [E_n, Jx_mJ]$ and $\sum_n [E_n, Jx_m J] T^{i(n)}_g E_n$ are norm convergent since $\|[E_n, Jx_m J]\|< 2^{-n}$ for sufficiently large $n$ and hence they also belong to ${\X}$. 

Using a similar argument with (\ref{eq:1}) and (\ref{item:com}) we deduce 
\[
p_{\ell} J p_{\ell} J \left( \sum_n E_n K_{g,m}^{i(n)} E_n - \sum_n E_n [T^{i(n)}_g, Jx_mJ]E_n \right) J p_{\ell} J p_{\ell}  \in {\X},
\]
for any fixed $\ell \geq 1$, and hence
\begin{equation}
\left( \sum_n E_n K_{g,m}^{i(n)} E_n - \sum_n E_n [T^{i(n)}_g, Jx_mJ]E_n \right)  \in \K_{\X}^{\infty, 1}.
\end{equation}

On the other hand, since $E_n^2=F_n-F_{n-1}$, we have from (\ref{item:aproxu}) that $\| E_n K_{g, m}^{i(n)} E_n \| < 2^{-n + 1}$ for large enough $n$ and hence $\sum_n E_n K_{g,m}^{i(n)} E_n$ is also norm convergent, and thus contained in $\X$ as well. This then shows that $[T_g^0, J x_m J] \in \K_{\X}^{\infty, 1}(M)$, for all $g \in \{ 0 \} \cup F$, and $m \geq 0$, verifying the claim that $T_g^0 \in \bS_{\X}(M)$.  Note that we also have $T_0^0\in (\bS_{\X}(M) )_+$ since $T_0^0$ is the SOT-limit of positive operators.

We define $S = \sum_n E_n S_{i(n)} E_n$. From (\ref{eq:1S}), (\ref{item:coma}) and (\ref{item:com}) we have 
\begin{align*}
&\left\| p_1 J p_1 J \left(1 + T_0^0 + \sum_{g \in F} ( T_g^0 - \alpha_g(T_g^0)  ) - S  \right) J p_1 J p_1 \right\|
\\
&\leq 2\sum_{g \in F} \sum_n \|  E_n - \alpha_g(E_n) \| + 2\left( \sum_n \| [p_1, E_n] \|  +  \sum_n \| [ J p_1 J, E_n ] \| \right) (2 + \| T_0 \| + 2 | F | ) \\
 &   \qquad +\|\sum_nE_n p_1 J p_1 J ( 1 + T_0^{i(n)} + \sum_{g \in F} ( T_g^{i(n)} - \alpha_g(T_g^{i(n)}) )  - S_{i(n)} ) J p_1 J p_1 E_n\| \\
&\leq 	1/16 +\|\sum_nE_n p_1 J p_1 J \left(1 + T_0^{i(n)} + \sum_{g \in F} ( T_g^{i(n)} - \alpha_g(T_g^{i(n)}))   - S_{i(n)} \right) J p_1 J p_1 E_n\| 
\leq 1/8.
\end{align*}

Note that since $E_n \geq 0$ and $\sum_n E_n^2 = 1$, the map $\phi: \ell^\infty( \B(L^2M) ) \to \B(L^2M)$ given by $\phi( ( z_n ) ) = \sum_n E_n z_n E_n$ is unital completely positive. In particular, we have that $\phi$ is a contraction and it follows that $\| S \| \leq \sup_n \| S_{i(n)} \| \leq 1/4$. Thus, from the previous inequalities we conclude that
\begin{equation}\label{eq:contr}
\left\| p_1 J p_1 J \left( 1 + T_0^0 + \sum_{g \in F} ( T_g^0 - \alpha_t(T_g^0)  ) \right) J p_1 J p_1 \right\| < 3/8.
\end{equation}

Since $\tau(p_1) \geq 1 - \varepsilon$, it follows easily by considering the disintegration into factors that there exist a central projection $z \leq z(p_1)$, with $\tau(z) > 1 - 2\varepsilon$, and two subprojections $q_2, q_3 \leq z$ with $q_2 + q_3 = z$ such that $q_2$ and $q_3$ are both Murray-von Neumann subequivalent to $p_1$. Take partial isometries $v_2, v_3$ such that $v_2^*v_2, v_3^*v_3 \leq p_1$ and $v_2v_2^* = q_2$, $v_3 v_3^* = q_3$.  We define the u.c.p.\ map $\phi: \B(L^2M) \to \B(L^2M)$ by $\phi(T) = J v_2 J T J v_2^* J + J v_3 J T J v_3^* J + z^\perp T z^\perp$. Note that $\phi$ maps $\bS_{\X}(M)$ into itself, and for all $T \in \bS_{\X}(M)$ we have $\phi(T) - T \in \K_\X^{\infty, 1}(M)$. 

It follows from (\ref{eq:contr}) that
\begin{equation}\label{eq:contfin}
\left\| p_1 z \left( 1 + \phi( T_0^0 ) + \sum_{g \in F} ( \phi( T_g^0) - \phi( \alpha_g(  T_g^0) ) ) \right) z p_1 \right\| < 3/4.
\end{equation}

If we had a $\Gamma$-invariant state $\varphi$ on $\bS_{\X}(M)$ such that $\varphi_{| M} = \tau$ and $\varphi_{| \K_{\X}^{\infty, 1}(M)} = 0$, then from (\ref{eq:contfin}), the Cauchy-Schwarz inequality, and the fact that $\phi(T) - T \in \K_\X^{\infty, 1}(M)$ for each $T \in \bS_\X(M)$, we would then have
\begin{align*}
1 - 3\varepsilon 
&\leq \tau(p_1 z ) 
= \varphi( p_1 z ) \\
&\leq 3/4 + \varphi \left(  p_1 z \left( \phi( T_0^0 ) + \sum_{g \in F} ( \phi( T_g^0) - \phi( \alpha_g(  T_g^0) ) ) \right) z p_1 \right) \\
&= 3/4 + \varphi \left( p_1 z \left( T_0^0  + \sum_{g \in F} (  T_g^0 -  \alpha_g(  T_g^0)  ) \right) p_1 z \right) \\
&\leq 3/4 + 2 | F | (1 - \varphi( p_1 z) ) 
\leq 3/4 + 6 | F | \varepsilon.
\end{align*}
Thus, we obtain a contradiction by choosing $\varepsilon > 0$ so that $(3 + 6 | F |)\varepsilon < 1/4$.
\end{proof}

The following slight variant of of the previous lemma will be of more use in the sequel.

\begin{lem}\label{lem:approxcentral2}
Using the same notation as in the previous lemma, the action $\Gamma \actson M$ is properly proximal if and only if $\Gamma \actson M$ is not weakly compact and there is no $M$-central $\Gamma$-invariant state $\varphi$ on $\widetilde{\bS}_{\X}(M)$ such that $\varphi_{| M } = \tau$ and $\varphi_{|  \left( \K_\X(M)_J^\sharp \right)^*} = 0$. 
\end{lem}
\begin{proof}
Using the previous lemma, suppose $\varphi: \widetilde{\bS}_{\X}(M) \to \mathbb C$ is an $M$-central $\Gamma$-invariant state such that $\varphi_{|M} = \tau$. By considering the natural $M$ and $JMJ$-bimodular u.c.p.\ map from $\B(L^2M)$ into $ \left( \K_\X(M)_J^\sharp \right)^*$ we see that if $\varphi$ does not vanish on $\left( \K_\X(M)_J^\sharp \right)^*$, then we get a state on $\B(L^2M)$ that is $\Gamma$-invariant, $M$-central, and normal on $M$, which shows that $\Gamma \actson M$ is weakly compact by Proposition 3.2 in \cite{OzPo10I}.
\end{proof}

\begin{lem}\label{lem:subalgebra}
Let $M$ be a finite von Neumann algebra and $B \subset M$ a von Neumann subalgera. Let $e_B: L^2M \to L^2B$ denote the orthogonal projection. Then $e_B \K^{\infty, 1}(M) e_B \subset \K^{\infty, 1}(B)$, and $e_B \bS(M) e_B \subset \bS(B)$. 
\end{lem}
\begin{proof}
	For any $T\in\K^{\infty,1}(M)$, denote by $\{T_n\}_{n\in\mathbb N}\subset \K(L^2M)$ a sequence that converges to $T$ in $\|\cdot\|_{\infty,1}$. 
	Notice that 
	\[
	\|e_B(T- T_n) e_B\|_{\infty,1}=\sup_{a,b\in (B)_1} \langle e_B (T-T_n) e_B \hat a, b\rangle\leq \|T-T_n\|_{\infty,1},
	\]
	and we conclude that $e_B T e_B\in \K^{\infty,1}(B)$.
	Since $JBJ$ commutes with $e_B$, we also have $e_B \bS(M) e_B \subset \bS(B)$.
\end{proof}

\begin{thm}\label{thm:tensorproduct}
Let $M_1$ and $M_2$ be separable finite von Neumann algebras and suppose we have trace-preserving actions $\Gamma \actson M_i$. Then the action $\Gamma \actson M_1 \ovt M_2$ is properly proximal if and only if the action $\Gamma \actson M_i$ is properly proximal for each $i = 1, 2$.
\end{thm}
\begin{proof}
If the action $\Gamma \actson M_1 \ovt M_2$ is weakly compact, then so is $\Gamma \actson M_i$, for each $i$. We may therefore restrict to the case when $\Gamma \actson M$ is not weakly compact. Also, if $\mathcal Z(M_i)^\Gamma$ is diffuse for some $i$, then the same is true for $\mathcal Z(M_1 \ovt M_2)^\Gamma$, and the result is easy. We may therefore also restrict to the case when $\mathcal Z(M_i)^\Gamma$ is completely atomic for each $i = 1, 2$. It is easy to see that a direct sum of trace-preserving actions is properly proximal if and only if each summand is properly proximal, and hence by restricting to each atom we may further reduce to the case when $\Gamma \actson Z(M_i)$ is ergodic for each $i = 1, 2$. 

Suppose $\Gamma \actson M$ is not properly proximal.  By Lemma~\ref{lem:approxcentral2} there is an $M$-central $\Gamma$-invariant state $\varphi$ on $\widetilde {\bS}(M)$ with $\varphi_{|M}=\tau$ and $\varphi_{| (\K^{\infty,1}(M)^{M \sharp M}_J)^*}= 0$. Let $A_i \in \K(L^2M_1) \subset \B(L^2M_1)$ be an increasing quasi-central approximate unit. We consider the von Neumann algebra $(\B(L^2M)_J^{{M \sharp M}} )^*$ and let $A \in (\B(L^2M)_J^{{M \sharp M}})^*$ be the weak$^*$-limit of the increasing net $A_i \otimes 1$. We then define an $M_2$-bimodular $\Gamma$-equivariant completely positive map $\Phi: \bS(M_2) \to (\B(L^2M)_J^{{M \sharp M}})^*$ by $\Phi(T) = A^{1/2} (1 \otimes T) A^{1/2}$.

Note that for $x_1 \in M_1$ we have 
\[
[\Phi(T), Jx_1 J] = \lim_i [A_i, J x_1 J] \otimes T = 0,
\]
and for $x_2 \in M_2$ we have 
\[
[\Phi(T), Jx_2J] = \lim_iA_i \otimes [T, J x_2 J] \in (\K^{\infty,1}(M)^{M \sharp M}_J)^*.
\] 
Since $(\K^{\infty,1}(M)^{M \sharp M}_J)^*$ is a strong $JMJ$-bimodule, we have $[\Phi(T), J x J] \in(\K^{\infty,1}(M)^{M \sharp M}_J)^*$ for $x\in M$ in general.

Therefore we have $\Phi: \bS(M_2) \to \widetilde{\bS}(M)$. If $\varphi(A) \not= 0$, then composing $\Phi$ with $\varphi$ gives an $M_2$-central $\Gamma$-invariant positive linear functional on $\bS(M_2)$ that restricts to $\varphi(A) \tau$ on $M_2$. As we have $\Phi: \mathbb K^{\infty, 1}(M_2) \to (\K^{\infty,1}(M)^{M \sharp M}_J)^*$ this then shows that $\Gamma \actson M_2$ is not properly proximal. 

If we let $B_j \in \K(L^2M_2) \subset \B(L^2M_2)$ be an increasing quasi-central approximate unit and if we let $B \in (\B(L^2M)_J^{{M \sharp M}})^*$ be a weak$^*$-limit of the increasing net $1 \otimes B_j$, then just as above if $\varphi(B) \not= 0$ it would then follow that $\Gamma \actson M_1$ is not properly proximal. 

We now suppose $\varphi(A) = \varphi(B) = 0$ and consider the $M_1$-bimodular $\Gamma$-equivariant u.c.p.\ map $\Psi: \bS(M_1) \to (\B(L^2M)_J^{{M \sharp M}} )^*$ given by $\Psi(T) = (1 - A)^{1/2} (T \otimes 1) (1 - A)^{1/2}$. Similar to the calculation above, we have that $[\Psi(T), JxJ] = 0$ for all $x \in M_1 \otimes_{\rm alg} M_2$ and hence $[\Psi(T), JxJ] = 0$ for all $x \in M$. Therefore, $\Psi(T) \in  \widetilde{\bS}(M)$, and since $\varphi(A) = 0$ it then follows that $\varphi \circ \Psi$ defines an $M_1$-central $\Gamma$-equivariant state on $\bS(M_1)$ that restricts to the trace on $M_1$. Moreover, since $\varphi(B) = 0$ we have that $\varphi \circ \Psi_{| \K^{\infty, 1}(M_1)} = 0$ and hence $\Gamma \actson M_1$ is not properly proximal. 

Conversely, suppose $\Gamma \actson M_1$ is not properly proximal. Let $\varphi: \bS(M_1) \to \mathbb C$ be an $M_1$-central $\Gamma$-invariant state satisfying $\varphi_{| M_1} = \tau$ and $\varphi_{| \K^{\infty, 1}(M_1)} = 0$

We consider the u.c.p. map 
\[
{\rm Ad}_{e_1}: \B(L^2 (M_1 \ovt M_2) ) \to \B(L^2(M_1))
\] 
given by ${\rm Ad}_{e_1}(T) = e_1 T e_1$ where $e_1: L^2(M_1 \ovt M_2) \to L^2(M_1)$ is the orthogonal projection. By Lemma~\ref{lem:subalgebra} ${\rm Ad}_{e_1}$ maps $\bS(M_1 \ovt M_2 )$ into $\bS(M_1)$ and maps $\K^{\infty, 1}(M_1 \ovt M_2)$ into $\K^{\infty, 1}(M_1)$. We obtain a state $\psi: \bS(M_1 \ovt M_2) \to \mathbb C$ by setting $\psi = \varphi \circ {\rm Ad}_{e_1}$. 

Note that $\psi_{|M} = \tau \circ E_{M_1} = \tau$. Also note that since ${\rm Ad}_{e_1}$ is $M_1$-bimodular and $\Gamma$-equivariant we see that $\psi$ is $M_1$-central. If $u \in \mathcal U(M_2)$, then we have $e_1 u JuJ = e_1$, and since $\varphi_{| \K^{\infty, 1}(M_1)} = 0$ it follows that for $T \in \bS(M_1 \ovt M_2)$ we have 
\begin{align*}
\psi(u T) 
&= \varphi( e_1 u T e_1) 
= \varphi( e_1 (u JuJ) T (Ju^*J) e_1 ) 
= \varphi( e_1 T (J u^* J) e_1) \\
&= \varphi( e_1 T (J u^* J)  (u JuJ) e_1) 
= \psi(Tu).
\end{align*}
Hence, $M_2$ is also in the centralizer of $\psi$, so that the centralizer of $\psi$ contains the entirety of $M_1 \ovt M_2 = W^*(M_1, M_2)$. Since $\psi_{|\K^{\infty, 1}(M_1 \ovt M_2)} = 0$ this then shows that $\Gamma \actson M_1 \ovt M_2$ is not properly proximal. 
\end{proof}

\begin{rem}
	Let $M$ be a properly proximal separable \RN{2}$_1$ factor. Then the amplification $M^t$ is also properly proximal for any $t>0$. Indeed, it follows from the previous theorem that $M\overline\otimes M \cong M^{t}\overline\otimes M^{1/t}$ is properly proximal and thus $M^t$ is also properly proximal.
	Also, if $N\subset M$ is a finite index subfactor in a separable \RN{2}$_1$ factor $M$, and if $N$ is properly proximal, then so is $\langle M, e_N \rangle$ since it is an amplification of $N$. Combining with Proposition~\ref{prop:co-amenable}, we then see that $M$ is properly proximal as well. 
\end{rem}

\begin{thm}
Let $\Gamma \actson (X, \mu)$ be a probability measure preserving action, and let $\pi: (X, \mu) \to (Y, \nu)$ be a factor map with $(Y, \nu)$ diffuse. If $\Gamma \actson (X, \mu)$ is properly proximal, then so is $\Gamma \actson (Y, \nu)$. 
\end{thm}
\begin{proof}
Suppose $\Gamma \actson (X, \mu)$ is properly proximal. We set $M = L^\infty(X, \mu)$ and $M_1 = L^\infty(Y, \nu)$, and we view $M_1 \subset M$ via the embedding $\pi^*(f) = f \circ \pi$. If we consider the map ${\rm Ad}_{e_1} : \B(L^2M ) \to \B(L^2(M_1))$ as in the second half of the proof of Theorem~\ref{thm:tensorproduct}, then we see from the proof of Theorem~\ref{thm:tensorproduct} that the only property used to show that $\Gamma \actson M_1$ is properly proximal is that $M = W^*(M_1, M_1' \cap M)$, which obviously holds here since $M$ is abelian. 
\end{proof}

\begin{prop}\label{prop:crossedproduct implies action}
Suppose $\Gamma \aactson{\alpha} (M, \tau)$ is a trace-preserving action on a finite von Neumann algebra. If $M\rtimes\Gamma$ is properly proximal, then the action $\Gamma\actson M$ is properly proximal.
\end{prop}
\begin{proof}
This is also similar to the second half of the proof of Theorem~\ref{thm:tensorproduct}.
Suppose $\Gamma\actson M$ is not properly proximal, and let $\varphi: \bS(M) \to \C$ be a $\Gamma$-invariant $M$-central state such that $\varphi_{| M }$ is normal and $\varphi_{|\K^{\infty, 1}(M)} = 0$. 

Consider the u.c.p. map $\Ad_{e_M}:\B(L^2 (M\rtimes\Gamma))\to\B(L^2M)$ given by $\Ad_{e_M}(T)=e_M T e_M$, where $e_M: L^2M \overline\otimes \ell^2\Gamma\to L^2M$ is the orthogonal projection.
By Lemma~\ref{lem:subalgebra} we have $\Ad_{e_M}: \bS(M\rtimes\Gamma)\to \bS(M)$ and and $\Ad_{e_M}: \K^{\infty, 1}(M \rtimes \Gamma) \to \K^{\infty, 1}(M)$. We then consider the state $\psi=\varphi\circ\Ad_{e_M}$ on $\bS(M\rtimes\Gamma)$.

Note that $\psi_{|M}=\varphi\circ E_{M}$ is normal and $\psi$ is $M$-central since $\Ad_{e_M}$ is $M$-bimodular.
For any $t\in\Gamma$, denote by $u_t\in \mathcal U(M\rtimes\Gamma)$ and $\alpha_t^0\in \mathcal U(\B(L^2M))$ the corresponding unitaries.
Notice that $e_M u_t Ju_tJ=u_t Ju_tJ e_M$ and $u_t Ju_tJ_{\mid L^2M}=\alpha_t^0$. 
Since $\varphi_{|\K^{\infty, 1}(M)}=0$, it follows that for any $T\in \bS(M\rtimes \Gamma)$ we have
\begin{align*}
\psi(u_t T)&=\varphi(e_M u_t T e_M)=\varphi(e_M (u_t Ju_tJ) T(Ju_tJ)^* e_M)=\varphi(\alpha_t^0 e_M T (Ju_tJ)^* e_M)
\\	
&=\varphi(e_M T (Ju_tJ)^* e_M\alpha_t^0)=\varphi(e_M T (Ju_tJ)^* u_t Ju_tJ e_M)=\psi(T u_t).
\end{align*}
Therefore, $M\rtimes\Gamma$ is contained in the centralizer of $\psi$. Since $\psi_{| \K^{\infty, 1}(M \rtimes \Gamma)} = 0$ this then shows that $M\rtimes\Gamma$ is not properly proximal.
\end{proof}

We note here that non-strong ergodicity and the existence of weakly compact factors are not the only obstructions for an action proper proximality. For $d \geq 3$ the group $SL_d(\mathbb Z) \ltimes (\mathbb Z[\frac{1}{p}])^d$ is not properly proximal by \cite{IsPeRu19}, and since $SL_d(\mathbb Z)$ is properly proximal it then follows from Theorem~\ref{thm:crossedproduct} that the dual action $SL_d(\mathbb Z) \actson \widehat{( \mathbb Z[\frac{1}{p}])^d}$ is not properly proximal. However, this action has stable spectral gap (and hence is strongly ergodic and has no weakly compact diffuse factor) since it is a weak mixing action of a property (T) group. 

We also mention that proper proximality of $M\rtimes\Gamma$ does not imply that $\Gamma$ is properly proximal, as $\Lambda\wr\Gamma$ is always properly proximal for nontrivial $\Lambda$ and nonamenable $\Gamma$ \cite{DiKuEl21}.

\begin{thm}\label{thm:crossedproduct}
	Let $\Gamma\aactson{\alpha} (M,\tau)$ be a trace-preserving action of a countable group $\Gamma$ on a separable finite von Neumann algebra $M$ with a normal faithful trace $\tau$.
	If $\Gamma\actson (M,\tau)$ and $\Gamma$ are both properly proximal, then $M\rtimes\Gamma$ is properly proximal.
\end{thm}
\begin{proof}
	This is similar to the first half of the proof of Theorem~\ref{thm:tensorproduct}. We first note that if $\mathcal Z(M)^\Gamma$ is diffuse, then the result is easy. As above we may then restrict to the case when $\mathcal Z(M)^\Gamma$ is completely atomic, and then by considering atoms we may restrict to the case when $\Gamma \actson \mathcal Z(M)$ is ergodic. 
	
	Suppose then that $\Gamma \actson \mathcal Z(M)$ is ergodic and $M \rtimes \Gamma$ is not properly proximal, then by Lemma~\ref{lem:approxcentral2} there exists an $M$-central $\Gamma$-invariant state $\varphi$ on $\tilde \bS(M)$ such that $\varphi_{| M } = \tau$ and $\varphi_{|( \K^{\infty, 1}(M \rtimes \Gamma)_J^\sharp )^*} = 0$.	
	
	Let $A_i \in c_0(\Gamma) \subset \B(\ell^2\Gamma)$ be an increasing approximate unit that is almost $\Gamma$-invariant. Let $A$ denote the weak$^*$-limit of the increasing net $1 \otimes A_i$ when viewed inside of $\left( \B(L^2M \ovt \ell^2 \Gamma)_J^\sharp \right)^*$, where we identify $L^2(M \rtimes \Gamma)$ with $L^2M \ovt \ell^2 \Gamma$ by mapping $\widehat{x u_t}$ to $\hat x \otimes \delta_t$ for each $x \in M$ and $t \in \Gamma$. Note that under this identification we have that $M$ is represented as $M \otimes 1$, and if $x \in M$, then $J x J$ is represented by the operator $\sum_{t \in \Gamma} J_M \alpha_t(x) J_M \otimes P_t$ where $P_t$ denotes the rank-one projection onto $\mathbb C \delta_t \subset \ell^2 \Gamma$. In particular, we see that $M$ and $JMJ$ both commute with $\mathbb C \otimes \ell^\infty \Gamma$. 
	
	We then consider the u.c.p.\ map $\Phi: \B(L^2M) \to \left( \B(L^2M \ovt \ell^2 \Gamma)_J^\sharp \right)^*$ given by $\Phi(T) = A^{1/2} (T \otimes 1) A^{1/2}$. Since $A_i \in c_0(\Gamma)$ it then follows that $\Phi$ is both $M$ and $JMJ$-bimodular. Also since $A_i$ is almost $\Gamma$-invariant it follows that $\Phi$ is $\Gamma$-equivariant, where the action on $\left( \B(L^2M \ovt \ell^2 \Gamma)_J^\sharp \right)^*$ is given by ${\rm Ad}(u_t)$. We also have that the range of $\Phi$ maps into the subspace of invariant operators under the action of ${\rm Ad}(1 \otimes \rho_t)$. Hence, it follows that the restriction of $\Phi$ to $\bS(M)$ gives a $M$-bimodular $\Gamma$-equivariant u.c.p.\ map into $\bS(M \rtimes \Gamma)$.
	
	If $\varphi(A) \not= 0$, then considering $\varphi \circ \Phi$ we obtain an $M$-central $\Gamma$-invariant positive linear functional on $\bS(M)$ that restricts to $\varphi(A) \tau$ on $M$. Moreover, since $\varphi_{|( \K^{\infty, 1}(M \rtimes \Gamma)_J^\sharp )^*} = 0$ we then have that $\varphi \circ \Phi_{| \K^{\infty, 1}(M)} = 0$, which shows that $\Gamma \actson M$ is not properly proximal. 
	
	Otherwise, if $\varphi(A) = 0$, then we may consider the $\Gamma$-equivariant u.c.p.\ map $\Psi: \ell^\infty \Gamma \to \left( \B(L^2M \ovt \ell^2 \Gamma)_J^\sharp \right)^*$ given by $\Psi(f) = 1 \otimes (1 - A_n)^{1/2} M_f  (1 - A_n)^{1/2}$. As above, we see that the restriction of $\Psi$ to $\bS(\Gamma)$ maps into $\bS(M \rtimes \Gamma)$, and then the state $\varphi \circ \Psi$ defines a $\Gamma$-invariant state on $\bS(\Gamma)$, which shows that $\Gamma$ is not properly proximal. 	
\end{proof}

\section{Proper proximality and deformation/rigidity theory}

In this section we show how proper proximality for finite von Neumann algebras naturally fits in to the realm of Popa's deformation/rigidity theory (see \cite{Po07B}). We focus on two main sources of deformations. Malleable deformations in the sense of Popa \cite{Po06D}, and u.c.p.\ semigroups arising from closable derivations \cite{Sa90, DaLi92}. 

\subsection{Proper proximality via malleable deformations}

The following result shows how one can obtain proper proximality via Popa's malleable deformations. Malleable deformations arise in a number of different contexts in the theory of von Neumann algebras. We refer the reader to \cite{SaHaHoSi21}, and the references therein, for recent results regarding malleable deformations.

\begin{prop}\label{prop:deformation}
	Let $\Gamma \actson M$ be a trace-preserving action on a separable finite von Neumann algebra, that preserves a boundary piece $\X \subset \B(L^2M)$ and such that $\Gamma \actson \mathcal Z(M)$ is ergodic. Let $\Gamma \actson \tilde M$ be another trace-preserving action on a finite von Neumann algebra and suppose we have a tracial embedding $M \subset \tilde M$, such that the action of $\Gamma$ on $\tilde M$ extends the action on $M$. Suppose $\{\alpha_n\}_{n \geq 1}\subset\Aut(\tilde M)$ is a sequence of trace-preserving $\Gamma$-equivariant automorphisms such that $\| \alpha_n(x) - x \|_2 \to 0$ for all $x \in M$. 
		If $E_M\circ \alpha_n \in \K_\X^L(M, L^2M)$ for any $n \geq 1$, and if $\Gamma \actson M$ is not properly proximal relative to $\X$, then there exists a $\Gamma$-invariant $M$-central state on $(JM J)' \cap \B( L^2\tilde M\ominus L^2 M)$ that restricts to the trace on $M$, where $\Gamma$ acts on $\B( L^2\tilde M\ominus L^2 M)$ by conjugation with the Koopman representation. 
\end{prop}
\begin{proof}
We let $V_n$ denote the isometry from $L^2M$ to $L^2 \tilde M$ given by $V_n(\hat{x}) = \widehat{\alpha_n(x)}$ for $x \in M$. We also let $\phi_n: \B(L^2 \tilde M) \to \B(L^2M)$ be the $\Gamma$-equivariant u.c.p.\ map given by $\phi_n(T) = V_n^* T V_n$, for $T \in \B(L^2 \tilde M)$. 

If $x \in M$, then we have $\| xV_n  -  V_n x \|_{\infty, 2} = \sup_{a \in (M)_1} \|  x \alpha_n(a) - \alpha_n(xa) \|_2 \leq \| x - \alpha_n(x) \|_2$. It therefore follows that for all $x \in M$ and $T \in \B(L^2 \tilde M)$ we have 
\begin{equation}\label{eq:bimod1}
\| \phi_n(T x) - \phi_n(T) x \|_{\infty, 1}, \| \phi_n(xT) - x \phi_n(T) \|_{\infty, 1} \leq \| T \| \| x - \alpha_n(x) \|_2.
\end{equation}
We similarly have 
\begin{equation}\label{eq:bimod2}
\| \phi_n(T  J x  J) - \phi_n(T) J x J \|_{\infty, 1}, \| \phi_n(  J x  JT) - JxJ \phi_n(T) \|_{\infty, 1} \leq \| T \| \| x - \alpha_n(x) \|_2.
\end{equation}

We view each $\phi_n$ as a u.c.p.\ map into the von Neumann algebra $(\B(L^2M)^{\sharp}_J)^*$ and let $\phi: \B(L^2 \tilde M) \to (\B(L^2M)^{\sharp}_J)^*$ denote a point ultraweak limit point of $\{ \phi_n \}_n$.  From (\ref{eq:bimod1}) and (\ref{eq:bimod2}) we have that $\phi$ is bimodular with respect to both $M$ and $JMJ$. Moreover, since $e_M \circ \alpha_n \in \K_\X^L(M, L^2M)$ it follows that $\phi_n( e_M )  \in \K_\X^{\infty, 1}(M)$ for each $n \geq 1$ and hence $\phi(e_M) \in (\K_\X(M)_J^\sharp)^*$.

If $\Gamma \actson M$ were not properly proximal relative to $\X$, then by Lemma~\ref{lem:approxcentral2}, there exists a $\Gamma$-invariant $M$-central state $\varphi$ on $\widetilde {\bS}_\X(M)$ with $\varphi_{\mid M}=\tau$ and such that $\varphi$ vanishes on $(\K_\X(M)_J^\sharp)^*$. Since $\phi$ is $JMJ$-bimodular it follows that $\phi: JMJ' \cap \B(L^2 \tilde M) \to \widetilde {\bS}_\X(M)$. Since $\phi$ is $M$-bimodular and $\Gamma$-equivariant it follows that $\varphi \circ \phi$ is $M$-central and $\Gamma$-invariant, and we also see that $\varphi \circ \phi_{|M} = \tau$. And since $\varphi$ vanishes on $(\K_\X(M)_J^\sharp)^*$ and $\phi(e_M) \in (\K_\X(M)_J^\sharp)^*$ it then follows that the state $\varphi \circ \phi \circ {\rm Ad}(e_M^\perp)$ verifies the conclusion of the proposition. 
\end{proof}

As an application of the previous proposition, we give here examples of properly proximal actions stemming from Gaussian processes associated to an orthogonal representations. See, e.g., \cite{KeLi16} for some details on Gaussian actions. Let $\mathcal H$ be a real Hilbert space, the Gaussian process gives a tracial abelian von Neumann algebra $A_\cH$, together with an isometry $S: \mathcal H \to L^2_{\mathbb R}(A_\cH)$ so that orthogonal vectors are sent to independent Gaussian random variables, and so that the spectral projections of vectors in the range of $S$ generate $A_\cH$ as a von Neumann algebra. 

In this case, the complexification of the isometry $S$ extends to a unitary operator from the symmetric Fock space $\mathfrak S(\mathcal H) = \mathbb C\Omega\oplus \bigoplus_{n=1}^\infty (\cH\otimes\mathbb C)^{\odot n}$ into $L^2(A_\cH)$.  If $\mathcal H = \mathcal H_1 \oplus \mathcal H_2$, then conjugation by the unitary implementing the canonical isomorphism $\mathfrak  S (\mathcal H_1 \oplus \mathcal H_2) \cong \mathfrak S(\mathcal H_1) \ovt \mathfrak S(\mathcal H_2)$ implements a canonical isomorphism $A_{\mathcal H_1 \oplus \mathcal H_2} \cong A_{\mathcal H_1} \ovt A_{\mathcal H_2}$.

If $V: \mathcal K \to \mathcal H$ is an isometry, then we obtain an isometry $V^{\mathfrak S} : \mathfrak S(\mathcal K) \to \mathfrak S(\mathcal H)$ on the level of the symmetric Fock spaces, and conjugation by this isometry gives an embedding of von Neumann algebras ${\rm Ad}(V^{\mathfrak S}): A_{\mathcal K} \to A_{\mathcal H}$. If $V$ were a co-isometry the conjugation by $V^{\mathfrak S}$ implements instead a conditional expectation from $A_{\mathcal H}$ to $A_{\mathcal K}$. In particular, if $U \in \mathcal O(\mathcal H)$ is an orthogonal operator, then we obtain a trace-preserving $*$-isomorphism $\alpha_U = {\rm Ad}( U^{\mathfrak S} ) \in {\rm Aut}(A_{\mathcal H})$. If $\pi: \Gamma \to \mathcal O(\mathcal H)$ is an orthogonal representation, then the Gaussian action associated to $\pi$ is given by $\Gamma \ni t \mapsto \alpha_{\pi(t)} \in {\rm Aut}(A_{\mathcal H})$. The corresponding Koopman representation is then canonically isomorphic to the symmetric Fock space representation $\pi^{\mathfrak S}: \Gamma \to \mathcal U( \mathfrak S(\mathcal H) )$.

Every contraction $A: \mathcal K \to \mathcal H$ can be written as the composition of an isometry $V: \mathcal K \to \mathcal L$ followed by a co-isometry $W: \mathcal L \to \mathcal H$ so that conjugation by $A^{\mathfrak S}: \mathfrak S(\mathcal K) \to \mathfrak S(\mathcal H)$ implements a u.c.p.\ map $\phi_A: \mathbb B(L^2 A_{\mathcal H}) \to \mathbb B(L^2  A_{\mathcal K})$, which maps $A_{\mathcal H}$ into $A_{\mathcal K}$. Note that the association $A \mapsto \phi_A$ is continuous as a map from the space of contractions endowed with the strong operator topology into the space of u.c.p.\ maps endowed with the topology of point-ultraweak convergence.

If $A$ is self-adjoint this can be realized explicitly by considering the orthogonal matrix on $\mathcal H \oplus \mathcal H$ given by $\tilde A = \begin{pmatrix}
A & -\sqrt{1 - A^2} \\
\sqrt{1 - A^2} & A 
\end{pmatrix}$, so that $A = V^* \tilde A V$, where $V\xi = \xi \oplus 0$, for $\xi \in \mathcal H$. 

We recall that an orthogonal representation $\pi: \Gamma \to \mathcal O(\mathcal H)$ is amenable in the sense of Bekka \cite{Be90} if there exists a state $\varphi$ on $\B(\mathcal H \otimes \mathbb C)$ that is invariant under the action given by conjugation by $\pi(t)$. If $\pi$ is a nonamenable representation, then neither is $\pi \otimes \rho$ for any representation $\rho$ \cite[Corollary 5.4]{Be90}. In particular, if $\pi$ is a nonamenable representation, then neither is the restriction of $\pi^{\mathfrak S}$ to $\mathfrak S(\mathcal H) \ominus \mathbb C\Omega \cong \mathcal H \ovt \mathfrak S(\mathcal H)$ as it is isomorphic to $\pi \otimes \pi^{\mathfrak S}$.

\begin{thm}\label{thm:gaussian}
	Let $\pi:\Gamma\to\mathcal O(\cH)$ be a nonamenable representation on a separable real Hilbert space, then the associated Gaussian action $\Gamma\aactson{{\alpha_\pi}} A_\cH$ is properly proximal.
\end{thm}
\begin{proof}
Note that since $\pi$ is a nonamenable representation we have that $\Gamma \actson A_\cH$ is ergodic (see, e.g., \cite{PeSi12}). Denote by $P_k$ the orthogonal projection from $L^2 A_\cH$ onto $\mathbb C\Omega\oplus \bigoplus_{n=1}^k (\cH\otimes\mathbb C)^{\odot n}$ for each $k\geq 1$ and set $P_0$ to be the projection onto $\mathbb C\Omega$. 
Consider the hereditary $C^*$-algebra $\X_F \subset \B(L^2 A_\cH)$ generated by $\{P_k\}_{k\geq 0}$ in $\B(L^2 A_\cH)$, i.e., $T \in \X_F$ if and only if $\lim_{k \to \infty} \| T - TP_k \| = \lim_{k \to \infty} \| T - P_k T \| = 0$. 

To see hat $\X_F$ is a boundary piece it suffices to show that the multiplier algebra contains any von Neumann subalgebra of the form $A_{\mathbb R \xi}$ where $\xi \in \cH \setminus \{ 0 \}$. Considering the canonical isomorphisms $\mathfrak  S (\mathcal H) \cong \mathfrak S(\mathbb R \xi) \ovt \mathfrak S(\mathcal H \ominus \mathbb R \xi)$ and $A_{\mathcal H} \cong A_{\mathbb R \xi} \ovt A_{\mathcal H \ominus \mathbb R \xi}$ we may decompose $P_k$ as $P_k = \sum_{j = 0}^k Q_j \otimes R_{k - j} \leq Q_k \otimes 1$ where $Q_j$ denotes the projection onto $\mathbb C\Omega \oplus \bigoplus_{n = 1}^j ( \mathbb R \xi \otimes \mathbb C)^{\odot n}$ and $R_{k - j}$ denotes the projection onto $( ( \mathcal H \ominus \mathbb R \xi ) \otimes \mathbb C)^{\odot k - j}$. 

Since $\{ Q_j \}_{j \geq 0}$ gives an increasing family of finite rank projections such that $\vee_j Q_j = {\rm id}_{\mathfrak S(\mathbb R \xi)}$ we see that for each $a \in A_{\mathbb R\xi}$, $j \geq 0$, and $\varepsilon > 0$ there exists $k > j$ so that $\| Q_j a Q_{k - j}^\perp \| < \varepsilon$. If $i \leq j$ we have $(1 \otimes R_{j - i}) P_k^\perp = Q_{k - j + i}^\perp \otimes R_{j - i}$ and hence 
\[
\| P_j (a \otimes 1) P_k^\perp \|
=  \max_{0 \leq i \leq j} \| Q_i a Q_{k - j + i}^\perp \| 
\leq \| Q_j a Q_{k - j}^\perp \| 
< \varepsilon.
\] 
From this fact it is then easy to see that $(a \otimes 1) \in M(\mathbb X_F)$, and so we indeed have $A_{\mathbb R \xi} \subset M(\mathbb X_F)$ for each $\xi \in \mathcal H \setminus \{ 0 \}$.

For $t \geq 0$ we consider the orthogonal matrix $\begin{pmatrix}
{\cos(\pi t/2)} & -{\sin(\pi t/2)} \\
{\sin(\pi t/2)} & {\cos(\pi t/2)} 
\end{pmatrix} \in \mathcal O(\mathcal H \oplus \mathcal H)$ and let $\alpha_t \in {\rm Aut}(A_{\cH \oplus \cH})$ be the associated automorphism.

Note that $\alpha_t$ commutes with the action of $\Gamma$ on $A_{\cH}$, and that $\lim_{t \to 0} \| x - \alpha_t(x) \|_2 = 0$ for each $x \in A_{\cH}$. We may compute $E\circ \alpha_t$ explicitly as $E \circ \alpha_t = P_0 + \sum_{n = 1}^\infty \cos^n(\pi t/2) (P_n - P_{n - 1})$, so that $E \circ \alpha_t \in \X_F$ for each $0 < t < 1$. Using the isomorphism $L^2( A_{\cH \oplus \cH} ) \cong L^2 (A_{\cH} ) \ovt L^2(A_{\cH}) \cong \mathfrak S( \cH) \ovt \mathfrak S( \cH)$ we obtain an isomorphism of $A_{\cH}$-modules $L^2 ( A_{\cH \oplus \cH} ) \ominus L^2(A_{\cH \oplus 0} ) \cong L^2(A_{\cH \oplus 0}) \ovt ( \mathfrak S( \mathcal H) \ominus \mathbb C \Omega )$. So that we have an isomorphism 
\[
A_{\cH}' \cap \B(L^2 ( A_{\cH \oplus \cH} ) \ominus L^2(A_{\cH \oplus 0} ) ) \cong A_{\cH} \ovt \B( \mathfrak S(\mathcal H) \ominus \mathbb C \Omega),
\]
where $\Gamma$ acts on the latter space by $\alpha_{\pi(t)} \otimes {\rm Ad}( \pi(t)^{\mathfrak S} )$. If we had a $\Gamma$-invariant state on $A_{\cH} \ovt \B( \mathfrak S(\mathcal H) \ominus \mathbb C \Omega)$, then restricting it to $\B( \mathfrak S(\mathcal H) \ominus \mathbb C \Omega)$ would show that the restriction of the representation $\pi^{\mathfrak S}$ to $\mathfrak S(\mathcal H) \ominus \mathbb C \Omega$ is amenable, which would then imply that $\pi$ is an amenable representation, giving a contradiction. Thus, we conclude that there is no $\Gamma$-invariant state on $A_{\cH} \ovt \B( \mathfrak S(\mathcal H) \ominus \mathbb C \Omega)$ and it then follows from Proposition~\ref{prop:deformation} that $\Gamma \actson A_{\cH}$ is properly proximal relative to $\X_F$.

We now take a $\Gamma$-almost invariant approximate identity $A_n \in \K(\mathcal H)$ with $\| A_n \| < 1$ for each $n \geq 0$, and let $\alpha_{\tilde A_n}$ be the automorphisms of $A_{\cH \oplus \cH}$ as given above. As in the proof of Proposition~\ref{prop:deformation} we let $V_n: L^2(A_{\cH}) \to L^2(A_{\cH \oplus \cH})$ be given by $V_n(\hat{x}) = \widehat {\alpha_{\tilde A_n}(x)}$, for $x \in A_{\cH}$. We also denote by $\phi_n: \B(L^2 ( A_{\cH \oplus \cH} ) ) \to \B(L^2 (A_{\cH}))$ the u.c.p.\ map $\phi_n(T) = V_n^* T V_n$, for $T \in \B(L^2 ( A_{\cH \oplus \cH} ) )$. 

We view each $\phi_n$ as a u.c.p.\ map into the von Neumann algebra $( {\X_F}^\sharp_J )^*$, which we view as a corner of $( \B(L^2(A_{\cH}))^\sharp_J)^*$. We then let $\phi: \B(L^2 ( A_{\cH \oplus \cH} ) ) \to ( {\X_F}^\sharp_J )^*$ be a point ultraweak cluster point of these maps. Since $A_n \to 1$ in the strong operator topology it follows that $\| \alpha_{\tilde A_n}(x) - x \|_2 \to 0$ for each $x \in A_{\cH \oplus \cH}$. Also, since each $A_n$ is compact, and since $\| A_n \| < 1$, it follows that $A_n^{\mathfrak S} \in \K( \mathfrak S(\mathcal H))$ and the same computation as above then shows that $E_{A_{\mathcal H \oplus 0}} \circ \alpha_{\tilde A_n} \in \K(L^2M)$. The proof of Proposition~\ref{prop:deformation} then shows that $\phi$ is $A_{\cH}$-bimodular and satisfies $\phi( e_{A_{\cH \oplus 0}}) \in (\K(M)_J^\sharp)^*$. 

If we take $t \in \Gamma$, then as $\| \pi(t) A_n \pi(t^{-1}) - A_n \| \to 0$ we have that $\lim_{n \to \infty} \| ( \alpha_{\pi(t) \oplus \pi(t)} V_n  - V_n \alpha_{\pi(t)} ) P_k \| = 0$ for each $k \geq 1$, and hence for each $T \in  \B(L^2 ( A_{\cH \oplus \cH} ) )$ we have
\[
\lim_{n \to \infty} \| P_k ( \alpha_{\pi(t)} \phi_n(T) \alpha_{\pi(t^{-1})} - \phi_n(  \alpha_{\pi(t) \oplus \pi(t)} T  \alpha_{\pi(t^{-1}) \oplus \pi(t^{-1})} ) ) P_k \| = 0.
\] 
Since $\{ P_k \}_k$ gives an approximate identity for $\X_F$ it then follows by passing to a limit that $\phi$ is $\Gamma$-equivariant.

Suppose now that $\Gamma\actson A_\cH$ is not properly proximal, then by Lemma~\ref{lem:approxcentral}
there exists an $A_\cH$-central and $\Gamma$-invariant state $\varphi$ on $\tilde \bS(M)$. We may further assume it vanishes on $((\K^{\infty,1}_{A_\cH})^{\sharp })^*$. As $\Gamma\actson A_\cH$ is properly proximal relative to $\X_F$ we must then have that $\varphi$ is supported on $((\X_F)_J^\sharp)^*$, but then by considering the composition of $\varphi$ with the restriction of $\phi$ to $A_{\cH}' \cap \B( L^2(A_{\cH \oplus \cH})  \ominus L^2(A_{\cH \oplus 0}) )$ we obtain a $\Gamma$-invariant state on this space, which would give a contradiction, as noted above. 
\end{proof}

\begin{proof}[Proof of Theorem~\ref{thm:soliderg}]
From the proof of Theorems~\ref{prop:deformation} and \ref{thm:gaussian} we show, in fact, that if $\pi: \Gamma \to \mathcal O(\mathcal H)$ is an orthogonal representation, then there are $A_{\mathcal H}$-bimodular $\Gamma$-equivariant maps $\psi: A_{\mathcal H} \ovt \B( \mathfrak S(\mathcal H) \ominus \mathbb C \Omega) \to \tilde \bS_{\X_F}( A_{\mathcal H})$ and $\phi: A_{\mathcal H} \ovt \B( \mathfrak S(\mathcal H) \ominus \mathbb C \Omega) \to  ((\X_F)_J^\sharp)^* \cap \tilde \bS( A_{\mathcal H})$. If $\mathcal R$ is a subequivalence relation of the orbit equivalence relation associated to the action $\Gamma \actson A_{\mathcal H}$, then just as above, if the action $\mathcal N_{L\mathcal R}(A_{\mathcal H}) \actson A_{\mathcal H}$ is not properly proximal, then we would obtain a $\mathcal N_{L\mathcal R}(A_{\mathcal H})$-invariant state $\varphi$ on $A_{\mathcal H} \ovt \B(  \mathfrak S(\mathcal H) \ominus \mathbb C \Omega)$ such that $\varphi_{|A_{\mathcal H}}$ is normal with support $p \in A_{\mathcal H}^\Gamma$. Under the isomorphism $\mathcal N_{L\mathcal R}(A_{\mathcal H})/ \mathcal U(A_{\mathcal H}) \cong [\mathcal R]$, where $[\mathcal R]$ denotes the full group of the equivalence relation we then see that $\varphi$ is $[\mathcal R]$-invariant for the natural action of $[\mathcal R]$ on $A_{\mathcal H} \ovt \B(  \mathfrak S(\mathcal H) \ominus \mathbb C \Omega)$.

By a standard argument using Day's trick (e.g., as in \cite{Be90}) this would then give a net of vectors $\xi_i \in L^2(A_{\mathcal H}) \ovt HS( \mathfrak S(\mathcal H) \ominus \mathbb C \Omega)$ such that $\langle a \xi_i, \xi_i \rangle = \frac{1}{\tau(p)}\tau(ap)$ for $a \in A_{\mathcal H}$ and such that $\{ \xi_i \}_i$ is asymptotically $[\mathcal R]$-invariant. If $\pi \prec \lambda$, then this, in turn, gives an asymptotically $[\mathcal R]$-invariant net of vectors in $L^2(A_{\mathcal H}) \ovt \ell^2 \Gamma$, which would show that $\mathcal R_{p}$ is an amenable equivalence relation \cite[Theorem 6.1.4]{AnDeRe00}. A simple maximality argument then gives the result. 

More generally, if some tensor power $\pi^{\otimes k}$ satisfies $\pi^{\otimes k} \prec \lambda$, then by considering self-tensor powers of the net $\{ \xi_i \}_i$ it would again follow that $\mathcal R_{p}$ is an amenable equivalence relation. 
\end{proof}

\subsection{Proper proximality via closable derivations}

Let $M$ be a tracial von Neumann algebra and $\mathcal H$ an $M$-$M$ correspondence that has a real structure, i.e., there exists an antilinear isometry $\mathcal J: \mathcal H \to \mathcal H$ such that $\mathcal J( x \xi y ) = y^* \mathcal J(\xi) x^*$ for all $x, y \in M$, $\xi \in \mathcal H$. A closable real derivation is an unbounded closable linear map $\delta: L^2 M \to \mathcal H$, such that the domain $D(\delta)$ is an ultraweakly dense unital $*$-subalgebra of $M \subset L^2M$, and such that $\delta$ preserves the real structure ($\delta(x^*) = \mathcal J( \delta(x))$ for $x \in D(\delta)$) and satisfies Leibniz's formula
\[
\delta(xy) = x \delta(y) + \delta(x)y \ \ \ \ \ x, y \in D(\delta).
\]

A result of Davies and Lindsay in \cite{DaLi92} shows that $D(\overline{\delta}) \cap M$ is then again a $*$-subalgebra and $\overline \delta_{| D(\overline \delta) \cap M}$ again gives a closable real derivation. We recycle the following notation from \cite{Pe09, OzPo10II} 
\[
\Delta = \delta^* \overline{ \delta},  \ \ \ \ \ \zeta_\alpha = \left( \frac{\alpha}{\alpha + \Delta} \right)^{1/2},  \ \ \ \ \ \tilde \delta_\alpha = \frac{1}{\sqrt{\alpha}} \overline \delta \circ \zeta_\alpha, \ \ \ \ \ \tilde \Delta_\alpha = \frac{1}{\sqrt{\alpha}} \Delta^{1/2} \circ \zeta_\alpha, \ \ \ \ \ \theta_\alpha = 1 - \tilde \Delta_\alpha.
\] 
We note that from \cite{Sa90, Sa99} (cf. \cite{Pe09}) we have that $\zeta_\alpha$ and $\theta_\alpha$ are $\tau$-symmetric u.c.p.\ maps.

For $x, y \in D(\overline \delta) \cap M$ we set
\[
\Gamma(x^*, y) = \Delta^{1/2}(x^*) y + x^* \Delta^{1/2}(y) - \Delta^{1/2}(x^* y),
\]
and recall from \cite[p. 850]{OzPo10II} that we have
\begin{equation}\label{eq:carre}
\| \Gamma(x^*, y) \|_2 
\leq 4 \| x \|^{1/2} \| \delta(x) \|^{1/2} \| y \|^{1/2} \| \delta(y) \|^{1/2}. 
\end{equation}

We also recall \cite[Lemma 4.1]{OzPo10II} (cf. \cite[Lemma 3.3]{Pe09}) that for $x, a \in M$ we have the following approximate bimodularity property
\begin{equation}\label{eq:approxbimod1}
\| \zeta_\alpha(x) \tilde \delta_\alpha(a) - \tilde \delta_\alpha(xa) \|
\leq 10 \| a \| \| x \|^{1/2} \| \tilde \delta_\alpha(x) \|^{1/2},
\end{equation}
and 
\begin{equation}\label{eq:approxbimod2}
\| \tilde \delta_\alpha(a) \zeta_\alpha(x) - \tilde \delta_\alpha(ax) \|
\leq 10 \| a \| \| x \|^{1/2} \| \tilde \delta_\alpha(x) \|^{1/2}.
\end{equation}

If $x \in M$ is contained in the domain $D(\overline \delta) \cap M$, then we have the following estimate, allowing us to replace the term $\zeta_\alpha(x)$ from (\ref{eq:approxbimod1}) and (\ref{eq:approxbimod2}) with $x$.

\begin{lem}\label{lem:almostbimod}
For $x \in M \cap D(\overline \delta)$, $a \in M$, and $\alpha > 1$ we have
\[
\| x \tilde \delta_\alpha(a) - \tilde \delta_\alpha(xa) \| \leq \alpha^{-1/4} (2 \| \delta(x) \|^{1/2} + 6\| x \|^{1/2} ) \| \delta(x) \|^{1/2} \| a \|. 
\]
and
\[
\| \tilde \delta_\alpha(a)x - \tilde \delta_\alpha(ax) \| \leq \alpha^{-1/4} (2 \| \delta(x) \|^{1/2} + 6\| x \|^{1/2} ) \| \delta(x) \|^{1/2} \| a \|.  
\]
\end{lem}
\begin{proof}
We follow the strategy from \cite{Pe09} and \cite{OzPo10II}. We have
\[
x \tilde \delta_\alpha(a) = \alpha^{-1/2} \delta( x \zeta_\alpha(a) ) - \alpha^{-1/2} \delta(x) \zeta_\alpha(a) =: A_1 - A_2.
\]
Note that $\| A_2 \| \leq \alpha^{-1/2} \| \delta(x) \| \| a \|$. Let $\delta = V \Delta^{1/2}$ be the polar decomposition. Then, one has
\begin{align*}
V^* A_1 &= x \tilde \Delta_\alpha(a)  &&+ \alpha^{-1/2} \Delta^{1/2}(x) \zeta_\alpha(a) &&- \alpha^{1/2} \Gamma(x, \zeta_\alpha(a) ) \\
&=: B_1 &&+ B_2 &&- B_3.
\end{align*}
We have $\| B_2 \|_2 \leq \alpha^{-1/2} \| \delta(x) \| \| a \|$, and by (\ref{eq:carre}) we have $\| B_3 \|_2 \leq 4 \alpha^{-1/4} \| a \| \| x \|^{1/2} \| \delta(x) \|^{1/2}$. 
We have $B_1 - \tilde \Delta_\alpha(xa) = x \theta_\alpha(a) - \theta_\alpha(xa)$, and since $\theta_\alpha$ is a $\tau$-symmetric u.c.p.\ map, we may use the Stinespring representation $\theta_\alpha(y)  = W^* \pi(y) W$ to estimate 
\begin{align*}
\| x \theta_\alpha(a) - \theta_\alpha(xa) \|_2^2
&= \| \theta_\alpha(a^* x^*) - \theta_\alpha(a^*) x^* \|_2^2 \\
&= \| W^* \pi(a^*) ( \pi(x^*) W -  W \pi(x^*) ) \hat{1} \|_2^2 \\ 
&\leq \| a \|^2 \| ( \pi(x^*) W -  W \pi(x^*) ) \hat{1} \|_2^2 \\
&= 2 \| a \|^2 \tau(( x - \theta_\alpha(x) ) x^*) \\ 
&\leq 2 \| a \|^2 \| x \| \| x - \theta_\alpha(x) \|_2 \\
&= 2 \| a \|^2 \| x \| \| \tilde \delta_\alpha(x) \|
\leq 2 \alpha^{-1/2} \| a \|^2 \| x \| \| \delta(x) \|.
\end{align*}
Combining these estimates gives the first part of the lemma, and the second part follows since $\delta$ is a real derivation. 
\end{proof}

We now show that the argument in Proposition~\ref{prop:deformation} above can be adapted to deformations associated to closable derivations.

\begin{prop}\label{prop:derivationprox}
Let $M$ be a finite factor, $\mathcal H$ an $M$-$M$ correspondence with a real structure, and $\delta: L^2M \to \mathcal H$ a closable real derivation. Let $A$ denote the $C^*$-algebra generated by $D(\overline{\delta}) \cap M$ and suppose that $(M^{\rm op})' \cap \mathbb B(\mathcal H)$ has no $A$-central state $\psi$ such that $\psi_{| A} = \tau$. If $\X \subset \B(L^2M)$ is a boundary piece and $\zeta_\alpha \in \K_\X^L(M, L^2M)$ for each $\alpha > 0$, then $M$ is properly proximal relative to $\X$.
\end{prop}
\begin{proof}
For $\alpha > 0$ let $V_\alpha: L^2M \to L^2M \oplus \mathcal H$ denote the map given by $V_\alpha(\hat{x}) = \zeta_{\alpha}(x) \oplus \tilde \delta_\alpha (x)$. For $x \in M$ we compute 
\[
\| V_\alpha(\hat{x}) \|^2 = \tau( \rho_\alpha(x) x^* ) + \frac{1}{\alpha} \tau(\Delta \circ \rho_\alpha(x) x^*) = \| x \|_2^2,
\] 
where $\rho_\alpha = \frac{\alpha}{\alpha + \Delta}$. Hence $V_\alpha$ is an isometry. By Lemmas 3.5 in \cite{OzPo10II} and from (\ref{eq:approxbimod1}) above we have 
\begin{align}
\| V_\alpha x - \zeta_\alpha(x) V_\alpha \|_{\infty, 2}^2
&= \sup_{\| a \| \leq 1} \| \zeta_\alpha(xa) - \zeta_\alpha(x)\zeta_\alpha(a) \|_2^2 + \| \tilde \delta_\alpha(xa) - \zeta_\alpha(x) \tilde \delta_\alpha(a) \|^2 \nonumber \\ 
&\leq 4 \| x \| \| x - \zeta_\alpha(x) \|_2 + 100 \| a \| \| \tilde \delta_\alpha(x) \| \label{eq:approxbim1} \\
&\leq 104 \| x \| \| \tilde \delta_\alpha(x) \|, \nonumber
\end{align}
where the last inequality follows from the fact that $1 - \sqrt{t} \leq \sqrt{1 - t}$ for $0 \leq t \leq 1$ as in \cite[p. 850]{OzPo10II}. It similarly follows that 
\begin{equation}\label{eq:approxbim2}
\| V_\alpha J x J - (J \oplus \mathcal J) \zeta_\alpha(x) (J \oplus \mathcal J) V_\alpha \|_{\infty, 2}^2 \leq 104 \| x \| \| \tilde \delta_\alpha(x) \|.
\end{equation}

Moreover, if $x \in D(\overline{\delta}) \cap M$, then we may use instead Lemma~\ref{lem:almostbimod} to conclude that for $\alpha \geq 1$ we have
\begin{equation}\label{eq:approxbim3}
\| V_\alpha x - x V_\alpha \|_{\infty, 2}
\leq \alpha^{-1/4} ( 2 \| \delta(x) \|^{1/2} + 8 \| x \|^{1/2} ) \| \delta(x) \|^{1/2}.
\end{equation}

As in the proof of Proposition~\ref{prop:deformation} we define u.c.p.\ maps $\phi_\alpha: \B(L^2 M \oplus \mathcal H) \to \B(L^2M)$ by $\phi_\alpha(T) = V_\alpha^* T  V_\alpha$, and we let $\phi: \B(L^2 M \oplus \mathcal H) \to (\B(L^2M)^{\sharp}_J)^*$ be a point ultraweak limit point of $\{ \phi_\alpha \}_{\alpha}$. 
If we denote by $e$ the orthogonal projection from $L^2M \oplus \mathcal H$ down to $L^2M$, then by hypothesis we have $e V_\alpha = \zeta_\alpha \in \K_{\X}^L(M, L^2M)$ and hence $\phi(e) \in ((\K^{\infty,1}_\X(M))_J^{\sharp})^*$.

If $M$ is not properly proximal relative to $\X$, then by Lemma~\ref{lem:approxcentral2}, there exists an $M$-central state $\varphi$ on $\widetilde {\bS}_\X(M)$ with $\varphi_{\mid M} = \tau$ is normal and such that $\varphi$ vanishes on $((\K^{\infty,1}_\X(M))_J^{\sharp})^*$. 

For $T \in (M^{\rm op})' \cap \B(L^2 M \oplus \mathcal H)$ and $x \in M$ we have that $[\phi(T), JxJ]$ is an ultraweak limit point of $[\phi_\alpha(T), JxJ]$ and from (\ref{eq:approxbim2}) we see that this is an ultraweak limit point of $\phi_\alpha( [T, (J \oplus \mathcal J) \zeta_\alpha(x) (J \oplus \mathcal J)] ) = 0$. Therefore we see that $\phi$ maps into $\widetilde {\bS}(M) \subset \widetilde {\bS}_\X(M)$. Also, by (\ref{eq:approxbim3}) we see that $\phi$ is bimodular for the $C^*$-closure $A$ of $D(\overline \delta) \cap M$. Thus, $\varphi \circ \phi \circ {\rm Ad }(e^\perp)$ gives an $A$-central state on $(M^{\rm op})' \cap \B(L^2 M \oplus \mathcal H)$.
\end{proof}

We remark that since $\zeta_\alpha$ is positive as an operator on $L^2M$ we have that $\zeta_\alpha \in \K_{\X}^L(M, L^2M)$ if and only if $\zeta_\alpha \in \K_\X(M)$, and this is if and only if $\frac{\alpha}{\alpha + \Delta} \in \K_\X(M)$. 

The details of the following example for which the previous proposition applies can be found in \cite{Pe09I} (cf. also \cite{ChPe13}). Suppose $\Gamma$ is a group, $\pi: \Gamma \to \mathcal O(\mathcal K)$ is an orthogonal representation, and $c: \Gamma \to \mathcal K$ is a $1$-cocycle. Then $\mathcal K \ovt \ell^2 \Gamma$ is naturally an $L\Gamma$-correspondence and we obtain a closable real derivation $\delta: \mathbb C\Gamma \to \mathcal K \ovt \ell^2 \Gamma$ by setting $\delta(u_t) = c(t) \otimes \delta_t$ for $t \in \Gamma$. It is then not hard to see that $(L\Gamma^{\rm op})' \cap \B( \mathcal K \ovt \ell^2 \Gamma) = \B(\mathcal K) \ovt L\Gamma$ has a $C^*_\lambda \Gamma$-central state if and only if $\pi$ is an amenable representation in the sense of Bekka. This then gives another way to show proper proximality for group von Neumann algebras associated with groups that have property (HH) from \cite{OzPo10II}.

The following lemma will allow us to apply Proposition~\ref{prop:derivationprox} when the $M$-$M$ correspondence $\mathcal H$ is weakly contained in the coarse correspondence.

\begin{lem}\label{lem:coarse}
Let $M$ be a finite factor, $\mathcal H$ an $M$-$M$ correspondence with a real structure that is weakly contained in the coarse correspondence and suppose $\delta: L^2M \to \mathcal H$ is a closable real derivation with $A = D( \overline{\delta}) \cap M$. If $(M^{\rm op})' \cap \B(\mathcal H)$ has an $A$-central state $\psi$ such that $\psi_{|A} = \tau$, then $M$ is amenable.
\end{lem}
\begin{proof}
As $M$-$M$ correspondences we have $L^2( (M^{\rm op})' \cap \B(\mathcal H) ) \cong \mathcal H \oovt{M} \overline{\mathcal H}$ \cite[Proposition 3.1]{Sa83} and so if there existed an $A$-central state on $(M^{\rm op})' \cap \B(\mathcal H)$ that restricts to the trace on $A$, then by Day's trick \cite{Co76II} there would exist a net of unit vectors $\{ \xi_i \}_i \subset \mathcal H \oovt{M} \overline{\mathcal H}$ such that $\| [a, \xi_i ] \| \to 0$ and $\langle a \xi_i, \xi_i \rangle \to \tau(a)$ for each $a \in A$. 

If $a_1, \ldots, a_n \in A$, and $b_1, \ldots, b_n \in A$, we then have
\begin{align*}
 \| \sum_{k = 1}^n a_k  b_k^* \|_2 
 &= \lim_{i \to \infty} \| ( \sum_{k = 1}^n a_k   b_k^* ) \xi_i \| \\
 &= \lim_{i \to \infty} \| \sum_{k = 1}^n a_k  \xi_i b_k^* \| \\
 &\leq \| \sum_{k =1}^n a_k \otimes Jb_kJ \|_{\B(L^2M \ovt L^2M)}  \\
 &\to \| \sum_{k =1}^n a_k \otimes Jb_kJ \|_{\B(L^2M \ovt L^2M)} 
\end{align*}
where the inequality follows since $\mathcal H \oovt{M} \overline{\mathcal H}$ is weakly contained in the coarse correspondence.

If we now take unitaries $u_1, \ldots, u_n \in M$, then for $\alpha > 0$ we have $\zeta_{\alpha}(u_1 ), \ldots, \zeta_\alpha(u_n  ) \in D(\overline \delta) \cap M$ and hence
\begin{align*}
\| \sum_{k = 1}^n \zeta_\alpha(u_k )  \zeta_\alpha(  u_k^* ) \|_2 
&\leq \| \sum_{k = 1}^n \zeta_\alpha(u_k ) \otimes J \zeta_\alpha(u_k  ) J \|_{\B(L^2M \ovt L^2M)} \\
&\leq \| \sum_{k =1}^n u_k \otimes Ju_k J \|_{\B(L^2M \ovt L^2M)}.
\end{align*}
Since $\zeta_\alpha(x) \to x$ strongly as $\alpha \to \infty$, we then conclude that
\[
n = \lim_{\alpha \to \infty} \| \sum_{k = 1}^n \zeta_\alpha(u_k  )   \zeta_\alpha(  u_k^*)  \|_2 
\leq \| \sum_{k =1}^n u_k   \otimes J  u_k J \|_{\B(L^2M \ovt L^2M)},
\]
and hence $M$ is amenable by \cite{Co76}.
\end{proof}

We refer to \cite{Ca21}, and the references therein, for background and our notation for free orthogonal quantum groups. 

\begin{cor}
If $N \geq 3$, then the von Neumann algebra $L_\infty(O_N^+)$ associated to the free orthogonal quantum group $O_N^+$ is properly proximal. 
\end{cor}
\begin{proof}
By \cite[Section 7.1]{Ca21} (cf. \cite{FiVe15}) there exists a $L_\infty(O_N^+)$-$L_\infty(O_N^+)$ correspondence $\mathcal H$ with a real structure that is weakly contained in the coarse corrspondence, and a closable real derivation $\delta: L_\infty(O_N^+) \to \mathcal H$ so that $\frac{\alpha}{\alpha + \delta^* \overline \delta}$ is compact as an $L^2$-operator for each $\alpha > 0$. The result then follows from Lemma~\ref{lem:coarse} and Proposition~\ref{prop:derivationprox}.
\end{proof}

\providecommand{\bysame}{\leavevmode\hbox to3em{\hrulefill}\thinspace}
\providecommand{\MR}{\relax\ifhmode\unskip\space\fi MR }
\providecommand{\MRhref}[2]{%
  \href{http://www.ams.org/mathscinet-getitem?mr=#1}{#2}
}
\providecommand{\href}[2]{#2}

\end{document}